\documentclass[12pt]{amsart} 
\usepackage{ifthen,color,xypic,amssymb}

\newtheorem{thm}{Theorem}[section] 
\newtheorem{corol}[thm]{Corollary}
\newtheorem{lemma}[thm]{Lemma} 
\newtheorem{prop}[thm]{Proposition}
\newtheorem*{conj}{Conjecture}
\theoremstyle{definition}
\newtheorem{defin}[thm]{Definition}
\newtheorem{note}[thm]{Note}
\theoremstyle{remark}
\newtheorem{remark}[thm]{Remark}
\numberwithin{equation}{section}

\newcommand\what[1]{\widehat{#1}}
\newcommand\op[1]{\operatorname{#1}}
\DeclareMathOperator{\Ker}{Ker}
\DeclareMathOperator{\Img}{Im}
\DeclareMathOperator{\rk}{rk}
\DeclareMathOperator{\Tr}{Tr}
\DeclareMathOperator{\tr}{tr}

\newcommand\sEnd{{\mathcal{E}xt}}
\DeclareMathOperator{\Ext}{Ext}
\DeclareMathOperator{\Hom}{Hom}
\DeclareMathOperator{\td}{td}
\DeclareMathOperator{\ch}{ch}

\DeclareMathOperator{\Supp}{Supp}

\newcommand\iso{\kern.35em{\raise3pt\hbox{$\sim$}\kern-1.1em\to}\kern.3em}

\newcommand\fp{\times_{ B}}
\newcommand\U{{\mathcal U}}
\newcommand\rest[2]{{#1}_{\vert #2}}

\newcommand\Fc{{\mathcal F}}
\newcommand\Oc{{\mathcal O}}
\newcommand\Pc{{\mathcal P}}
\newcommand\M{{\mathcal M}}
\newcommand\Qc{{\mathcal Q}}

\newcommand\Z{{\mathbb Z}}
\newcommand\bbC{{\mathbb C}}
\newcommand\bbQ{{\mathbb Q}}

\newcommand\V{{\mathcal V}}

\newcommand\Rca{{\mathcal R}}

\newcommand\Nc{{\mathcal N}}
\newcommand\Lcl{{\mathcal L}}
\newcommand\Ec{{\mathcal E}}
\newcommand\Ac{{\mathcal A}}
\newcommand\Ps{{\mathbb P}}
\newcommand\X{{\widehat X}}
\newcommand\G{{\mathcal G}}

\newcommand\CA{{\mathcal C}({\mathcal A})}
\newcommand\Hc{{\mathcal H}}
\newcommand\bS{{\mathbf\Phi}}

\newcommand\Ic{{\mathcal I}}
\newcommand\Jc{{\mathcal J}}

\newcommand\bC{{\mathbf  C}}

\newcommand\bM{{\mathbf  M}}
\newcommand\bJ{{\mathbf  J}}

\newcommand\mg{{\mathfrak m}}
\newcommand\FM{Fourier-Mukai transform}
\newcommand\gif{geometric integral functor}
\newcommand\GIF{Geometric integral functor}
\newcommand\FMF{Fourier-Mukai functor}

\newcommand{\marginnote}[1]{\ifthenelse{\isodd{\thepage}}{\normalmarginpar}
{\reversemarginpar}\marginpar{\fbox{\parbox{18mm}{\sloppy\footnotesize #1}}}}

\begin{document}
\title{FOURIER MUKAI TRANSFORMS AND APPLICATIONS TO STRING THEORY}
\author{Bj\"orn Andreas}
\email{andreas@mathematik.hu-berlin.de}
\address{Institut f\"ur Mathematik, Humboldt Universit\"at zu
Berlin, D-10115 Berlin, Germany}
\author{Daniel Hern\'andez Ruip\'erez}
\email{ruiperez@usal.es}
\address{Departamento de Matem\'aticas and Instituto Universitario de F\'{\i}sica Fundamental y Matem\'aticas (IUFFYM), Universidad de Salamanca, Plaza
de la Merced 1-4, 37008 Salamanca, Spain}
\date{\today} 
\thanks {B.A. is supported by DFG Schwerpunktprogramm (1096) ``String Theory im Kontext
von Teilchenphysik, Quantenfeldtheorie, Quantengravitation, Kosmologie und Mathematik''. 
D.H.R. is supported by DGI research project BFM2003-00097 ``Transformadas Geom\'etricas Integrales y 
Aplicaciones'' and by JCYL research project  SA114/04 ``Aplicaciones de los functores integrales a la Geometr\'{\i}a  y a la F\'{\i}sica''.
}
\subjclass{14J60, 14J32, 18E30, 81T30, 83E30}
\keywords{stable sheaves and vector bundles, semistable sheaves and
vector bundles, moduli, elliptic fibrations, elliptic surfaces, \gif s,
\FM s, compactified Jacobians, spectral covers, D-branes, T-duality, heterotic strings}
\begin{abstract} 
We give an introductory review of Fourier-Mukai transforms and their application to various aspects of moduli problems,
string theory and mirror symmetry. We develop the necessary mathematical 
background for Fourier-Mukai transforms such as aspects of derived categories and integral 
functors as well as their relative version which becomes important for making precise the notion of  fiberwise T-duality on elliptic Calabi-Yau threefolds. We discuss various applications of the 
Fourier-Mukai transform to D-branes on Calabi-Yau manifolds as well as homological mirror symmetry and the construction of  vector bundles for heterotic string theory. 

\end{abstract}
\maketitle 

{\small
\tableofcontents
}

\section{Introduction}

The interplay between geometry and physics has a long story. Traditionally differential geometry played a fundamental role in many physical aspects such as general relativity or gauge theory. Today various methods and objects of algebraic geometry are required for formulating and understanding string theory.
In particular the physical notion of ``duality'' has lead in mathematics to unexpected connections between the geometry of different spaces. Mirror symmetry was an example of this but many more remain to be explored. On the physics side one hopes to obtain a better understanding of nonperturbative aspects of the way string theory describes the real world. So both mathematics and physics appear to benefit greatly from the duality correspondences. These sparked various new developments on both sides, in
particular, a new geometrical understanding of mirror symmetry has begun recently to evolve. Some aspects of this evolution, like Kontsevich proposal of homological mirror symmetry but also the construction of holomorphic bundles from spectral data require the use of \FM s and inspired the present review paper. 

In this paper we will review some aspects of derived categories, \FM s and their relative version, that is, their formulation for families rather than for single varieties.The main advantage of the relative setting is that base-change properties (or parameter dependencies) are better encoded into the problem; the drawback is of course the increasing abstraction and technical machinery we need. In most of this review we will be concerned with the relative Fourier-Mukai transform for elliptic fibrations. However, to begin let us briefly recall a few aspects of Fourier-Mukai transforms to put the subject in context and describe how the relative setting enters into string theory.

The Fourier-Mukai transform was introduced in the study of abelian varieties
by Mukai and can be thought of as a nontrivial algebro-geometric analogue of the Fourier transform. Since its original introduction, the Fourier-Mukai transform turned out to be a useful tool
for studying various aspects of sheaves on varieties and their moduli spaces, and as a natural consequence, to learn about the varieties themselves.
In recent years it was found that the Fourier-Mukai transform also enters into string theory. The most prominent example is given by Kontsevich's homological mirror symmetry conjecture \cite{Kon95}. The conjecture predicts (for mirror dual pairs of Calabi-Yau manifolds) an equivalence between the bounded derived category of coherent sheaves and the Fukaya category. 
The conjecture implies a correspondence between self-equivalences of the derived category and certain symplectic self-equivalences of the mirror manifold.

Besides their importance for geometrical aspects of mirror symmetry, the Fourier-Mukai transform 
has been shown to be also important for heterotic string compactifications on ellipitic fibrations.
The motivation for this came from the conjectured correspondence between the heterotic string
and the so called F-theory which both rely on elliptically fibered Calabi-Yau manifolds. To give 
evidence for this correspondence, an explicit description of stable holomorphic vector bundles 
was required and inspired the seminal work of Friedman, Morgan and Witten \cite{FM, FMW97, FMW99}.  They showed how to construct vector bundles in terms of two geometrical objetcs: a hypersurface in the Calabi-Yau manifold together with a line bundle on it (called the spectral data). Various aspects and refinements of this construction have been studied in subsequent work \cite{Don98,AsDo,DLO,DoOv,Bri98,And,Cu98,DiacIon,HMP02,AnHR03,AnHR04}. Moreover, a physical way to understand this bundle construction can be given if one views holomorphic vector bundles as D-branes and uses the fact that D-branes are mapped under T-duality to new D-branes (of different dimensions) which can then be made mathematically precise in terms of a relative \FM.

More generally, D-branes can be interpreted as objects of the derived category, one then expects the Fourier-Mukai transform (or its relative version) to act on the spectrum of D-branes. This suggests that the Fourier-Mukai transform is actually a symmetry of string theory. Furthermore, the study of D-branes on Calabi-Yau manifolds inspired numerous mathematical questions, for instance, the search for new Fourier-Mukai partners \cite{Muk81,Or97,BM2,Kaw02b,HLOY03,Ue04}, the reconstruction of the underlying variety from the dervied category \cite{BO01}, the study of the self-equivalence/monodromy correspondence in the context of mirror symmetry and the search for a new categorical stability notion \cite{Bri02pp,Bri03pp} which has been motivated by Douglas $\Pi$-stability of D-branes \cite{Do01a,Do01b,Do02}.

\medskip
The paper is organized as follows: In section \ref{s:dercat} we review the definition and basic facts about the derived category and derived functors. We have tried to avoid technicalities as much as possible. In section \ref{s:gifs} we define Fourier-Mukai functors and some examples of them, together with their composition and IT and WIT conditions. We also deal with Fourier-Mukai functors for families, and we state the base-change properties of those functors. Since the main example of a relative Fourier-Mukai functor is the one defined for elliptic fibrations, we have devoted section \ref{s:elliptic} to those fibrations, their Weierstrass models and particular aspects of the  Fourier-Mukai transforms for them. We also define spectral covers and so prepare further applications in Section \ref{s:applmirror}. The computation of the topological invariants of the Fourier-Mukai transform is given in section \ref{s:topinv}. Section \ref{s:applmod} is devoted to the application of the Fourier-Mukai transform to certain moduli problems, like the determination of the moduli of relatively stable sheaves on an elliptic fibrations or the isomorphism of the moduli of absolutely stable sheaves on an elliptic surface with a (generically) integral system over a Hilbert scheme. In Section \ref{s:applmirror} we discuss some applications of the Fourier-Mukai transform to string theory. After giving a brief introduction to D-branes, we discuss the action of the Fourier-Mukai transform on the spectrum of D-branes on elliptically fibered Calabi-Yau threefolds
(in physical terms this reflects the adiabatic character of T-duality). Then we outline a procedure 
(which relies on the comparison of central charges associated to A-type, respectively, B-type 
D-brane configurations) that allows to make explicit Kontsevich's proposed self-equivalence/monodromy correspondence. We refer to monodromies in the
moduli space of the complexified K\"ahler form which have to be identified by mirror symmetry with the complex structure moduli space of the mirror manifold. If the complex structure varies while keeping fixed the K\"ahler structure, the isotopy classes of Lagrangian submanifolds vary as well and then loops in the complex structure moduli space produce monodromies on the classes of Lagrangian submanifolds. 
The last application is devoted to heterotic string theory on elliptically fibered Calabi-Yau threefolds. It is shown how the relative Fourier-Mukai transform can be used to construct vector bundles (via the spectral cover approach) which satisfy the topological consistency conditions of heterotic string theory. We have also included a short subsection about the influence of mirror symmetry in the problem of the reconstruction of a variety out of its derived category.
Finally there is an appendix that collects some basic definitions and results regarding pure sheaves and Simpson stability.

\section{Aspects of derived categories}\label{s:dercat}
 
 Derived categories were introduced in the sixties in the framework of homological algebra and Grothendieck duality of coherent sheaves. One of the first accounts on the subject is a 1963 Verdier's booklet  reproduced in \cite{Ver} and expanded in Verdier's Ph.D. thesis (1967), which has been now reedited \cite{Ver96}. Since then, many readable accounts have been written, for instance \cite{Illu90} or quite recently \cite{Th00}.
 
\subsection{What is the derived category?}
Let us start with a complex algebraic variety $X$ whose structure ring sheaf we denote by $\Oc_X$\footnote{Definitions and properties in this section are true for schemes over an arbitrary algebraically closed field. Most of them remain true even for more general schemes.}. 
By sheaves we will always understand \emph{sheaves of $\Oc_X$-modules} and 
we will denote them by calligraphic  letters, like $\Fc$, as far as possible. 

A sequence 
$$
F \equiv \dots \xrightarrow{d_{i-2}} \Fc^{i-1}\xrightarrow{d_{i-1}}  \Fc^i \xrightarrow{d_{i}} \Fc^{i+1}\xrightarrow{d_{i+1}} \dots\qquad (i\in\Z)
$$
of morphisms of sheaves where the composition of any two consecutive maps is zero, is known as  a \emph{complex} of  sheaves. The morphisms $d_i$ are also known as \emph{differentials} of the complex $F$. The complex is \emph{bounded below} if it starts at some place $i$ (that is, the sheaves $\Fc^j$ are zero for $j<i$), \emph{bounded above} if it ends at some place $i$, and simply \emph{bounded} if it has only a finite number of non-vanishing sheaves. 

A complex has  \emph{cohomology sheaves} defined as 
$$
\Hc^i(F)=\Ker d_i/\Img d_{i-1}\,.
$$
Morphisms of complexes $\phi\colon F\to G$ are defined as collections $\phi_i\colon \Fc^i\to \G^i$ of morphisms commuting with the differentials.
A complex morphism induces morphisms $\Hc^i(\phi)\colon \Hc^i(F)\to \Hc^i(G)$. We then say that $\phi$ is a \emph{quasi-isomorphism} whenever all the induced morphisms $\Hc^i(\phi)$ are isomorphisms.

The derived category of  sheaves is a category built from complexes in such a way that quasi-isomorphisms become isomorphisms. 
This is accomplished as follows: We first identify two complex morphisms $\phi,\psi\colon F\to G$ when they are homotopically equivalent. This gives rise to a category, the category of complexes up to homotopies. 

The second step  is to ``localize'' by (classes of) quasi-isomorphims. This localization is a fraction calculus for categories: just think of the composition of morphisms as a product. Then quasi-isomorphisms verify the conditions for being a multiplicative system (that is, a nice set of denominators), namely, the identity is a quasi-isomorphism and the composition of two quasi-isomorphisms is a quasi-isomorphism.  Now, a fraction is nothing but a diagram of (homotopy classes of) complex morphisms
\begin{equation}
\xymatrix{ \bar F\ar[d]_{\phi}\ar[dr]^{\psi} & \\
F & G
}
\label{e:morph}
\end{equation}
where $\phi$ is a \emph{quasi-isomorphism}. We also have a notion of equivalence of fractions, we say that two ``fractions''
$$
\xymatrix{ \bar F\ar[d]_{\phi}\ar[dr]^{\psi} & \\
F & G
}\qquad
\xymatrix{ \bar F'\ar[d]_{\phi'}\ar[dr]^{\psi'} & \\
F & G
}
$$
are equivalent when there exist quasi-isomorphisms
$$
\xymatrix{ & \tilde F\ar[dl]_{\gamma}\ar[dr]^{\gamma'} & \\
\bar F & & \bar F'
}
$$
such that  $\phi\circ\gamma=\phi'\circ\gamma'$ and $\psi\circ\gamma=\psi'\circ\gamma'$. 

The \emph{derived category of complexes of sheaves} $D(X)$ is then defined as the category whose objects are complexes of quasi-coherent sheaves and whose morphisms are ``fractions'' like (\ref{e:morph}) where two equivalent fractions give rise to the same morphism.

By the very definition, complex quasi-isomorphisms become isomorphisms in $D(X)$ and quasi-isomorphic complexes (that is, complexes $F$ and $G$ such that there exists a diagram like \ref{e:morph} where both $\phi$ and $\psi$ are quasi-isomorphisms), become isomorphic.

Some other derived categories can be defined in a similar way:
\begin{enumerate}
\item The derived category $D^+(X)$ of bounded below complexes of sheaves,
$$
 F \equiv \Fc^m\to\Fc^{m+1}\to\dots \qquad (m\in\Z)\,.
$$
There is a functor $D^+(X)\to D(X)$, that is an equivalence between $D^+(X)$ and the subcategory of objects in  $D(X)$ defined by complexes which are quasi-isomorphic to bounded below complexes.  One can see that they are exactly those complexes $F$ whose homology sheaves $\Hc^i(F)$ are zero for all $i<m$  for a certain $m\in\Z$.
\item The derived category $D^-(X)$ of bounded above complexes,
$$
 F \equiv  \dots\to\Fc^{m-1}\to\Fc^{m}  \qquad (m\in\Z)\,.
$$
As above, $D^-(X)$ is equivalent to the subcategory of all objects in  $D(X)$ defined by complexes  quasi-isomorphic to bounded above complexes, or what amount to the same, to complexes  $F$ whose homology sheaves $\Hc^i(F)$ are zero for all $i\ge m$  for a certain $m\in\Z$.
\item The derived category $D^b(X)$ of bounded complexes of sheaves,
$$
 F \equiv  \Fc^m\to\Fc^{m+1}\to\dots\to\Fc^{s-1}\to\Fc^{s}  \qquad (m\le s\in\Z)\,.
$$
Again $D^b(X)$ is equivalent to the subcategory of all complexes in  $D(X)$ quasi-isomorphic to bounded complexes, or all complexes  $F$ whose homology sheaves $\Hc^i(F)$ are zero for all $i\notin [m,s]$ for some $m\le s\in\Z$.
\item The corresponding categories $D_{qc}(X)$, $D_{qc}^+(X)$, $D_{qc}^-(X)$, $D_{qc}^b(X)$ defined as above by using complexes of quasi-coherent sheaves. They are isomorphic to the subcategories of $D(X)$, for instance $D_{qc}^b(X)$ is equivalent to the subcategory of $D(X)$ defined by complexes of quasi-coherent sheaves whose cohomology sheaves are coherent and zero above and below certain finite indexes.
\item Finally, the corresponding categories $D_c(X)$, $D_c^+(X)$, $D_c^-(X)$, $D_c^b(X)$ defined as above using complexes of coherent sheaves instead of complexes of quasi-coherent sheaves. It turns out that they are isomorphic to the subcategories of $D_{qc}(X)$ and $D(X)$. For instance, $D_c^+(X)$ is equivalent to the subcategory of $D_{qc}^+(X)$ defined by bounded below complexes of quasi-coherent sheaves whose cohomology sheaves are coherent, and also to the subcategory of $D(X)$ of complexes whose cohomology sheaves $\Hc^i(F)$  are coherent and zero for all $i<m$ for $m\in\Z$.
\end{enumerate}

The derived categories we have defined are \emph{triangulated categories}. We are not giving the definition of what a triangulated category is. We just say that part of the notion of triangulated category is the existence of a translation functor. In the case of $D(X)$ (and of any of the other derived categories) that functor is
\begin{align*}
D(X)&\xrightarrow{\tau} D(X) \\
F&\mapsto \tau(F)=F[1]
\end{align*}
where for a complex $F$ and an integer number $i$, the complex $F[i]$ is the complex given by $F[i]^n=F^{n+i}$, that is, is the complex $F$ shifted $i$-places to the left.

\subsection{Derived functors in derived categories}

Derived functors are the ``derived category notion'' that corresponds to ``cohomolgy''. We know that sheaf cohomology groups are computed with the aid of resolutions. 
If we have a resolution
$$
0\to \Fc\to R\simeq\Rca^0\to\Rca^1\to\dots
$$
of a sheaf $\Fc$ by injective sheaves $\Rca^i$, the cohomology groups of $\Fc$ are defined as the cohomology groups
$$
H^i(X,\Fc)=H^i(\Gamma(X, R))
$$
of the complex
$$
\Gamma(X, R)\simeq \Gamma(X,\Rca^0)\to\Gamma(X,\Rca^1)\to\dots
$$
One proves that the definition of $H^i(X,\Fc)$ is well-posed, that is, it is independent of the injective resolution $R$.  Working in the derived category, we can see that this is equivalent to saying that whenever $R$ and $\bar R$ are injective resolutions of a sheaf $\Fc$, then the complexes $\Gamma(X, R)$ and $\Gamma(X, \bar R)$ are quasi-isomorphic, that is, $\Gamma(X, R)\simeq \Gamma(X, \bar R)$ in the derived category.

In that way we can associated to $\Fc$ a single object $R\Gamma(X,\Fc):=\Gamma(X, R)$ of the derived category.  This suggests that the derived category is the natural arena for cohomology constructions, such as derived functors. 

\subsubsection*{The derived direct image} 
Assume for instance that $f\colon X\to Y$ is a morphism of algebraic varieties (or schemes if you prefer so). As for cohomolgy, the higher direct images are defined as the cohomology sheaves $R^i f_\ast (\Fc)=\Hc^i( f_\ast(R))$, where
$$
0\to \Fc\to \Rca^0\to\Rca^1\to\dots
$$
is a resolution of $\Fc$ by injective sheaves of $\Oc_X$-modules. As for the cohomology groups,  we can generalize to the derived category the construction of the higher direct images by defining the \emph{right derived functor} of the direct image as the functor
\begin{align*}
Rf_\ast \colon D^+(X)& \to D^+(Y) \\
F&\mapsto Rf_\ast(F):= f_\ast (R)
\end{align*}
where $R$ is a bounded below complex of injective sheaves quasi-isomorphic to $F$. Such a complex $R$ always exists. 

In this way we can derive many functors. Sometimes, as we have seen, we can extend  functors defined for sheaves and taking values on sheaves as well; besides, we can also derive, that is, extend to the derived category,  functors defined only for complexes and taking values in complexes of sheaves.

Let us go back to the right derived direct image $Rf_\ast\colon D^+(X)\to D^+(Y)$. Under very mild conditions,\footnote{$f$ has to be quasi-compact locally of finte type so that the direct image of a quasi-coherent sheaf is still quasi-coherent.} $Rf_\ast$ maps complexes with quasi-coherent cohomology to complexes with quasi-coherent cohomolgy, thus defining a functor
$$
Rf_\ast\colon D_{qc}^+(X)\to D_{qc}^+(Y)\,,
$$
that we denote with the same symbol. When $f$ is proper, so that the higher direct images $R^if_\ast \Fc$ of a coherent sheaf $\Fc$ are coherent as well (cf. \cite[Thm.3.2.1]{EGAIII-1} or \cite[Thm. 5.2]{Hart77}, \cite{FAC} in the projective case), we also have a functor
$$
Rf_\ast\colon D_{c}^+(X)\to D_{c}^+(Y)\,.
$$

Finally, when the cohomological dimension of $f$ is finite, and this happens for instance when the dimensions of the fibers are bounded, then $Rf_\ast$ maps complexes with bounded cohomology to complexes with bounded cohomology, thus defining a functor
$$
Rf_\ast\colon D_{c}^b(X)\to D_{c}^b(Y)\,.
$$

Moreover, in this case the derived direct image can be extended to functors 
$$
Rf_\ast\colon D_{qc}(X)\to D_{qc}(Y)\,,\qquad Rf_\ast\colon D_{c}(X)\to D_{c}(Y)
$$
between the whole derived categories, which actually map $D_{qc}^b(X)$ to $D_{qc}^b(Y)$ and $D_{c}^b(X)$ to $D_{c}^b(Y)$. This follows essentially because every complex of sheaves, even if infinite on both sides, is still quasi-isomorphic to a complex of injective sheaves. 

This procedure is quite general, and applies with minor changes to other situations. We list a few relevant cases:

\subsubsection*{The derived inverse image}  
Again $f\colon X\to Y$ is a morphism of algebraic varieties. Once we know how the classical definition of higher direct images can be defined in terms of the derived category, we then also know how to do for the higher inverse images:
The higher inverse images are defined as the cohomology sheaves 
$$
L_i f_\ast (\Fc)=\Hc^{-i}( f^\ast(P))\,,
$$ where
$$
\dots\to \Pc^{-1}\to\Pc^0\to \Fc\to 0
$$
is a resolution of $\Fc$ by locally free sheaves. Then we define  the \emph{left derived functor} of the inverse image as the functor
\begin{align*}
Lf^\ast \colon D^-(X)& \to D^-(Y) \\
F&\mapsto Lf^\ast(F):= f^\ast (P)
\end{align*}
where $P$ is a bounded above complex of locally free sheaves quasi-isomorphic to $F$ (it always exists). It is very easy to check that $Lf^\ast$ defines also functors $Lf^\ast \colon D_{qc}^-(X) \to D_{qc}^-(Y)$ and $Lf^\ast \colon D_{c}^-(X) \to D_{c}^-(Y)$.

In some cases $Lf^\ast$ defines a functor
$$
Lf^\ast \colon D_c(X) \to D_c(Y)\,
$$
that maps $D_c^b(X)$ to $D_c^b(Y)$. One is when every coherent sheaf $\G$ on $Y$ admits a \emph{finite} resolution by coherent locally free sheaves, a condition that is equivalent to the smoothness of $Y$ (by the Serre criterion,cf. \cite{Ser65}). In such a case, every object in $D_c^b(Y)$ can be represented as a bounded complex of coherent locally free sheaves\footnote{The complexes that are quasi-isomorphic to a bounded complex of coherent locally free sheaves are known as \emph{perfect} complexes.}.

Note that, when $f$ is \emph{of finite Tor-dimension}, that is, when for every coherent sheaf $\G$ on $Y$ there are only a finite number of non-zero derived inverse images $L^if^\ast (\G)=Tor_i^{f^{-1}\Oc_Y}(f^{-1}(\G),\Oc_X)$. This happens of course when $f$ is flat, because no Tor can arise. 
In this case, if $F$ is a bounded complex of coherent sheaves and $P$ is a bounded above complex of locally free sheaves quasi-isomorphic to $F$, then $f^\ast (F)$ is quasi-isomorphic to $f^\ast (P)$, so that $L^\ast F=f^\ast F$ in $D^b(X)$.

\subsubsection*{Deriving the tensor product} 
Recall that given two complexes $F$, $G$ of sheaves, the tensor product complex $F\otimes G$ is defined by
$$
(F\otimes G)^n=\oplus_{p+q=n} F^p\otimes G^q
$$
with the differential $d$ that acts on $F^p\otimes G^q$ as  $d_F\otimes 1+(-1)^p1\otimes d_G$.

Take a fixed bounded above complex of  sheaves on an algebraic variety $X$. If $G$ is another  bounded above complex of  sheaves and $P$, $\bar P$ are two bounded above complexes of locally free sheaves quasi-isomorphic to $G$, then 
the simple complexes
$$
F\otimes P\,,\quad F\otimes\bar P
$$
associated to the bicomplexes whose $(p,q)$ terms are respectively $F^p\otimes P^q$ and $F^p\otimes \bar P^q$, are quasi-isomorphic. We can then define a left-derived functor
\begin{align*}
F\otimes \colon D^-(X)& \to D^-(X) \\
G&\mapsto F\underline\otimes G:= F\otimes P
\end{align*}
where $P$ is a bounded above complex of locally free sheaves quasi-isomorphic to $G$ (it always exists). One can now prove that if we fix the complex $P$ and consider two different bounded above complexes $Q$, $\bar Q$ of locally free sheaves quasi-isomorphic to $F$, then the simple complexes $$
Q\otimes P\,,\quad \bar Q\otimes P
$$
associated to the bicomplexes whose $(p,q)$ terms are, respectively, $Q^p\otimes P^q$ and $\bar Q^p\otimes P^q$, are quasi-isomorphic. In this way we can define the \emph{total left-derived functor} of the tensor product as the functor
\begin{align*}
\otimes^L \colon D^-(X)\times D^-(X)& \to D^-(X) \\
(F,G)&\mapsto F\otimes^L G:= Q\otimes P
\end{align*}
where $P$ is a bounded above complex of locally free sheaves quasi-isomorphic to $G$ and $Q$ is a bounded above complex of locally free sheaves quasi-isomorphic to $F$.

The total derived functor of the tensor product can be defined as well as a functor $\otimes^L \colon D_{qc}^-(X)\times D_{qc}^-(X) \to D_{qc}^-(X)$ and also as a functor $\otimes^L \colon D_{c}^-(X)\times D_{c}^-(X) \to D_{c}^-(X)$, analogously to what happens for the inverse image. Also as in this case, sometimes the derived tensor product can be extended as a functor between bigger derived categories; for instance, if $F$ is  quasi-isomorphic to a bounded complex of coherent locally free sheaves (a perfect complex), then one can define
 \begin{align*}
F \otimes^L \colon D(X)& \to D(X) \\
G&\mapsto F\otimes^L G
\end{align*}
and similarly for $D_{qc}(X)$ and $D_{c}$. These functors preserve the categories $D^{b}(X)$, $D_{qc}^b(X)$ and $D_{c}^b(X)$.

\subsubsection*{Deriving the homomorphisms}  
We have here two  types of complexes of homomorphisms, the global and the local ones. 

Let us consider first the global case. Recall that given two complexes $F$, $G$ of sheaves, the complex of (global) homomorphisms is the defined as the complex of abelian groups given by
$$
\Hom_X^n(F,G)=\prod_{p} \Hom_X(F^p,G^{p+n})
$$
with differential $df=f\circ d_F+(-1)^{n+1}d_G\circ f$. Here the key point is that if $I$,$\bar I$ are two quasi-isomorphic bounded below complexes of injective sheaves, then the complexes of homomorphisms $\Hom_X^\bullet(F,I)$ and $\Hom_X^\bullet(F,\bar I)$ are quasi-isomorphic as well whatever the complex $F$ is. Also, if we fix a  bounded below complex $I$ of injective sheaves and $F$, $\bar F$ are quasi-isomorphic, then $\Hom_X^\bullet(F,I)$ and $\Hom_X^\bullet(\bar F,I)$ are still quasi-isomorphic. We can then define a right derived functor 
\begin{alignat*}{2}
R\Hom_X \colon D(X)^0&\times D^+(X)&& \to D(Ab) \\
(F&,G)&&\mapsto R\Hom_X(F,G):=\Hom_X^\bullet(F,I)
\end{alignat*}
where $I$ is any bounded below complex of injective sheaves quasi-isomorphic to $G$. Here $D(Ab)$ stands for the derived category of abelian groups. We can also define a left derived functor
\begin{alignat*}{2}
L\Hom_X\colon D^-(X)^0&\times D(X)&& \to D(Ab) \\
(F&,G)&&\mapsto L\Hom_X(F,G):=\Hom_X^\bullet(P,G)
\end{alignat*}
where $P$ is any bounded above complex of locally-free sheaves quasi-isomorphic to $F$. Both functors coincide over  $D^-(X)^0\times D^+(X)$.

One can define the Ext groups for objects $F$, $G$ of the derived category, they are defined as the groups
$$
\Ext_X^i(F,G):=H^i(R\Hom_X(F,G)) \quad\text{for every $i\in\Z$}\,,
$$
and they are defined whenever the second right hand has a sense. There is a nice formula which allow to compute those Ext's and most of their properties:
\begin{prop} Assume that either $F$ is in $D(X)$ and $G$ in $D^+(X)$ or $F$ is in $D^-(X)$ and $G$ in $D(X)$. Then one has
$$
\Ext_X^i(F,G)=\Hom_{D(X)}(F,G[i]) \quad\text{for every $i\in\Z$}\,.
$$
\label{p:exti}
\qed\end{prop}

We can derive the homomorphism sheaves as well. The procedure is the same, considering now the complex of sheaves
$$
\mathcal{H}om_{\Oc_X}^n(F,G)=\prod_{p} \mathcal{H}om_{\Oc_X}(F^p,G^{p+n})
$$
with differential $df=f\circ d_F+(-1)^{n+1}d_G\circ f$ as above. The result here is a right derived functor 
\begin{alignat*}{2}
R\mathcal{H}om_{\Oc_X} \colon D(X)^0&\times D^+(X)&& \to D(X) \\
(F&,G)&&\mapsto R \mathcal{H}om_{\Oc_X}(F,G):=\mathcal{H}om_{\Oc_X}^\bullet(F,I)
\end{alignat*}
where $I$ is any bounded below complex of injective sheaves quasi-isomorphic to $G$.

On readily checks that the derived homomorphism sheaves preserve the categories $D_{qc}(X)$ and $D_c(X)$ in the sense one naturally thinks of.  In some more precise terms it induces functors
$$
R\mathcal{H}om_{\Oc_X} \colon D_{qc}(X)^0\times D_{qc}^+(X) \to D_{qc}(X) \,,\quad
R\mathcal{H}om_{\Oc_X} \colon D_{c}(X)^0\times D_{c}^+(X) \to D_{c}(X). 
$$

We can now define the ``dual'' of an object $F$ in any of the derived categories $D(X)$, $D_{qc}(X)$ or $D_c(X)$. It is the object $R\mathcal{H}om_{\Oc_X}(F,\Oc_X)$. Of course even if $F$ reduces to a sheaf $\Fc$, the dual $R\mathcal{H}om_{\Oc_X}(\Fc,\Oc_X)$ may fail to be a sheaf. It is represented by a complex whose $(-i)$-th cohomology sheaf is the Ext-sheaf $\mathcal{E}xt^i_{\Oc_X}(\Fc,\Oc_X)$.

The relationship between derived homomorphism sheaves and groups is very easy, and it is a particular case of what is known as Grothendieck theorem on the composition functor. One simply has that
$$
R Hom(F,G)=R\Gamma(X,R\mathcal{H}om_{\Oc_X}(F,G))
$$
when $F$, $G$ reduce to single sheaves $\Fc$, $\G$, the above equality means that there exists a spectral sequence whose term $E_2^{pg}$ is $H^p(X, \mathcal{E}xt^q_{\Oc_X}(\Fc,\G))$ converging to $E_\infty^{p+q}=\Ext^{p+q}_X(\Fc,\G)$.

Many formulae like the above one, relating different ordinary derived functors by means of spectral sequences, can be also stated in a very clean way by the aid of derived categories. We have base change theorems, projection formulae and many others. The reader is referred for instance to \cite{Hart66}.

\subsubsection*{Chern classes in derived category}

When $X$ is a smooth projective variety, any complex $M$ in the derived category $D(X)$ is isomorphic to a bounded complex 
$$
M\simeq E\equiv \Ec^s\to\Ec^{s+1}\to\dots\to \Ec^{s+n}
$$
 of locally free sheaves (i.e. it is \emph{perfect} or \emph{of finite Tor-dimension}). 

The Chern characters of $M$  are then defined by 
$$
\ch_j(M )=\sum_i (-1)^i \ch_j(\Ec^i)\in A^j(X)\otimes \bbQ\,,
$$
 where $A^j(X)$ is the $j$-th component of the Chow ring (when $k = \bbC$, the group  $A^j(X)\otimes \bbQ$ is the algebraic part of the rational cohomology group $H^{2j} (X,\bbQ)$). This definition is well posed since it is independent of the choice of the bounded complex $E$ of locally free sheaves.  

By definition the rank of $M$ is the integer number 
\begin{equation}
\rk(M )=\ch_0(M ).
\label{e:rank}
\end{equation}

We shall see however a different definition of rank of a sheaf, namely the polarized rank (Definition \ref{d:polrank}), and we will find the relationship between both.

\section{\GIF s and \FM s}\label{s:gifs}
We now define Fourier-Mukai transforms, or more generally \emph{geometrical integral functors}. 

We will always refer to algebraic varieties as a synonymous for the more technical ``shemes of finite type over an algebraically closed field $k$''.  You may think that $k$ is the field of complex numbers if you feel more comfortable; however, the characteristic of $k$ does not play any roll  and then can be arbitrary.  A sheaf on $X$ is always assumed to be coherent.

For a scheme $X$ we denote by $\Oc_X$ it structure sheaf and by $\Oc_{X,x}$ or $\Oc_x$ the local ring of $\Oc_X$ at $x$. The ideal sheaf of $x$  will be denoted by $\mg_x$ and $\kappa(x)=\Oc_x/\mg_x$ is the residue field at the point.

By the sake of simplicity we will simply write $D(X)$ for the bounded derived category of coherent sheaves. 
\subsection{\GIF s}

Let $X$, $\X$ be proper algebraic varieties; 
the projections of the cartesian product $X\times \X$ onto the factors $X$, $\X$ are denoted, respectively, by $\pi$, $\hat\pi$. We can put this information into a diagram
$$
\xymatrix{
& X\times \X \ar[ld]_{\pi}
\ar[dr]^{\hat\pi} \\
X & & \X}
$$

Let $E$ be an object in the derived category 
$D(X\times \X)$. We shall call it a ``kernel'' and define a \emph{\gif}  between the
derived  categories by
$$
\bS^{E}\colon D(X)\to D(\X) \,,\quad F\mapsto \bS^{E}(F)=
R\hat\pi_{\ast}(L\pi^\ast F\otimes E)
$$ 
(the tensor product is made in the derived category). That is, we first pull back $F$ to $X\times\X$, then twist with the kernel $E$ and push forward to $\X$.

We shall call \emph{\FMF s} those \gif s that are \emph{equivalences of categories} between $D(X)$ and $D(\X)$ and \emph{\FM s}, the \FMF s whose kernel $E$ is a single sheaf.

If $\bS^{E}\colon D(X)\to D(\X)$ is a \FMF, it preserves the homomorphism groups. Since the Ext groups in derived category are defined as 
$$
\Ext_X^i(F,G)=\Hom_{D(X)}(F,G[i]
$$
 by Proposition \ref{p:exti}, we get

\begin{prop} \emph{(Parseval formula)} $\bS^{E}\colon D(X)\to D(\X)$ is a \FMF, then for every $F$ and $G$ in $D(X)$ one has
$$
\Ext_X^i(F,G)\simeq \Ext_{\X}^i(\bS^{E}(F),\bS^{E}(G))
$$
for all indexes $i$.
\label{p:parseval}
\qed
\end{prop}

\subsubsection*{WIT and IT conditions}

An important feature of \gif s  is that they are \emph{exact}
as  functors of triangulated categories. In more familiar terms we can say that for any exact sequence $0\to\Nc \to \Fc \to \G \to 0$
of coherent sheaves in $X$ we obtain an exact sequence
\begin{equation}
\dots \to \bS^{i-1} (\G)\to \bS^i(\Nc)\to \bS^i(\Fc) \to \bS^i(\G) \to \bS^{i+1} (\Nc)\to \dots
\label{1:exactseq}
\end{equation}
where we have written $\bS=\bS^{E}$ and $\bS^i(F)=\Hc^i(\bS(F))$.

\begin{defin}\label{1:WIT} Given an \gif\ $\bS^{E}$,
a complex $F$ in $D(X)$  satisfies the WIT$_i$  condition
(or is WIT$_i$) if there is a coherent sheaf $\G$ on $\X$ such that 
 $\bS^{E}(F)\simeq\G[i]$ in  $D(\X)$, where $\G[i]$ is the associated complex concentrated in degree $i$.
We say that $F$  satisfies the IT$_i$  condition if in
addition $\G$ is locally free. 
\end{defin}

When the kernel $E$ is simply a sheaf $\Qc$ on $X\times \X$ flat over $\X$,  by cohomology and base change theorem \cite[III.12.11]{Hart77}  one has
\begin{prop}
A coherent sheaf $\Fc$ on $X$ is IT$_i$ if and only if $H^j(X, \Fc\otimes \Qc_\xi) = 0$ for all $\xi\in \X$ and for all $j \ne i$, where $\Qc_\xi$ denotes the restriction of $\Qc$ to $X\times\{\xi\}$. Furthermore, $\Fc$ is WIT$_0$ if and only if it is IT$_0$.
\label{1:WITvsIT} 
\qed\end{prop}

The acronym ``IT" stands for ``index theorem'',  while ``W'' stands for ``weak''. This terminology comes from Nahm transforms for connections on tori in complex differential geometry. 

A systematic and comprehensive treatment of \gif s and \FM s is to appear in the book \cite{BBHJ}.


\subsubsection*{The original \FM}

Mukai introduced the first \FM\ in the framework of abelian varieties, we refer to \cite{Mum74} or \cite{LB} for very readable accounts on abelian varieties. Abelian varieties are simply proper algebraic groups; however, properness implies commutativity which explains the terminology. From a differential geometric viewpoint, a complex abelian variety $X$ of dimension $g$ is a complex torus, $X\iso \bbC^g/\Lambda$, with $\Lambda\iso \Z^g$ being a lattice.

The play of the second variety $\X$ is played by the ``dual'' abelian variety. This is described algebraically as the variety parametrizing line bundles of degree zero on $X$, or analytically as $\X=\bbC^g/\Lambda^\vee$, where $\Lambda^\vee$ is the dual lattice. 

The kernel is the Poincar\'e line bundle $\Pc$ on $X\times\X$. This is the universal line bundle of degree zero, and  it is characterized by the property that its restriction to $X\times\{\xi\}$ where $\xi\in\X$ is precisely the line bundle $\Pc_\xi$ on $X$ defined by $\xi$.  This characterization determines $\Pc$ only up to tensor products by inverse images of line bundles on $\X$.  It is customary to normalize $\Pc$ such that  its restricition to $\{0\}\times \X$ is trivial (here $0$ denotes the origein of the abelian variety $X$). 

 Mukai's seminal idea \cite{Muk81} was to use the normalized Poincar\'e bundle $\Pc$ 
to define an integral functor between the derived categories
$$
\bS^{\Pc} \colon D(X) \to D(\X)
$$
which turns out to be an equivalence of triangulated categories, or in our terminology, a \FM.

\subsubsection*{Compostion of \gif s}
The composition of two \gif s is still a \gif\ whose kernel  can be  expressed as a kind of ``convolution product'' of the two original
kernels.  If $\tilde X$ is a  third 
proper variety, let us take a kernel $\widetilde E$ in $D(\X\times\tilde X)$ and consider the \gif\ 
$$
\bS^{\widetilde E}\colon D(X)\to D(\X) \,,\quad G\mapsto \bS^{\widetilde E}(G)=
R\tilde\pi_{\ast}(L\hat\pi^\ast G\otimes\widetilde E)
$$ 
We now consider the diagram 
$$\xymatrix{
& X\times \X \times \tilde X \ar[ld]_{\pi_{1,2}} \ar[d]^{\pi_{2,3}}
\ar[dr]^{\pi_{1,3}} \\
X\times \X & \X\times \tilde X & X\times \tilde X}$$

\begin{prop} The composition of the \gif s $\bS^{E}\colon D(X)\to D(\X)$ and $\bS^{\widetilde E}\colon D(\X)\to D(\widetilde X)$ is the \gif\ 
$D(X)\to D(\widetilde X)$ with kernel $R \pi_{13,\ast}(L\pi_{12}^\ast E\otimes L\pi_{23}^\ast \widetilde E)$, that is,
$$
\bS^{\widetilde E}(\bS^{E}(F))= R\pi_{2,\ast}(\pi_1^\ast(F)\otimes R \pi_{13,\ast}(L\pi_{12}^\ast E\otimes L\pi_{23}^\ast \widetilde E))\,.
$$
\label{p:composition}
\end{prop} 
The proof (see \cite{Muk81} for the original Fourier-Mukai transform or the book \cite{BBHJ})  is an standard exercise in derived category (base-change and projection formula).
%


\subsubsection*{\FMF s}\label{ss:fmfs}

We now give a few elementary examples of \gif s and \FMF s:
\begin{enumerate}
\item
Let $E$ be the complex in $D(X\times X)$ defined by the 
the structure sheaf  $\Oc_{\Delta}$ of the diagonal $\Delta\subset X\times X$.
Then it is easy to check that $\bS^{E}\colon D(X)\to D(X)$ is isomorphic to the
identity functor on $D(X)$. 

If we shift degrees by $n$ taking $E=\Oc_{\Delta}[n]$ (a complex with only the sheaf $\Oc_\Delta$ placed in degree $n$), then $\bS^{E}\colon D(X)\to D(X)$ is the degree shifting functor $\G\mapsto \G[n]$.
\item
More generally,  given a proper morphism $f\colon X\to \X$, 
by taking as $E$ the structure sheaf of the
graph $\Gamma_f\subset X\times Y$, one has isomorphisms
of functors $\bS^{E}\simeq R f_\ast$ as functors  $D(X)\to D(\X)$ and 
$\bS^{E}\simeq f^\ast$ now as functors $D(\X)\to D(X)$.
\item  Take $\X=X$ and let $\Lcl$ be a line bundle on $X$. If $E=\hat\pi^\ast \Lcl$, then $\bS^{E}(G)=G\otimes\Lcl$ for any $G$ in $D(X)$. 
\end{enumerate}

Actually we can not find examples of equivalences of derived categories other than \FMF s. This is due to the following Orlov's crucial result \cite{Or97}.
\begin{thm} Let $X$ and $\X$ be smooth projective varieties. Any fully faithful functor $D(X)\to D(\X)$ is a \gif. In particular, any equivalence of categories $D(X)\iso D(\X)$ is a \FMF.
\label{th:orlov}
\end{thm}

It is very interesting to characterize for which kernels $E$ in $D(X\times\X)$ the corresponding \gif \  $\bS^{E}$ is an equivalence of categories, or in our current terminology, a \FMF.

For simplicity, we consider only the case of kernels that reduce to a single sheaf $\Qc$.  This covers an important number of relevant situations, including the original \FM\ whose kernel is the Poincar\'e line bundle $\Pc$ on the product of an abelian variety $X$ and it dual variety $\X$. 

Bridgeland-Maciocia \cite{Mac96} extracted the properties of $\Pc$ that make the \gif\ $\bS^{\Pc}\colon D(X) \to D(\X)$ into a \FM\ and introduced strongly simple sheaves.
\begin{defin} A coherent sheaf $\Pc$ on $X\times \X$ is \emph{strongly simple over $X$} if it is flat over $X$ and satisfies the following two conditions
\begin{enumerate}
\item $\Ext_{D(\X)}^i (\Pc_{x_1},\Pc_{x_2}) = 0$ for every $i\in \Z$ whenever $x_1\ne x_2$;
\item $\Pc_x$ is simple for every $x\in X$, i.e., all its automorphisms are constant multiples of the identity,  $\Hom_{\Oc_X} (\Pc,\Pc) =k$. 
\end{enumerate}
\end{defin}

The relevant result   is
\begin{prop} Let $X$, $\X$ be smooth proper algebraic varieties and let $\omega_{\X}$ be the canonical line bundle of $\X$. If $\Pc$ is a sheaf on $X\times \X$ strongly simple over $X$, then $\bS^{\Pc}\colon D(X)\to D(\X)$ is  a \FM\ if $\dim X=\dim \X$ and $\Pc_x\otimes\omega_{\X}\simeq
\Pc_x$ for all $x\in X$. Moreover, $\Pc$ is strongly simple over $\X$ as well and $\bS^{\Pc}\colon D(\X)\to D(X)$ is also a \FM.
\label{p:invertibility}
\qed
\end{prop}

The inverse \gif\ of $\bS^{\Pc}\colon D(X)\to D(\X)$ is the \FMF\ $\bS^{\Qc}\colon D(\X)\to D(X)$ whose kernel is the complex
$$
\Qc= \Pc^\ast \otimes \pi^\ast \omega_X [n]
$$
where $n=\dim X$.
The proof of the above Proposition is based on the description of the composition of two \gif s given by Proposition \ref{p:composition} plus a technical argument that ensures that the convolution $\Pc * \Qc$ is actually the structure sheaf $\Oc_\Delta$ of the diagonal on $X\times X$ (see \cite{Bri99} or \cite{BBHJ}).

Earlier results about the invertibility of certain \gif s can be proved in much simpler way using Proposition \ref{p:invertibility}. We mention here just the two who historically came first:
\begin{enumerate}
\item The \emph{original \FM} already described in a precedent section. It is very easy to see that the Poincar\'e line bundle $\Pc$ is strongly simple over $\X$. Since both the canonical bundles of $X$ and $\X$ are trivial, Proposition \ref{p:invertibility} gives that $\Pc$ is also strongly simple over $X$ and that the original \FM
$$
\bS^{\Pc}\colon D(X)\to D(\X)
$$
is indeed a \FM, that is, an equivalence of categories.
\item The \FM\ for \emph{reflexive K3 surfaces} (\cite{BBH97a}). Here $X$ is a K3 surface with a polarization $H$ and a divisor $\ell$ such that $H^2=2$, $H\cdot\ell=0$,
$\ell^2=-12$, and $\ell+2H$ is not effective. We take $\X$ as the fine moduli space of stable sheaves (with respect to $H$) of rank $2$, $c_1=\ell$ and Euler characteristic equal to $-1$. This topological invariants are selected so that earlier theorems of Mukai (\cite{Muk87a}) ensure that $\X$ is another K3 surface (which turns out to be isomorphic with $X$). The kernel of the \gif is the universal rank $2$ sheaf $\Pc$ on $X\times\X$ (suitably normalized).

As in the former example, one can readily prove, using that $\Pc_{\xi}$ is a rank 2 stable vector bundle of degree zero on $X$ for every point $\xi\ n\X$, that $\Pc$ is strongly simple over $\X$. Again the canonical bundles of $X$ and $\X$ are trivial, so that again by Proposition \ref{p:invertibility} one has that the \gif\ 
$$
\bS^{\Pc}\colon D(X)\to D(\X)
$$
is a \FM.
\end{enumerate}

\subsection{Relative \gif s}

From the old days of Grothendieck, algebraic geometers use to consider the ``relative'' situation, that is, they study problems for families rather than for single varieties.  As we said before, we can then better encode base-change properties into the problem. 

We can do that for \gif s as well.  To this end, we consider two morphisms $p\colon X\to B$, $\hat p\colon\X\to B$ of algebraic varieties. We shall define a 
relative \gif\ in this setting by means of a ``kernel'' $E$ in the derived category $D(X\times_B\X)$, just by mimicking the ``abosulte'' definition we already gave.

And since we want \gif s for families, we don't content ourselves with this setting and go beyond allowing further changes in the base space $B$, that is, we consider  base-change  morphisms $g\colon S\to B$ we denote all objects obtained by base
change to $S$ by a subscript $S$, like $X_S=S\fp X$ etc. In particular, the kernel $E$ defines an object $E_S=Lg^\ast E\in D((X\fp\X)_S)=D(X_S\times_S\X_S)$.

There is then a diagram
$$
\xymatrix{\save-<2.3truecm,0pt>*{(X\fp\X)_S\simeq }\restore{}
X_S\times_S\X_S \ar[d]^{\pi_S}\ar[r]^{\hat\pi_S}&\X_S
\ar[d]_{\hat p_S}\\ X_S \ar[r]^{p_S}& S}
$$ 
and the \emph{relative \gif} associated to $E$ is the functor between the
derived  categories of quasi-coherent sheaves given by
$$
\bS^{E_S}\colon D(X_S)\to D(\X_S) \,,\quad F\mapsto \bS^{E_S}(F)=
R\hat\pi_{S\ast}(L\pi_S^\ast F\otimes E_S)
$$ 
(the tensor product is made in the derived category). When $\hat p$ is a flat morphism, $\pi_S$ is flat as well and we can simply write $\pi_S^\ast F$ instead of $L\pi_S^\ast F$. 

We should not be scared by this new definition, because we immediately note that the  the relative \gif\  with respect to $E\in D(X\times_B\X)$ \emph{is nothing but the absolute \gif\  with kernel $i_\ast E\in D(X\times\X)$, where $i\colon X\fp\X\hookrightarrow X\times\X$ is the immersion}. The gain is that we can state neatly the following base-change property: 
\begin{prop} Let $F$ be an object in $D(\X_S)$. For every morphism
$g\colon S'\to S$ there is an isomorphism
$$ Lg_{\X}^\ast(\bS^{E_S}(F))\simeq \bS^{E_{S'}}(L g_X^\ast F)
$$ in the derived category $D^-(\X_{S'})$,  where
$g_X\colon X_{S'}\to X_S$, $g_{\X}\colon\X_{S'}\to\X_S$ are the
morphisms induced by $g$.
\label{p:basechange}
\end{prop}
\begin{proof} We give the proof as an easy example of  standard properties of derived categories. We have
$$
Lg_{\X}^\ast(\bS^{E_S}(F))=Lg_{\X}^\ast(R\hat\pi_{S\ast}(L\pi_S^\ast F\otimes E_S))=
R\hat\pi_{S'\ast}(Lg_{X\times_B\X}^\ast(L\pi_S^\ast F\otimes E_S))
$$
by base-change in the derived category, and then
$$
Lg_{X\times_B\X}^\ast(L\pi_S^\ast F\otimes E_S)=L\pi_{S'}^\ast(L g_X^\ast F)\otimes Lg_{X\times_B\X}^\ast(E_S)= L\pi_{S'}^\ast(L g_X^\ast F)\otimes E_{S'} \,.
$$
\end{proof}

Due to this property we shall very often drop the subscript $S$ and 
refer only to $X\to B$. 

\subsubsection*{Base change and WIT$_i$-conditions} 
In this paragraph we assume that the kernel is a single sheaf $\Pc$ flat over $B$ and that both $p$ and $\hat p$ are proper flat morphisms of relative dimension $n$.  We study the relationship between the WIT condition for a sheaf $\Fc$ on $X$ with respect to the relative \gif\ defined by $E$ in $D(X\fp\X)$ and the WIT condition for the restrictions $\Fc_s$ of $\Fc$ to the fibers $X_s=p^{-1}(s)$  ($s\in B$) with respect to the restriction.

Let us write $\bS=\bS^{E}$ and $\bS_s=\bS^{E_s}$. 

\begin{corol} Let $\Fc$ be a sheaf on $X$, flat over $B$. 
\begin{enumerate}
\item The formation of $\bS^n(\Fc)$ is compatible with base change,  that
is, one has $\bS^n(\Fc)_{s}\break \simeq
\bS_{s}^n(\Fc_{s})$,  for every point $s\in B$.
\item There is a convergent spectral sequence
$$
E_2^{-p,q}= {\mathcal T}or_p^{\Oc_S}(\bS^q(\Fc),\kappa(s)) \implies \bS^{q-p}(\Fc_s)\,.
$$
\item Assume that
$\Fc$ is WIT$_i$ and let $\what\Fc=\bS^i(\Fc)$ be its Fourier-Mukai
transform. Then for every $s\in B$ there are isomorphisms
$$ {\mathcal T}or_j^{\Oc_S}(\what\Fc,\kappa(s))\simeq \bS_s^{i-j}(\Fc_s)\,,\quad j\le i
$$ of sheaves over $\X_s$. In particular $\what\Fc$ is flat over $B$ if 
and only if the restriction $\Fc_s$ to the fiber $X_s$ is WIT$_i$ for 
every point $s\in B$.
\end{enumerate} 
\label{c:basechange}
\end{corol} 
\begin{proof} (1) follows from the fact that the highest direct image is compatible with base change (\cite{Hart77}). (2) is a consequence of Proposition \ref{p:basechange} and implies (3).
\end{proof}
\begin{corol} Let  $\Fc$ be a  sheaf on $X$, flat over $B$. There exists
an open subscheme $V\subseteq B$ which is the largest subscheme 
$V$ fulfilling one of the following equivalent conditions hold:
\begin{enumerate}
\item $\Fc_{V}$ is WIT$_i$ on $X_{V}$ and the \gif\ 
$\what\Fc_{V}$ is flat over $V$.
\item The sheaves $\Fc_s$ are WIT$_i$ for every point $s\in V$.
\end{enumerate}
\label{c:wit1locus}
\qed\end{corol}

\section{Elliptic fibrations}\label{s:elliptic}

An elliptic fibration is a proper flat
morphism $p\colon X\to  B$ of schemes whose fibers are Gorenstein curves of arithmetic genus 1.  We also assume that 
$p$ has a \emph{section}
$\sigma\colon B\hookrightarrow X$ taking values in the smooth locus
$X'\to B$ of $p$.  The generic fibers are then smooth elliptic curves whereas some singular fibers are allowed. If the base $B$ is a smooth curve, elliptic fibrations were studied and classified by Kodaira \cite{kod}, who described all the types of singular fibers that may occur, the so-called Kodaira curves. When the base is a smooth surface, more complicated configuration of singular curves can occur (see Miranda \cite{Mir83}). 

When $X$ is a Calabi-Yau threefold, the presence of the section imposes constraints to the base surface $B$; it is known that it  has to be of a particular kind,
namely $B$ has to be a Del Pezzo surface (a surface whose anticanonical divisor $-K_B$ is ample), a Hirzebruch surface (a rational ruled surface), a Enriques surface (a minimal surface with $2 K_B$ numerically equivalent to zero) or a blow-up of a
Hirzebruch surface (see for instance \cite{DLO}\ or \cite{MV96b}).

We denote by  $\Theta=\sigma(B)$  the image of the section, by $X_t$ the fiber of $p$ over $t\in B$ and by $i_t\colon X_t\hookrightarrow X$ the inclusion.
$\omega_{X/B}$  is the relative dualizing sheaf and we write $\omega=R^1p_\ast\Oc_X\iso (p_\ast\omega_{X/B})^\ast$, where the isomorphism is Grothendieck-Serre duality for $p$ (cf. \cite{De}).

The sheaf $\Lcl=p_\ast\omega_{X/B}$ is a line bundle  and
$\omega_{X/B}\iso p^\ast \Lcl$.  We write $\bar K=c_1(\Lcl)$. Adjunction formula for $\Theta\hookrightarrow X$ gives that
$\Theta^2=-\Theta\cdot p^{-1}\bar K$ as cycles on $X$.

\subsection{Weiestrass models and Todd classes}
We now recall some facts about the Weierstrass model for an elliptic fibration $p\colon X\to B$ with a section $\sigma$. If $B$ is a smooth curve, then from Kodaira's classification of possible singular fibers \cite{kod} one 
finds that the components of reducible fibers of $p$ which do not meet $\Theta$ form 
rational double point configurations disjoint from $\Theta$. Let $X\to {\bar{X}}$ be the result
of contracting these configurations and let $\bar{p}\colon {\bar{X}}\to B$ be the induced map.
Then all fibers of ${\bar p}$ are irreducible with at worst nodes or cusps as singularities. In this case 
one refers to ${\bar{X}}$ as the Weierstrass model of $X$.  The Weierstrass model can be constructed as follows:  the
divisor $3\Theta$ is relatively ample and if
$\Ec=p_\ast \Oc_X(3\Theta)\iso\Oc_B\oplus\omega^{\otimes
2}\oplus\omega^{\otimes3}$ and $\bar p\colon P=\Ps(\Ec^\ast)=\op{Proj}(S^\bullet(\Ec))\to B$ is the associated projective bundle, there is a projective morphism of
$B$-schemes $j\colon X\to P$ such that $\bar X=j(X)$. 

By the sake of simplicity we shall refer to a particular kind of elliptic fibrations, namely elliptic fibrations with a section as above whose fibers are all \emph{geometrically integral}. This means that the fibration is isomorphic with its Weierstrass model.
Now, special fibers can have at most one singular point, either a cusp or a simple node.  
Thus, in this case $3\Theta$ is relatively very ample and gives rise to a closed immersion
$j\colon X\hookrightarrow P$  such that $j^\ast\Oc_P(1)=\Oc_X(3\Theta)$.
Moreover $j$ is locally a complete intersection whose normal  sheaf is
\begin{equation}
\Nc(X/P)\iso p^\ast\omega^{-\otimes 6}\otimes\Oc_X(9\Theta)\,.
\label{eq:normal}
\end{equation} 
This follows by relative duality since
$$
\omega_{P/B}=\bigwedge \Omega_{P/B}\iso
\bar p^\ast\omega^{\otimes5}(-3)\,,
$$
 due to the Euler exact sequence 
$$
0\to\Omega_{P/B}\to \bar p^\ast \Ec(-1)\to\Oc_P\to 0\,.
$$
The morphism $p\colon X\to B$ is then a l.c.i. morphism in the sense of \cite[6.6]{Ful} and has a virtual relative
tangent bundle $T_{X/B}=[j^\ast T_{P/B}]-[\Nc_{X/P}]$ in the $K$-group
$K^\bullet (X)$. Even if  $T_{X/B}$ is not a true sheaf, it still has Chern classes; in particular, it has a Todd class which one can readily compute \cite{HMP02}. 
\begin{prop} The Todd class of the virtual tangent bundle $T_{X/B}$ is
$$
\td(T_{X/B})=1-\tfrac12\, p^{-1}\bar K+ \frac{1}{12}(12\Theta\cdot 
p^{-1}\bar K+13p^{-1}\bar K^2)-\frac 12 \Theta\cdot p^{-1} \bar K^2+\text{ terms of higher degree.}
$$\label{p:todd}
\qed\end{prop}

Since  $\td(B)=1+\frac12 c_1(B)+\frac1{12}(c_1(B)^2+c_2(B))+\frac1{24}c_1(B)c_2(B)+\dots$, we obtain the expression for the Todd class of $X$
\begin{equation}
\begin{aligned}
\td(X) & = 1+\frac12 p^{-1}(c_1(B)-\bar K) \\
& + \frac1{12}(12\Theta\cdot 
p^{-1}\bar K+13p^{-1}\bar K^2-3 p^{-1}(c_1(B)\cdot \bar K)+p^{-1}(c_1(B)^2+c_2(B))) \\
& + \frac1{24}[ p^{-1}(c_1(B)c_2(B))-p^{-1}(\bar K\cdot(c_1(B)^2+c_2(B)))+12 \Theta\cdot p^{-1}(\bar K\cdot c_1(B)) \\
& \phantom{+ \frac1{24}(}+ p^{-1}(c_1(B)\cdot \bar K^2)-6\Theta\cdot p^{-1}(\bar K^2\cdot c_1(B)) ] \\
&+\text{terms of higher degree.}
\end{aligned}
\label{eq:toddx}
\end{equation}
\subsection{Relative \gif s  for elliptic fibrations}

There is an algebraic variety $\hat p\colon\X \to B$ (the Altman-Kleiman compactification of the relative Jacobian) whose points parametrize torsion-free, rank one and degree zero sheaves on the fiber of $X\to B$. 
Moreover, the natural morphism of
$B$-schemes 
\begin{align*}
X & \to \X \\ 
x &\mapsto {\mathfrak
m}_x\otimes\Oc_{X_s}(e(s))
\end{align*}
 is an \emph{isomorphism}. Here  ${\mathfrak m}_x$ is the ideal sheaf of the point
$x$ in $X_s$. 

The variety $\X$ is a fine moduli space. This means that  there exists a coherent sheaf $\Pc$
on $X\fp\X$ flat over $\X$, whose restrictions to the fibers of $\hat\pi$ are torsion-free, and
of rank one and degree zero, such that for every morphism $f\colon S\to B$  and every sheaf $\Lcl$ on $X\fp S$ flat over $S$ and whose restrictions to the fibers of $p_{S}$ are torsion-free, of rank one and degree zero, there exists a unique morphism $\phi\colon X\to\X$ of $S$-schemes such that $\Lcl$ and $(1\times\phi)^\ast\Pc$ are isomorphic when restricted to every fiber of $X_S\to S$. 

The sheaf $\Pc$ is
defined up to tensor product by the pullback of a line bundle on $\X$, and
is called the \emph{universal Poincar\'e sheaf}.

Hereafter we  identify $X\iso\X$. Now $\Pc$ is a sheaf on $X\fp X$ that we can normalize by letting
\begin{equation}
\rest\Pc{\Theta\fp \X}\simeq \Oc_{X}\,.\label{e:norm}
\end{equation}
We shall
henceforth assume that $\Pc$ is normalized in this way so that
 \begin{equation}
\Pc=\Ic_\Delta\otimes \pi^\ast\Oc_X(\Theta)\otimes \hat\pi ^\ast\Oc_X(\Theta)\otimes q^\ast \omega^{-1}
\label{e:poin}
\end{equation}
where $\pi$, $\hat\pi$ and $q=p\circ\pi=\hat p\circ\hat\pi$ refer to the diagram
\begin{equation}
\xymatrix{
X\times_B X \ar[d]^{\pi}\ar[r]^{\hat\pi}\ar[rd]^q& X
\ar[d]_{\hat p}\\ X \ar[r]^{p}& B}
\label{e:ellipdiagram}
\end{equation}
and $\Ic_\Delta$ is the ideal sheaf of the diagonal immersion $X\hookrightarrow X\fp X$.

Here we consider an elliptic fibration $p\colon X\to B$ as above and the
associated ``dual'' fibration $\hat p\colon\X=X\to B$; we assume also that $X$ is smooth. We shall consider the 
relative \gif\ in this setting starting with the diagram \ref{e:ellipdiagram} and whose  kernel is the normalized relative
universal Poincar\'e sheaf $\Pc$ on the fibered product
$X\fp X$.

We then have a \gif
$$
\bS=\bS^{\Pc}\colon D(X_S)\to D(X_S) \,,\quad F\mapsto \bS(F)=
R\hat\pi_{S\ast}(\pi_S^\ast F\otimes\Pc_S)
$$
for every morphism $S\to B$.

Using our earlier invertibility  result  (Proposition \ref{p:invertibility}) or proceeding directly as in  \cite[Theorem 3.2]{BBHM98} (the latter was the first given proof), we easily obtain

\begin{prop} The \gif\ $\Phi$ is an equivalence of categories, (or a \FM). The inverse \FMF \ is  $\bS^{\Qc[1]}$, where
$$
\Qc= \Pc^\ast\otimes\pi^\ast p^\ast\omega^{-1}\,.
$$
\label{p:invert}
\qed\end{prop}

We shall denote by $\hat\bS$ the \FM\ $\bS^{\Qc}\colon D(X)\to D(X)$. The previous Proposition implies that if a sheaf $\Fc$ on $X$ is WIT$_i$ with respect to $\bS$ ($i=0,1$), then $\bS^i(\Fc)$ is WIT$_{1-i}$ with respect to $\hat\bS$ and $\hat\bS^{1-i}(\bS^i(\Fc))\simeq \Fc$. The analogous statement intertwining $\bS$ and $\hat\bS$ is also true.

\subsection{The spectral cover}
 
 We are going to see how the construction of vector bundles out of spectral data, first considered in \cite{Hit87} and \cite{BNR} can be easily described in the case of elliptic fibrations my means of the \FM\ we have defined. This construction was widely exploited by Friedman, Morgan and Witten \cite{FM, FMW97, FMW99} to construct stable bundles on elliptic Calabi-Yau threefolds. We shall come again to this point.

To start with, we think of things the other way round. We take a sheaf $\Fc$ of rank $n$ on an elliptic fibration $X\to B$ with certain properties on the fibers and construct its spectral data, namely, a pair $(C,\Lcl)$ where $C\hookrightarrow X$ is a closed subvariety projecting with finite fibers (generically of length $n$) onto $B$ (the \emph{spectral cover}) and a torsion free rank one sheaf $\Lcl$ on $C$ (in many case actually a line bundle), such that $\Fc$ can be recovered via the inverse \FM\ out of $(C,\Lcl)$.

We take $\Fc$ as a good parametrization of semistable sheaves of rank $n$ and degree 0 on the fibers of $X\to B$. Here good means that $\Fc$ is flat over the base $B$. The reason for doing so is twofold; first we know  the structure of the semistable sheaves of rank $n$ and degree 0 on a fiber as we report in Proposition \ref{p:ss}, second we can easily compute the \FM\ of torsion free rank one sheaves of degree cero on $X_s$ (Proposition \ref{p:cor1}).

Let us then fix a fiber $X_s$ of the elliptic fibration  ($s\in B$).  We denote by $\bS_s$ the \FM\ on the fiber with kernel $\Pc_s$.

The structure of the semistable sheaves of rank $n$ and degree 0 on $X_s$ is due to Atiyah 
\cite{At57} and Tu \cite{Tu94a} in the smooth case and to
Friedman-Morgan-Witten \cite{FMW99} for Weierstrass curves and locally free sheaves. The result we need is
\begin{prop} Every torsion-free semistable sheaf of rank  $n$ and degree
0  on $X_s$ is S-equivalent (see Appendix \ref{A:Simpson}) to a direct sum of torsion-free
 rank  $1$ and degree 0 sheaves:
$$
\Fc\sim \bigoplus_{i=0}^r (\Lcl_i\oplus\overset{n_i}\dots\oplus\Lcl_i)\,.
$$
\label{p:ss}
\qed\end{prop} If $X_s$ is smooth all the sheaves $\Lcl_i$ are line
bundles. If $X_s$ is singular, at most one of them, say  $\Lcl_0$, is
nonlocally-free; the number
$n_0$ of factors isomorphic to $\Lcl_0$ can be zero.

Now, let $\Lcl$ be a rank-one, zero-degree, torsion-free  sheaf
on $X_s$.

\begin{prop} $\Lcl$ is WIT$_1$ and $\bS^1_s(\Lcl)=\kappa(\xi^\ast)$, where
$\xi^\ast=[\Lcl^\ast]$ is the point of
$X_s\simeq \X_s$ defined by $\Lcl^\ast$. 
\qed\label{p:cor1}\end{prop}

%


%

We derive a few consequences of the two previous results. First is that  a zero-degree torsion-free sheaf of rank $n\ge 1$ and semistable on a fiber
$X_s$ has to be WIT$_1$ because this is what happens when the rank is 1.  We can state something stronger:
\begin{prop} Let $\Fc$ be a zero-degree sheaf of rank $n\ge 1$ on a fiber
$X_s$. Then
$\Fc$ is torsion-free and semistable on $X_s$ if and only if it is
WIT$_1$.
\label{p:sstrans0} 
\qed\end{prop}

A second consequence is that the unique \FM\ $\bS_s^1(\Fc)$ is supported by a finite number of points, again because Proposition \ref{p:cor1} tell us so in the rank one case. 

If we go back to our elliptic fibration  $p\colon X\to B$ and our sheaf $\Fc_s$ moves in a flat family $\Fc$ on $X\to B$, the support
of $\Phi_s^1(\Fc_s)$ moves as well giving a finite covering $C\to
B$. We notice, however, that the fiber over $s$ of the support of
$\bS^1(\Fc)$ may fail to be equal to the support of $\bS_s^1(\Fc_s)$. To
circumvent this problem we consider the closed subscheme defined by the
0-th Fitting ideal of $\bS\Fc$ (see for instance \cite{Rim} for a summary of  properties of the Fitting ideals). 
The precise definition (see \cite{FM,FMW97,FMW99,AsDo, HMP02}) is

\begin{defin} Let $\Fc$ be a sheaf on $X$. The \emph{spectral cover}
of $\Fc$ is the closed subscheme
$C(\Fc)$ of $\X$ defined by the 0-the Fitting ideal $F_0(\bS^1(\Fc))$ of
$\bS^1(\Fc)$.
\label{d;spectralcover}
\end{defin} 

The fibered structure of the spectral cover is a consequence of:
\begin{lemma} Let $\Fc$ be a zero-degree torsion-free semistable sheaf of
rank $n\ge 1$ on a fiber $X_s$.
\begin{enumerate}
\item The 0-th Fitting ideal $F_0(\what\Fc)$ of
$\what\Fc=\bS^1_s(\Fc)$ only depends on the S-equivalence class of
$\Fc$.
\item One has $ F_0(\what\Fc)=\prod_{i=0}^r{\mathfrak m}_i^{n_i}$, where
$\Fc\sim \bigoplus_{i=0}^r (\Lcl_i\oplus\overset{n_i}\dots\oplus\Lcl_i)$ is
the S-equivalence given by Proposition \ref{p:ss} and
$\mathfrak{m}_i$ is the ideal of the point $\xi_i^\ast\in\X_s=X_s$ defined by
$\Lcl_i^\ast$. Then, $\op{length}(\Oc_{\X_t}/F_0(\what\Fc))\ge n$ with
equality if either $n_0=0$ or
$n_0=1$, that is, if the only possible nonlocally-free rank 1
torsion-free sheaf of degree 0 occurs at most once.
\end{enumerate}
\label{l:sstrans}
\qed\end{lemma} 

Then we have the structure of the spectral cover for a relatively semistable sheaf of degree zero on fibers.

\begin{prop} If $\Fc$ is relatively torsion-free and semistable of rank
$n$ and degree zero on $X\to B$, then the spectral cover
$C({\Fc})\to B$ is a finite morphism with fibers of degree
$\ge n$. If in addition $\Fc$ is locally free, then all the fibers of the spectral cover $C({\Fc})\to B$  have degree $n$.
\label{p:naive2}
\end{prop}
\begin{proof} Since the  spectral cover commutes with base changes,
$C({\Fc})\to S$ is quasi-finite with fibers of degree $\ge n$ by
Lemma \ref{l:sstrans}; then it is finite. The second statement follows from (2) of the same Lemma.
\end{proof}

We can also give information about the spectral cover in some other cases. Take for instance a sheaf $\Fc$ on $X$ flat over $B$ and of
fiberwise  degree zero. We don't need to assume that $\Fc_s$ is semistable for every point $s\in B$. If this is true only for all the points $s$ of a dense open subset $U$ (i.e. any non-empty open subset if $B$ is irreducible), then as a consequence of    
Corollary \ref{c:wit1locus} and Proposition \ref{p:sstrans0} we have that $\Fc$ is still globally WIT$_1$ (that is, $\bS^0(\Fc)=0$) even if for $s\notin U$ we have  $\bS_s^0(\Fc_s)\neq 0$.  In this case, \emph{the spectral cover $C(\Fc)$ contains the whole fiber $X_s$}.

Let's go back to the case of a relatively torsion-free sheaf $\Fc$ semistable of rank
$n$ and degree zero on $X\to B$. We then have that the unique \FM\ of $\Fc$ is of the form
$$
\bS^1(\Fc)=i_\ast \Lcl
$$
where $i\colon C(\Fc)\hookrightarrow X$ is the immersion of the spectral cover and $\Lcl$ is a sheaf on $C(\Fc)$. What can be said about $\Lcl$? 

A first look at  Proposition \ref{p:cor1} seems to say that $\Lcl$ has rank one at every point. And this is actually what happens though one has to be careful because the spectral cover can be pretty singular. If $C(\Fc)$ is irreducible and reduced, then one can see quite easily that $\Lcl$ is torsion-free of rank one. When $C(\Fc)$ is reducible (it can even have multiple components),  torsion-freeness has to be substituted with another notion; people familiar with moduli problems won't be surprised to hear that the relevant notion is the one of \emph{pure sheaf} of maximal dimension introduced by Simpson \cite{Simp96a}.  We have described for the reader's convenience the definition of pure sheaf, polarized rank and Simpson stability in Appendix \ref{A:Simpson}. Using the definitions given there we have

\begin{prop} Let $m=\dim B$ and let $\Fc$ be a relatively torsion-free and semistable of rank
$n$ and degree zero on $X\to B$. Assume that all fibers of the spectral cover $C(\Fc)\to B$ have degree $n$ (this happens for instance if $\Fc$ is locally free). Then the restriction $\Lcl$ of the unique \FM\ $\bS^1(\Fc)$ is a pure sheaf of dimension $m$ and polarized rank one on the spectral cover $C(\Fc)$.

Conversely, given a closed subscheme $i\colon C(\Fc)\hookrightarrow$ such that $C(\Fc)\to B$ is a finite covering of degree $n$ and a sheaf $\Lcl$ on $C(\Fc)$ of pure dimension $m$ and polarized rank one, the sheaf $i_\ast(\Lcl)$ is WIT$_0$ with respect to $\hat\bS$ and the sheaf $\Fc=\hat\bS(i_\ast\Lcl)$ is a sheaf on $X\to B$ relatively torsion-free and semistable sheaf of rank $n$ and degree zero.
\label{p:spectraldata}
\end{prop} 

The most interesting case is when the base $B$ is a \emph{smooth curve}, that is, when $X$ is an elliptic surface.   Let 
then $\Fc$ be a sheaf on $X$ flat over $B$ and fiberwise
of degree zero. Assume that the restriction of $\Fc$ to the generic 
fiber is semistable, so that it is 
$\Fc$ is WIT$_1$ by Corollary \ref{c:wit1locus}. 
\begin{prop}  Let $V\subseteq B$ be 
the  relative semistability locus of $\Fc$.
\begin{enumerate}
 \item  The spectral cover
$C({\Fc})\to B$ is flat  of degree $n$ over $V$; 
then $C({\Fc_{V}})$ is a  Cartier divisor of $X_{V}$.
\item If $s\notin V$ is a point such that $\Fc_s$ is unstable, then
$C(\Fc)$ contains the  whole fiber $\X_{s}$.
\end{enumerate}
Thus $C(\Fc)\to B$ is finite (and automatically flat of degree $n$) if 
and only if $\Fc_s$ is semistable for every $s\in B$.
\label{p:naive3} 
\end{prop}

\section{Topological invariants of the \FM}\label{s:topinv}

Let $E$ be an object of $D(X)$.  We can compute the topological invariants of the
\gif\ 
$\bS(E)= R\hat\pi_{\ast}(\pi^\ast E\otimes\Pc)$ by using
the Riemann-Roch theorem for $\hat\pi$. There is a technical point here: even if we assume that $X$ and $B$ are smooth, $X\fp X$ may be not. However, Fulton established the so-called \emph{singular Riemann-Roch} and it turns out that our morphism $\hat\pi$ is what is called a l.c.i. morphism. By (\cite{Ful}, Cor.18.3.1), we have
\begin{equation}
\ch(\bS(E))=\hat\pi_{\ast}[\pi^\ast(\ch E)\cdot\ch(\Pc)\td(T_{X/B})]\,.
\label{e:grr}
\end{equation}
 The Todd class $\td(T_{X/B})$ is readily determined from Proposition
\ref{p:todd} and the Chern character of $\Pc$ is computed from (\ref{e:poin}). 

The Chern character  $\ch(\bS(E))$  has as many components $\ch_i (\bS(E))$ as the dimension of  $X$. We give here the precise expressions of $\ch_i (\bS(E))$ in two low dimensional cases, when $X$ is a \emph{smooth elliptic surface} \cite{BBH97a,HMP02} and when $X$ is an \emph{elliptic Calabi-Yau} \cite{ACHY01,ACHY04}.  

We can also compute the effect of the \gif\ $\bS$ on the \emph{relative Chern character} of $E$. This is specified by the relative rank $n$ and the relative degree $d$, that is, the  Chern character of the restriction of $E$ to a fiber $X_s$\footnote{This restriction is the derived inverse image $Lj_s^\ast E$ where $j_s\colon X_s\hookrightarrow X$ is the immersion of the fiber. When $E$ is a sheaf flat over $s$, then $Lj_s^\ast E=\rest E{X_s}$}. This is independent on the fiber, so we can apply Grothendieck-Riemann-Roch on a smooth curve to obtain that the relative Chern character of the Fourier-Mukai
transform
$\bS_S(E)$ is $(d,-n)$, that is:
\begin{equation}
(\rk(\bS_S(E)), d(\bS_S(E))=(d(E), - \rk(E))
\label{e:heckechern}
\end{equation}

If we denote by $\mu_{rel}(E)=d/n$ the relative slope, we have

\begin{prop} If $\Fc$ is a WIT$_i$ sheaf on $X$ and $d\ne 0$, then $\mu_{rel}(\widehat\Fc)=-1/\mu(\Fc)$. Moreover
\begin{enumerate}
\item If $\Fc$ is WIT$_0$, then $d(\Fc)\ge 0$, and
$d(\Fc)=0$ if and only if $\Fc=0$.
\item If $\Fc$ is WIT$_1$, then $d(\Fc)\le 0$. 
\end{enumerate}
\label{p:relwit}
\end{prop}

\subsection{The case of elliptic surfaces}

Let us denote by $e$  the degree of the divisor $\bar K$ on $B$; we have
$\Theta\cdot p^\ast \bar K=e=-\Theta^2$ and $K_{X/B}=p^\ast \bar K\equiv e\,{\mathfrak f}$. 

By Proposition \ref{p:todd} the Todd class of the virtual relative tangent
bundle of $p$ is given by
\begin{equation}
\td(T_{X/B})=1-\tfrac12\, p^{-1}\bar K+ e\,w\,,
\label{eq:todd2}
\end{equation} where $w$ is the fundamental class of
$X$. 

Now, if $E$ is an object of $D(X)$, the Chern character of the
\FM\  $\bS(E)$ is given by
\begin{align*}
\ch(\bS(E))=&\pi_\ast[\pi^\ast(\ch E)\cdot(1-\delta_\ast(1)-\tfrac12\,
\delta_\ast(p^\ast \bar K)+e\,\delta_\ast(w))\cdot(1+\pi^\ast H-\tfrac12
e\,w)\\ &\cdot (1-\tfrac12 p^\ast \bar K+ew)]\cdot (1+\Theta-\tfrac12
w)\cdot(1+e\,{\mathfrak f})\,.
\end{align*}

Thus, the Chern characters of $\bS( E)$ are
\begin{equation}
\begin{aligned}
\ch _0(\bS( E))&=d\\
\ch _1(\bS( E))&=- c_1( E)+d\, p^\ast
\bar K+(d-n)\Theta+(c-\tfrac12\,ed+s)\,\mathfrak f\\
\ch _2(\bS( E))&=(-c-de+\tfrac12\,ne) w
\end{aligned}
\label{e:chern2}
\end{equation} 
where $n=\ch_0( E)$, $d=c_1( E)\cdot {\mathfrak f}$ is the relative
degree, $c=c_1( E)\cdot H$ and $\ch _2( E)=s\,w$.

Similar calculations can be done for the inverse \FM\ giving rise to the formulae
\begin{equation}
\begin{aligned}
\ch_0(\widehat{\bS}( E))&= d\\
\ch_1(\widehat{\bS}( E))&=(c_1( E))- n p^\ast \bar K-(d+ n)H+(s+ n e- c-\tfrac12\,e d){\mathfrak f}\\
\ch_2(\widehat{\bS}( E))&=-(c+ de+\tfrac12\,ne)w
\end{aligned}
\label{e:chern3}
\end{equation}

\subsection{The case of elliptic Calabi-Yau threefolds}

When $X$ is an elliptic Calabi-Yau threefold and $B$ is a smooth surface (with the restrictions mentioned in section \ref{s:elliptic}), the formula given in Proposition \ref{p:todd}  for the Todd class of the relative tangent bundle takes the form
\begin{equation}
\td (T_{X/B})=1- \frac 12 c_1+ \frac 1{12}(13c_1^2+12\Theta
c_1)- \frac 12 \Theta c_1^2
\label{eq:todd22}
\end{equation}
with  $c_1= p^\ast c_1(B)=- p^\ast(K_B)$ and equation (\ref{e:poin}) is 
\begin{equation}
\Pc=\Ic_\Delta\otimes \pi^\ast\Oc_X(\Theta)\otimes \hat\pi ^\ast\Oc_X(\Theta)\otimes q^\ast \omega_B^{-1}
\label{e:poin2}
\end{equation}
where $\omega_B=\Oc_B(K_B)$ is the canonical line bundle of $B$ and $\delta$ is the diagonal immersion.  Note first that $\ch(\Ic_\Delta)=1-\ch(\delta_\ast \Oc_X)$. Singular Riemann-Roch gives
$$
\ch(\delta_\ast \Oc_X)\op\td (X\fp X)= \delta_\ast (\ch(\Oc_X) \op\td (X))
$$
where one has the expressions for $\op\td (X)$ (see \eqref{eq:toddx}) and $\op\td (X\times_B X)$ given by
\begin{align*}
\op\td (X)&= 1+ \frac 1{12}(c_2+11c_1^2+12\Theta c_1)\\
\op\td (X\fp X)&=\hat\pi^\ast\op\td (X)\pi_1^\ast\op\td (T_{X/B})
\end{align*}
with $c_2=p^\ast(c_2(B))$. The Chern character of the ideal sheaf $\Ic_\Delta$ is then given by
(with the diagonal class $\Delta=\delta_\ast(1)$)
$$
\ch({\Ic_\Delta})=1-\Delta- \frac 12 \Delta\cdot \hat\pi^\ast c_1+\Delta\cdot \hat\pi^\ast (\Theta\cdot c_1)+\frac 56\Delta\cdot \hat\pi^\ast (c_1^2)+ \frac 12\Delta\cdot \hat\pi^\ast(\Theta c_1^2)
$$
and one can compute $\ch(\Pc)$ form that expression and find the general formula for $\bS(E)$ form (\ref{e:grr}). 

We shall consider for simplicity objects $E$ in $D(X)$  whith Chern characters given by 
 \begin{equation}
 \begin{aligned}
\ch_0(E)&=n_E \\
\ch_1(E)&=x_E\Theta+ p^\ast S_E \\
\ch_2(E)& =\Theta p^\ast \eta_E+a_E \mathfrak f \\
\ch_3(E)& =s_E
\end{aligned}
\label{e:chE}
\end{equation}
where $\eta_E, S_E\in A^1(B)\otimes_\Z\bbQ$, $s_e\in A^3(X)\otimes_\Z\bbQ \simeq \bbQ$) anf $\mathfrak f\in A^2(X)\otimes_\Z\bbQ$ is the class of a fiber of $p$. This cover the most important applications.
Now (\ref{e:grr}) and the corresponding formula for the inverse \FM\ $\hat\Phi$ provide the first Chern characters of $\bS(E)$ and $\hat\bS(E)$. They are
\begin{equation}
\begin{aligned}
\ch_0(\bS(E))&=x_E\\
\ch_1(\bS(E))&=-n_E\Theta+ p^\ast\eta_E- \frac 12 x_E c_1\\
\ch_2(\bS(E))&=(\frac 12 n_E c_1-p^\ast S_E)\Theta+(s_E -\frac 12 p^\ast \eta_E c_1\Theta + \frac 1{12} x_E c_1^2\Theta)\mathfrak f\\
\ch_3(\bS(E))&=- \frac 16  n_E \Theta c_1^2-a_E + \frac 12 \Theta c_1
p^\ast S_E
\end{aligned}
\label{e:chVS}
\end{equation}
and
\begin{equation}
\begin{aligned}
\ch_0(\hat\bS(E))&=x_E \\
\ch_1(\hat\bS(E))&=-n_E \Theta+ p^\ast \eta_E+ \frac 12 x_E c_1\\
\ch_2(\hat\bS(E))&=(- \frac12 n_E c_1-p^\ast S_E)\Theta+(s_E + \frac 12 p^\ast \eta_E c_1\Theta
+ \frac 1{12}x_E c_1^2\Theta)\mathfrak f\\
\ch_3(\hat\bS(E))&=- \frac 16  n_E \Theta c_1^2-a_E- \frac 12 \Theta c_1
p^\ast S_E+x_E\Theta c_1^2
\end{aligned}
\label{e:chVSS}
\end{equation}

Many nice geometrical properties and applications to Physics can be derived from \FMF\ theory, and in particular from the formulae above.

\section{Applications  to moduli problems}\label{s:applmod}
In this section we report some applications of the \gif\ theory to moduli problems. The material is taken mostly from \cite{HMP02, BBHM02}.
\subsection{Moduli of relatively semistable sheaves on elliptic fibrations} 

The structure of relatively semistable
sheaves on an elliptic fibration $p\colon X\to B$ can be described in terms of the relative \FM\ introduced in section \ref{s:elliptic}. The techniques and results are of quite different nature when the degree on fibers of the sheaves are zero or not, so we consider separately the two cases. 

\subsubsection*{The case of relative degree equal to zero}

If we start with a 
single fiber $X_{s}$, then Proposition
\ref{p:ss} means that S-equivalence classes of semistable sheaves of rank $n$ and
degree 0 on $X_s$ are equivalent to families of $n$  torsion-free rank
one sheaves of degree zero, $\Fc_s\sim \bigoplus_{i=0}^r
(\Lcl_i\oplus\overset{n_i}\dots\oplus\Lcl_i)$. This gives a one-to-one 
correspondence 
\begin{equation}
\begin{aligned}
\overline{\M}(X_{s},n,0)&\leftrightarrow \op{Sym}^n X_s\\
\Fc_s&\mapsto n_0 x^\ast _0+\dots+n_r x^\ast_r\,,
\end{aligned}  \label{eq:sym}
\end{equation}
(where $x_i^\ast$ is the point of $X_s$ that corresponds to $[\Lcl^\ast_i]$ under $\varpi\colon X_s\iso \X_s$) 
 between the moduli space of  torsion-free and semistable sheaves of rank $n$
and degree 0 on $X_{s}$ and the $n$-th symmetric product of the compactified 
Jacobian $\X_s\simeq X_s$. The reason for taking duals comes from Proposition \ref{p:cor1} and
Lemma \ref{l:sstrans}: the skyscraper sheaf $\kappa([\xi^\ast_i])$ is
the \FM\  of $\Lcl_i$, and if $n_0=0$, then
$n_1 x^\ast _1+\dots+n_r x^\ast_r$ is the spectral cover $C(\Fc_s)$.

As we have seen, when $\Fc_s$ moves in a flat family, $C(\Fc_s)$ moves to give the spectral cover of $\Fc$ defined in \ref{d;spectralcover}.  We can then  extend (\ref{eq:sym}) to the whole elliptic
fibration $X\to B$ under some technical assumptions (namely, that $B$ is
\emph{normal of dimension bigger than zero} and the generic fiber is
 \emph{smooth}).  

There are  two different varieties that parametrize flat families of clusters of points on the fibers. The first one is the \emph{Hilbert scheme} $\op{Hilb}^n(X/ B)\to B$ of $B$-flat subschemes of $\X$ of fiberwise dimension 0 and length $n$. The second one is the \emph{relative symmetric
$n$-product} $\op{Sym}^n_{ B} X$ of  $X\to B$.  They are not isomorphic in general, only birational; actually, there is a Chow morphism
$\op{Hilb}^n(X/ B)\to \op{Sym}^n_{ B}X$ mapping a cluster of length $n$ to the $n$ points defined by the culster, which induces an isomorphism
$\op{Hilb}^n( X'/ B)\simeq
\op{Sym}^n_{ B} X'$, where $ X'\to B$ is the smooth locus of $p\colon X\to B$. 

 Let us denote by $\overline{\M}(n,0)$ the (coarse) moduli scheme of
torsion-free and semistable sheaves of rank $n$ and degree 0 on the
fibers of $X\to B$, by  $\M(n,0)$
 the open subscheme of $\overline{\M}(n,0)$ defined by those
sheaves on fibers which are 
$S$-equivalent to a direct sum of line bundles (see \cite{Simp96a} or Appendix \ref{A:Simpson}).
  
If $\Fc$ is a sheaf on
$X\to B$ defining a $B$-valued point of ${\bM}(n,0)$, then the spectral 
cover $C(\Fc)$ is flat of degree $n$ over $B$, 
and then defines a $B$-valued  point of
$\op{Hilb}^n( X'/ B)$ which depends only on the S-equivalence class of
$\Fc$.  One easily deduce that we can define in that way a morphism of
$B$-schemes
\begin{align*} {\bC}'\colon\M(n,0)& \to \op{Hilb}^n( X'/ B)\simeq \op{Sym}^n_{ B} X' \\
[\Fc] & \mapsto C(\Fc)
\end{align*}
where $[\Fc]$ is the point of $\overline{\M}(n,0)$ defined by $\Fc$.

We now have \cite[Theorem 2.1]{HMP02}
\begin{thm}  $ {\bC}'$ is an isomorphism $ {\bC}'\colon {\M}(n,0)\iso \op{Hilb}^n( X'/
B)\simeq \op{Sym}^n_{ B} X'$ and extends to an isomorphism of
$B$-schemes
${\bC}\colon \overline{\M}(n,0)\iso \op{Sym}^n_{ B}\X$.   For every
geometric point $\Fc\sim \bigoplus_i
(\Lcl_i\oplus\overset{n_i}\dots\oplus\Lcl_i)$ the image ${\bC}([\Fc])$ is
the point of $\op{Sym}^n_{ B}\X$ defined by
$n_1\xi_1^\ast+\dots+n_r\xi_r^\ast$.
\label{th:sym}
\qed\end{thm}

 We  denote by $J^n\to B$ the relative Jacobian of line bundles on
$p\colon X\to B$ fiberwise of degree $n$. We have an isomorphism  $\tau\colon J^n\iso J^0$ which is the translation $\tau(\Lcl)=\Lcl\otimes\Oc_{X}(-n\Theta)$,  and also the the natural
involution $\iota\colon J^0\iso J^0$ mapping a line bundle to its inverse. Let
$\gamma\colon J^n\iso J^0$ be the composition
$\gamma=\iota\circ\tau$, so that if $x_1+\dots+x_n$ is a
positive divisor in $X'_s$, then $\gamma[\Oc_{X_s}(x_1+\dots+x_n)]=
[\Lcl_1^\ast\otimes\dots\otimes\Lcl_n^\ast]$, where $[\Lcl_i]$ corresponds to $x$ under $\varpi\colon X\to\X$. 
\begin{thm} There is a commutative diagram of $B$-schemes
$$
\xymatrix{\M(n,0)\ar[r]^{\overset{\bC}\sim} \ar[d]_{\det}&
\op{Sym}^n_B(\X')
\ar[d]^{\phi_n}\\ J^0 &\what J^n \ar[l]_{\overset\gamma\sim}}
$$
 where $\det$ is the ``determinant'' morphism and $\phi_n$ is the Abel
morphism of degree $n$. 
\label{th:main}
\qed\end{thm}

The previous Theorem generalizes \cite[Theorem 3.14]{Fr95} and can be
considered as a global version of the results obtained in Section 4 of
\cite{FMW99} about the relative moduli space of locally free sheaves on
$X\to B$ whose restrictions to the fibers have rank $n$ and trivial
determinant. Theorem
\ref{th:main} leads to these results by using standard structure 
theorems for the Abel morphism.

Let $\Lcl_n$ be a universal line bundle over $q\colon
\X\fp  J^n \to J^n $. The Picard sheaf ${\Pc}_n=R^1
q_\ast(\Lcl_n^{-1}\otimes\omega_{\X/B})$ is a locally free sheaf of rank
$n$ and then defines a projective bundle $\mathbb
P({\Pc}_n^\ast)=\op{Proj}S^\bullet (\Pc_n)$.  We now have a diagram 
$$
\xymatrix{\M_U(n,0)\ar[r]^{\sim\;}\ar@{^{(}->}[d]&\op{Sym}^n_{U} X_{U}
\ar[rr]^{\sim}\ar@{^{(}->}[d]&&\mathbb
P(\rest{{\Pc_n^\ast}}{U})\ar@{^{(}->}[d] \\
\M(n,0)\ar[r]^{\overset{\mathbf{C}}\sim\;\;} \ar[d]_{\det}&
\op{Sym}^n_B(X')
\ar[d]_{\text{\tiny Abel}}\ar@{^{(}->}[rr]^{\text{\tiny \ \ dense}}& &
\mathbb P(\Pc_n^\ast)
\ar[lld]\\ J^0 & J^n
\ar[l]_{\overset\gamma\sim}& &}
$$
where $U\hookrightarrow B$ is the open subset supporting the smooth fibers of
$p\colon X\to B$ and $X_U=p^{-1}(U)$.  The immersions of the symmetric products into the projective bundles follow from the structure of the Abel morphism (cf.  \cite{AK80}).

\begin{corol}
$\widetilde\Pc_n\iso(\det)_\ast\Oc_{\M(n,0)}(\Theta_{n,0})$.
\label{c:moduliproj2}\qed
\end{corol}

Since $e^\ast(\widetilde\Pc_n)\iso (p_\ast
\Oc_X(nH))^\ast$,  we obtain the structure theorem proved in
\cite{FMW99}: Let $\M(n,\Oc_X)=(\det)^{-1}(\hat e(B))$ be the
subscheme of
 those locally free sheaves in $\M(n,0)$ with trivial determinant and
$\M_{\U}(n,\Oc_X)=\M(n,\Oc_X)\cap \M_{\U}(n,0)$.
\begin{corol}  There is a dense immersion of $B$-schemes
$\M(n,\Oc_X)\hookrightarrow \mathbb P({\V}_n)$, where
${\V}_n=p_\ast(\Oc_X(nH))$. Moreover, this morphism induces an
isomorphism of $\U$-schemes
$\M_{\U}(n,\Oc_X)\iso \mathbb P(\rest{{{\V}_n}}{\U})$.
\end{corol}

\subsubsection*{The case of  nonzero relative degree}

The study of relatively stable sheaves of positive degree on elliptic fibrations was done in \cite{BBHM02,Ca1}.  The first thing to do is to characterize WIT$_0$ sheaves.  We need a preliminary result (cf. also \cite{Bri98}), whose proof is given to show the techniques employed.

\begin{lemma} A coherent sheaf $\Fc$ on $X$ is WIT$_0$ if and only if 
$$
\Hom_X(\Fc, \Pc_\xi)=0
$$
for every $\xi\in \X\simeq X$.
\label{l:wit}
\end{lemma} 
\begin{proof} By Proposition (\ref{p:cor1}), $\Pc_\xi$ is WIT$_1$ and $\bS^1(\Pc_\xi)=\kappa(\xi^\ast)$, where $\xi^\ast$ is the point of $\widehat X_t$ corresponding to $\Pc_\xi^\ast$ ($t=p(x)$). Then, Parseval formula (Proposition \ref{p:parseval} implies that
$$
\Hom_X(\Fc, \Pc_\xi)\simeq \Hom_{D(\X)}(\bS(\Fc), \kappa(\xi^\ast)[1])\,.
$$
If $\Fc$ is not WIT$_0$, there is a point $\xi^\ast\in\widehat X$ such that a nonzero morphism $\bS^1(\Fc)\to \kappa(\xi^\ast)$ exists. This gives rise to a non-zero morphism $\bS(\Fc) \to \kappa(\xi^\ast)[1]$ in the derived category, so that $\Hom_X(\Fc, \Pc_\xi)\neq 0$. The converse is straightforward.
\end{proof}

\begin{prop} Let $\Fc$ be a relatively (semi)stable sheaf on $X$,
with
$d(\Fc)>0$. Then $\Fc$ is WIT$_0$ and its Fourier-Mukai transform $\widehat\Fc$ is
relatively (semi)stable.\label{p:ss2}
\qed\end{prop}

We don't give a complete proof (it can be founded in \cite{BBHM02}). The idea is to use the previous Lemma to show that the sheaf is WIT$_0$ and to apply the invertibility of the Fourier-Mukai transform to get a contradiction from the assumption that $\what\Fc$ could be destabilized.

\begin{corol} Let $\Fc$ be a torsion-free semistable sheaf on $X_t$ of degree $d>0$. Then
$H^1(X_t,\Fc\otimes\Pc_\xi)=0$ for every $\xi\in{\widehat X}_t$.
\qed\end{corol}

Proposition \ref{p:ss2} also gives the characterization of semi(stable) sheaves of relative negative degree. They are WIT$_1$ and their \FM s are also (semi)stable.

\medskip
As a side result of these results, we see that the relative Fourier-Mukai transform provides a
characterization of some moduli spaces of relatively stable bundles. 
Let $J_n\to
B$ be the relative Jacobian of invertible sheaves of relative degree 
$n$ and $\bar
J_n\to B$ the natural compactification of $J_n$ obtained by adding to 
$J_n$ the $
B$-flat coherent sheaves on $p\colon X\to B$ whose restrictions to 
the fibers of
are torsion-free, of rank one and degree
$n$ \cite{AK80}. Proposition \ref{p:ss2} gives (cf. \cite{BBHM02, Ca1}):
\begin{thm} Let $\Nc$ be an invertible sheaf on $X$ of relative 
degree $m$,
and let
$\mathcal M(n,$ $nm-1)$ be the moduli space of  rank $n$ relatively 
$\mu$-stable
sheaves on
$X\to B$ of degree $nm-1$. The Fourier-Mukai  transform induces an 
isomorphism of
$B$-schemes
\begin{eqnarray*}
\bar J_n& \stackrel{\bS^0\otimes\Nc}{\relbar\joinrel\longrightarrow}& 
\mathcal M
(n,nm-1)\\
\Lcl&\mapsto& \bS^0(\Lcl)\otimes\Nc\,.
\end{eqnarray*}
\label{t:cor3.5}
\qed\end{thm}

\subsection{Absolutely semistable sheaves on an elliptic surface}

In this section we  apply the theory so far developed to the study of
the moduli space of absolutely stable sheaves on an elliptic surface.

We relay on the computation of the Chern character of the
\FM s provided by \eqref{e:chern2}. This enable us to the study of the
preservation of stability. We shall see that stable sheaves on spectral
covers transform to absolutely stable sheaves on the surface and prove
that in this way one obtains an open subset of the moduli space of
absolutely stable sheaves on the surface.

In the whole section the base $B$ is a \emph{projective smooth curve}.

The elliptic surface is polarized by $H=a\Theta + b \mathfrak f$ for suitable positive integers $a$ and $b$.
Any effective divisor $i_t\colon C\hookrightarrow X$ is then polarized by the restricition $H_C=H\cdot C$ of $H$. Then, even if the curve $C$ is not integral, we can define rank, degree and Simpon stability for pure sheaves concentrated on $C$ as discused in Appendix \ref{A:Simpson}.

Moreover, a pure sheaf $\Qc$ of dimension 1, with support contained in $C$ is Simpson (semi)stable with respecto to $H_C$ if and only if $i_\ast\Qc$ is is Simpson (semi)stable with respecto to $H$

\subsubsection*{Preservation of absolute stability}

Let $\Fc$ be a sheaf on $X$ flat over $B$ with Chern
character $(n,\Delta,s)$ with $n>1$. Assume  that the restrictions of $\Fc$ to the fibers $X_s$ are semistable of degree 0. Then  $\Fc$ is WIT$_1$  
and the spectral cover $C(\Fc)$ is a Cartier divisor finite of degree $n$ over $B$ (Proposition \ref{p:naive3}). Moreover $\what\Fc$ is a sheaf of pure dimension 1 whose support is contained in $C(\Fc)$. 

\begin{prop} For any integer $a>0$ there is $b_0>0$ depending only on the numerical invariants $(n,c,s)$, such that for any $b>0$ the following is true: $\Fc$ is (semi)stable on $X$ with respect to $H=a\Theta+b\mathfrak f$ if and only if $\what\Fc$ is (semi)stable on $X$ with respect to $H=a\Theta+b\mathfrak f$, and then, if and only if $\Lcl=\rest{\what\Fc}{C(\Fc)}$ is (semi)stable as pure dimension 1 sheaf on the spectral cover $C(\Fc)$. 
\label{p:mainprop}
\end{prop}
\begin{proof} By \eqref{e:chern2}, one has
\begin{equation} 
\chi(X, \what\Fc(mH))= (nb-nae-as)m+ c-ne+\frac12 n c_1
\label{eq:hilbtransform}
\end{equation}
where $n=\ch_0(\Fc)$,  $c=c_1(\Fc)\cdot \Theta$, $\ch _2(\Fc)=s\,w$, $e=-\Theta^2$ and $c_1=c_1(B)$. The Simpson slope of $\what\Fc$ is 
$$
\mu(\what\Fc)= \frac{ c-ne+\frac12 n c_1}{nb-nae-as}\,.
$$
Let now 
\begin{equation} 0\to\G\to \what\Fc\to K\to 0
\label{eq:des}
\end{equation} be an exact sequence. Then $\G$ is concentrated on $C(\Fc)$, so that it is WIT$_0$ and the \FM\ $\bar\Fc=\what{\G}$ has relative degree 0 and it is WIT$_1$ by Proposition \ref{p:spectraldata}.  Reasoning as above, the Simpson slope of $\G$ is
$$
\mu(\G)= \frac{\bar c-\bar n e+\frac12 \bar n c_1}{\bar n b-\bar n ae-a\bar s}\,,
$$
where bars denote the topological invariants of $\bar\Fc$.
Moreover one has the exact sequence
$$
0\to \bar\Fc \to \Fc \to \what K \to 0\,.
$$
Assume that $\Fc$ is semistable with respect to $H$ and that  \eqref{eq:des} is a destabilizing sequence. Then $\mu(\G)>\mu(\what\Fc)$, 
which is equivalent to 
$$
(\bar n c - n \bar c) b + a(n\bar c-\bar n c +c\bar s-\bar c s+ e(\bar n s-n\bar s)+\tfrac12 c_1(n\bar s-\bar n s))>0
$$
Since the family of subsheaves of $\Fc$ is bounded,
there is a finite number of possibilities for the Hilbert polynomial of
$\bar\Fc$. Then, there is a finite number of possibilities for $n\bar c-\bar n c +c\bar s-\bar c s+ e(\bar n s-n\bar s)+\tfrac12 c_1(n\bar s-\bar n s)$ so that for fixed $a>0$ and $b\gg0$ the destabilizing condition is
$$
n\bar c-\bar n c>0\,.
$$
On the other hand, the semistability of $\Fc$ implies that 
$$
\frac{c_1(\bar\Fc)\cdot (a\Theta+b\mathfrak f)}{\bar n} \le \frac{c_1(\Fc)\cdot (a\Theta+b\mathfrak f)}{n}
$$
that is, $n\bar c-\bar n c\le 0$, which is a contradiction. The corresponding
semistability statement is proven analogously.

For the converse, assume that $\what\Fc$ is semistable on $X$ with respect to $H=a\Theta+b\mathfrak b$ for $b\gg 0$ and that 
$$
0\to \bar\Fc \to \Fc \to \Qc \to 0
$$
is a destabilizing sequence. We can assume that $\bar n< n$ and that $\Qc$ is torsion free and $H$-semistable; moreover one has $n c_1(\bar\Fc)\cdot H> \bar n c_1(\Fc)\cdot H$, that is, $n (a \bar c+ b \bar d)> \bar n ac$. 

The sheaf $\bar\Fc$ is WIT$_1$ so that $\bar d\le 0$ by Proposition \ref{p:relwit}. 
Assume first that $\bar d<0$ and let $\rho$ be the maximum of the integers
$nc_1(\tilde \Fc)\cdot \Theta-\rk(\tilde\Fc)c$ for all nonzero subsheaves $\bar\Fc$ of
$\Fc$. Then $n (a \bar c+ b \bar d)- \bar n ac\le a\rho +nb\bar d$ is
strictly negative for $b$ sufficiently large, which is absurd.  It follows that  $\bar d=0$ and the destabilizing condition is  
$$
n\bar c-\bar n c>0\,.
$$
Moreover $d(\Qc)=0$. Since $\Qc$ is torsion free,  for every $s\in B$ there is an exact sequence 
$$
0\to \bar\Fc_s \to \Fc_s \to \Qc_s \to 0
$$
so that $\Qc_s$ is semistable of degree 0. Then $\Qc$ is WIIT$_1$ and one has an exact sequence of \FM s:
$$
0\to \what{\bar\Fc}\to \what\Fc \to \what\Qc\to 0\,.
$$
Proceeding as above we see that the semistability of $\Fc$ for $b\gg0$ implies that $n\bar c-\bar n c\le0$, which is a contradiction.
\end{proof} 

Then absolute stability with respect to $a\Theta+b\mathfrak b$ is preserved for $b\gg0$ depending on $a$ and on the Chern character $(n,\Delta,s)$.  This was proved in a different way in \cite{HMP02}, similar results can be founded in \cite{JM,Yo01}.

Notice that \eqref{eq:hilbtransform} is deduced from the formula
$$
\chi(C(\Fc),\Lcl(mH))= \chi(X,\what\Fc(mH))= (C(\Fc)\cdot H) \,m+ \frac12 C(\Fc)\cdot c_1(X)+ \ch_2(i_\ast\Lcl)\\
$$
where $i\colon C(\Fc)\hookrightarrow X$ is the immersion, so that the polarized rank of $\Lcl$ is 1  in agreement with Proposition \ref{p:mainprop}, and its 
Euler characteristic is
$$
\chi(C(\Fc),\Lcl)=\frac12 C(\Fc)\cdot c_1(X)+ \ch_2(i_\ast\Lcl)=\frac12 n(c_1-e) + \ch_2(i_\ast\Lcl)
$$
We then see that given a Cartier divisor $i\colon C\hookrightarrow X$ flat of degree $n$ over $B$ and a pure dimension one sheaf $\Lcl$ of $C$ of polarized rank 1 and  we write $\ell=C\cdot\Theta$ and $r=\chi(C,\Lcl)$, then the numerical invariants of $\Fc=\widehat{\bS}(i_\ast\Lcl)$ are given according to \eqref{e:chern3} by
\begin{equation}
\rk(\Fc)=n\,,\quad d=0\,\quad c=ne+r-\tfrac12 nc_1\,,\quad s=\ell-ne
\label{e:chern4}
\end{equation}
 
\subsubsection*{Moduli of absolutely stable sheaves and compactified
Jacobian of the universal spectral cover}
In this subsection we shall prove that there exists a universal spectral
cover over a Hilbert scheme and that the \gif\  embeds
the compactified Jacobian of the universal spectral cover as an open 
subspace the
moduli space of absolutely stable sheaves on the elliptic surface (cf. \cite{HMP02}).
Most of what is needed has been proven
in the preceding subsection. 

We start by describing the spectral cover 
of a relatively semistable sheaf in terms of the isomorphism
 $\overline{\M}(n,0)\iso \op{Sym}^n_B \X$ provided by Theorem \ref{th:sym}.
There is a ``universal'' subscheme
$$ C\hookrightarrow \X\fp \op{Sym}^n_B \X
$$ defined as the image of the closed immersion $\X\fp
  \op{Sym}^{n-1}_B
  \X\hookrightarrow\X\fp \op{Sym}^n_B \X$,
  $(\xi,\xi_1+\dots+\xi_{n-1})\mapsto (\xi,\xi+\xi_1+\dots+\xi_{n-1})$.
  The natural morphism $g\colon C\to \op{Sym}^n_B \X$ is finite and 
  generically of degree $n$. Let $A\colon S\to \op{Sym}^n_B \X$ be a 
  morphism of $B$-schemes and let $C(A)=(1\times 
 A)^{-1}(C)\hookrightarrow
 \X_S$ be the closed subscheme of $\X_S$ obtained by pulling the
 universal subscheme back by the graph $1\times A\colon\X_S\hookrightarrow
  \X\fp \op{Sym}^n_B \X$  of $A$. There is a finite morphism 
  $g_A\colon C(A)\to S$ induced by $g$. 
  
  By Theorem \ref{th:sym}, a $S$-flat sheaf $\Fc$ on
  $X_S$ fiberwise torsion-free and semistable of rank $n$ and
degree 0 defines a morphism $A\colon S\to \op{Sym}^n_S(\X_S)$; we easily see
from Lemma \ref{l:sstrans} that
 \begin{prop} $C(A)$ is the spectral cover associated to $\Fc$, $C(A)=C({\Fc})$.
 \qed 
 \end{prop} 
  
 When $S=B$, $A$ is merely a section of $\op{Sym}^n_B
  \X\simeq\overline{\M}(n,0)\to B$. In this case, $C(A)\to B$ is flat 
  of degree $n$  because it is finite and $B$ is a smooth curve. 
 The 
  same happens when the base scheme is of the form  $S=B\times T$, where
  $T$ is an arbitrary scheme:
  \begin{prop} For every morphism $A\colon B\times T\to
  \op{Sym}^n_B \X$ of $B$-schemes, the spectral cover projection
$g_{A}\colon C(A)\to B\times T$ is flat of degree $n$.
  \label{p:relSpecFlat} 
  \qed\end{prop}

 If the section $A$ takes values in $\op{Sym}^n_B \X'\simeq{\M}(n,0)\to B$, then
 $g_A\colon C(A)\to B$ coincides with the spectral cover constructed in 
  \cite{FMW99}.

Let now
$\Hc$ be the Hilbert scheme of sections of the projection
$\hat\pi_n\colon \op{Sym}^n_B \X\to B$. If
$T$ is a 
$k$-scheme, a
$T$-valued point of $\Hc$ is a section $B\times T\hookrightarrow
\op{Sym}^n_B
\X\times T$ of the projection $\hat\pi_n\times 1\colon
\op{Sym}^n_B\X\times T\to B\times T$, that is, a morphism $B\times
T\to \op{Sym}^n_B\X$ of $B$-schemes. There is a universal section
$\Ac\colon B\times\Hc\to \op{Sym}^n_B\X$. It gives rise to a
\emph{``universal'' spectral cover} $\CA\hookrightarrow\X\times\Hc$.
By Proposition \ref{p:relSpecFlat}, the ``universal'' spectral cover
projection
$g_{\Ac}\colon \CA\to B\times\Hc$ is flat of degree $n$.
It is endowed with a relative polarization
$\Xi=H\times\Hc)$ where $H=a\Theta+b\mathfrak f$ for $a>0$ and $b\gg 0$. 

Let $\bar\bJ^r \to\Hc$ be the  functor of  sheaves of  pure
dimension one, polarized rank one, Euler characteristics $r$ and
semistable with respect to $\Xi$ on the fibers of the flat family of
curves
$\rho\colon\CA\to\Hc$.  

Let $\Hc_{\ell}$ be the subscheme of those points $h\in \Hc$ such that
$\rho^{-1}(h)\cdot\Theta=\ell$. The subscheme
$\Hc_{\ell}$ is a disjoint union of connected components of $\Hc$ and
then we can decompose $\rho$ as a union of projections
$\rho_{\ell}\colon\CA_{\ell}\to\Hc_{\ell}$. We
decompose $\bar\bJ^r$ accordingly into functors $\bar\bJ_{\ell}^r$. 

By Theorem 1.21 of \cite{Simp96a} (cf. Theorem \ref{t:fundamental}) there exists a coarse moduli scheme 
$\bar\Jc_{\ell}^r$ for $\bar\bJ_{\ell}^r$. It is projective over
$\Hc_{\ell}$ and can be  considered as a ``compactified'' relative
Jacobian of the universal spectral cover 
$\rho_{\ell}\colon\CA_{\ell}\to\Hc_{\ell}$. The open subfunctor
$\bJ_{\ell}^r$  of
$\bar\bJ_{\ell}^r$ corresponding to stable sheaves has a fine
moduli space $\Jc_{\ell}^r$ and it is an open subscheme of
$\bar\Jc_{\ell}^r$.

On the other side we can consider the coarse moduli scheme $\overline\M(a,b)$
torsion-free sheaves on $X$ that are semistable with respect to $a H+b\mu$ 
and have Chern character $(n,\Delta,s\,w)$ where $d=\Delta\cdot \mathfrak f=0$, and the values of $c=\Delta\cdot\Theta$ and $s$ are giving by \eqref{e:chern4}. We also have the corresponding moduli functor $\overline{\bM}(a,b)$ (see again 
\cite{Simp96a}). Let $\M(a,b)\subset \overline\M(a,b)$ the open subscheme 
defined by the stable sheaves. It is a fine moduli scheme for its 
moduli functor ${\bM}(a,b)$.

Given $a>0$, let us fix $b_0$ so that Proposition \ref{p:mainprop}
holds for $\ell$ and 
$(n,c,s)$, and take $b>b_0$.

\begin{lemma} The \gif\  induces morphisms of
functors 
$$
\what\bS^0\colon\bar\bJ_{p,\ell}^r \hookrightarrow
\overline{\bM}(a,b)\,,
\qquad
\what\bS^0\colon\bJ_{p,\ell}^r \hookrightarrow {\bM}(a,b)
$$ that are  representable by open immersions.
\label{l:mainfunc}
\end{lemma}
\begin{thm} The \gif\  gives a morphism
$\what\bS^0\colon\bar\Jc_{p,\ell}^r \to \overline\M(a,b)$ of
schemes that induces an isomorphism
$$
\what\bS^0\colon 
\Jc_{p,\ell}^r
\iso\M'_{p,\ell}(a,b)\,,
$$  where  $\M'_{p,\ell}(a,b)$ is the open subscheme of those sheaves
in $\M(a,b)$ whose spectral cover is finite over $S=B\times
T$ and verifies $\chi(C_t)=1-p$,
$C_t\cdot\Theta=\ell$ for every $t\in T$.
\label{th:main2} \qed
\end{thm}

\begin{remark} A similar result to Proposition \ref{p:mainprop} about preservation of absolute stability is also true for elliptic Calabi-Yau threefolds  \cite{AnHR03}. Similar results to Theorem \ref{th:main2} are true as well.
\label{r:preservation}
\end{remark}

\section{Applications to string theory and mirror symmetry}\label{s:applmirror}

We describe here a few applications of \FMF\ theory to string theory and mirror symmetry.

\subsection{Generalities on D-branes on Calabi-Yau manifolds}

The name D-brane is a contraction of Dirichlet brane. D-branes occur in type IIA, type IIB and
type I string theory as dynamical objects on which strings can end. The coordinates of
the attached strings satisfy Dirichlet boundary conditions in the directions normal to the brane
and Neumann boundary conditions in the directions tangent to the brane. Further, a D$_p$-brane
is $p$ dimensional, where $p$ is even for type IIA strings, odd for type IIB strings and $1,5$ or $9$ for  type I strings. 

One typically distinguishes between two types of D-branes on Calabi-Yau manifolds
called A-type or B-type D-branes. A-type D-branes occur 
in IIB string theory and B-type D-branes in IIA string theory.  All topological invariants of 
B-type D-brane are given by an element 
of a particular K-theory group \cite{MM97, Wit98}. The D-brane RR-charge can then be written \cite{MM97} as:
\begin{equation}
Q(E)=\ch(i_{!}E)\sqrt{\hat{A}(X)}\,.
\label{eq:charge}
\end{equation}
where $\hat{A}(X)$ denotes the $A$-roof genus; note that on a Calabi-Yau manifold
$\hat{A}(X)=\op\td (X)$. Furthermore, in order to give a supersymmetric configration the mass $M$ of a D-brane and its central charge have to satisfy the inequality $M\geq |Z(Q)|$. Here $Z(Q)$ 
is the Bogomolnyi-Prasad-Sommerfield ($BPS$)\footnote{Note that a BPS state is a state that is invariant
under a nontrivial subalgebra of the full supersymmetry algebra. BPS states carry conserved charges.
The supersymmetry algebra determines the mass of the state in terms of its charges. Formally speaking, the central charge is an operator (or constant) that appears
on the right-hand side of a Lie algebra and commutes with all operators in the algebra (for example
Virasoro algebra or supersymmetry algebras). In $N=2$ supersymmetric Seiberg-Witten theory it has been shown that the central charge becomes geometrically; depending there 
on the periods of a particular elliptic curve and electric and magnetic charges. For details
see for example \cite{Wit96a}. This concept carries over for $N=2$ type II string theories on Calabi-Yau threefolds; the central charge depending now on the periods of the holomorphic three-form and
the RR-charges.} central charge is defined in terms of the prepotential $F$ (a function on the complexified K\"ahler moduli space) and $Q(E)$.

\begin{defin} (B-type at large volume). A holomorphic D-brane on a Calabi-Yau manifold $X$ is given by a triple $(C, E, \nabla)$ where $C$ is a holomorphic submanifold
of $X$  and $\nabla$ is a holomorphic connection on $E$ and so $E$ a holomorphic vector bundle.
\end{defin}
If $C=X$, that is, if the D-brane is wrapped over $X$ then supersymmetry requires that $\nabla$ has to
satisfy the hermitian Yang-Mills equations and thus $E$ has to be a $\mu$-stable vector bundle. If the 
branes are wrapped around holomorphic submanifolds of $X$ then the hermitian Yang-Mills equations
must be replaced by a generalization of the Hitchin equations \cite{HaMo98}. More precisely, the gauge fields
which are polarized transverse to $C$ are replaced by ``twisted'' scalars $\Phi$; these are one forms
in the normal bundle of $C$ in $X$ \cite{BVS96}.
\vskip 0.2cm

The above point of view can be generalized if one takes into account that a holomorphic vector bundle defined on a holomorphic submanifold $C$ defines a coherent sheaf ${\G}=i_\ast E$ (with $i\colon C\hookrightarrow X$ being the inclusion map).

Moreover, Kontsevich's homological mirror symmetry conjecture,  
tachyon condensation \cite{Sen} and the fact that $K$-theory \cite{Wit98} classifies D-brane charges has led to the proposal that Kontsevich's mirror conjecture could be physically realized via off-shell states in the open string $B$-model; this suggested to consider D-branes as objects in the derived category of coherent sheaves \cite{Do01a,Do01b} are then to represent those objects the D-brane/anti-D-brane configurations (whereas maps between objects is the derived categories are represented by tachyons and localization on quasi-isomorphisms is expected to be realized by renormalization-group flow). By all progress, there are still many open problems in identifying open string B-model boundary states with objects in $D^b(X)$, for instance, it is hard to confirm that the localization on quasi-isomorphisms is actually realized physically by renormalization-group flow. For a recent 
review on these developments see \cite{As04a}. Formulating a general definition of D-branes
is still an open problem, however, a preliminary definition can be given as follows:

\begin{defin} (B-type). A B-type BPS brane on a Calabi-Yau manifold $X$ is a $\Pi$-stable
object in the bounded derived category $D^b(X)$ on $X$. 
\end{defin}

At the large volume limit where D-branes are represented by coherent sheaves, we can choose any of the available stability notions for them, like slope stability, Gieseker stability or in the generalized Simpson approach. However, we know that as we move away from the large volume limit, D-branes are no longer represented by sheaves but rather they are objects of the derived category. Indeed, the transformation mirroring certain symplectic automorphisms on the sLag side (like Kontsevich monodromies) are conjecturally the automorphisms of the derived category, and there are evidences in this direction \cite{ACHY04, DOPW, ACHY01, Ho01}. Then a stable sheaf is transformed in an object of the derived category and we need a notion of stability for those new objects. 
Douglas made the first attempt to define a notion of stability for D-branes, called $\Pi$-stability  \cite{Do01b,Do01a, Do02} so that stable branes correspond to BPS states.
$\Pi$-stability has been originally introduced in \cite{DFR} as a generalization of $\mu$-stability;
in particular, it has been shown that in the large volume limit $\Pi$-stability reduces to $\mu$-stability, 
respectively, to the $\theta$-stability notion at the orbifold point at which the description of D-branes
involves supersymmetric gauge theories constructed from quivers. Since $\Pi$-stability 
depends on the periods of the Calabi-Yau manifold, it can be used to predict the lines of marginal
stability in the compactification moduli space, at which the $BPS$ spectrum is expected to jump.
This has been analyzed in \cite{AsDo02,As04}. 

There have been many attempts to make this notion rigourous from a mathematical point of view. It seems that one cannot define when an object is stable, but define all the stable objects as a whole.  In other words, one can define certain special subcategories of the derived category whose objects would correspond to the BPS branes. The key references are Bridgeland papers \cite{Bri03pp,Bri02pp}  where the notion of stability for a triangulated category $D$ is established and it is shown that there is a complex manifold $\op{Stab} (D)$ parametrizing stability conditions on $D$.  When $D=D^b(X)$ is the derived category of a Calabi-Yau manifold, then $\op{Stab} (D^b(X))$ is a finite dimensional complex manifold on which $\op{Aut} (X)$ acts naturally. The points of $\op{Stab} (D^b(X))$ correspond to t-structures on $D^b(X)$ together with some extra data defined by Harder-Narasimhan filtrations. 

The space $\op{Stab} (D^b(X))/ \op{Aut} (X)$ is proposed by Bridgeland as the first approximation to the stringy K\"ahler moduli space. 
The problem of computing this space its very hard; in the case of K3 surfaces, Bridgeland has computed $\op{Stab} (D^b(X))$ \cite{Bri03pp}. However, the group $\op{Aut} (X)$ is still unknown. Very recently stability conditions for projective spaces and del Pezzo surfaces has been also founded \cite{Macri}. The very interesting case of open Calabi-Yau threefolds is in progress.

\begin{defin} (A-type). A special Lagrangian D-brane in a Calabi-Yau manifold $Y$
is given by a triple $(\Sigma, E, \nabla)$ where $\Sigma$ is a special Lagrangian submanifold
of $Y$ and $E$ a flat vector bundle on $\Sigma$ with a flat connection $\nabla$.
\end{defin}

\begin{remark} 
Note that $\Sigma$ is said to be special Lagrangian if the following conditions hold:
\begin{align*}
\rest {\omega}{\Sigma}&=0\\
\op {Re}(e^{i\theta}\rest{\Omega}{\Sigma})&=0
\end{align*}
where $\Omega$ is the holomorphic three-form, $\omega$ the K\"ahler form 
and $\theta$ is an arbitrary phase. Equivalently to the second equation, one can require that $\Omega$ pulls back to a constant multiple of the volume element on $\Sigma$.
\end{remark}
\begin{remark}
A more precise definition of A-type branes has to take as well into account destabilizing quantum effects arising from open string tadpoles \cite{As04a}. Also one expects A-type branes to be objects of the Fukaya category of $Y$, and so B-type branes should be 
mapped to A-type branes, using the homological mirror symmetry conjecture of Kontsevich, Furthermore, for A-type branes the central charge is the integral of the holomorphic three-form over the Lagrangian submanifold $\Sigma$.  
\end{remark}
\begin{remark}
A given A-type or B-type D-brane can be deformed by deforming the submanifold
and the bundle with its connection such that the deformations respect the BPS-condition, i.e., a given D-brane stays supersymmetric (the submanifold is holomorphic or special 
Lagrangian and the bundle holomorphic or flat, respectively). So the space of all continues deformations of the triple, say ($C, E, \nabla$), is the D-brane moduli space, denoted by
${\M}_{\rm hol}$ for B-type branes and ${\M}_{\rm slag}$ for A-type branes. 
If we map a A-type brane in ${\M}_{\rm slag}$ to its submanifold $\Sigma$, we can define 
a fibration $p\colon {\M}_{\rm slag}(\Sigma, E, \nabla)\to{\M}_{\rm slag}(\Sigma, Y)$, where 
${\M}_{\rm slag}(\Sigma, Y)$ denotes the space of all continues deformations of the submanifold
$\Sigma$. Having defined this fibration, one can ask whether the space of one-forms on $\Sigma$
can be identified with the tangent space of ${\M}_{\rm slag}(\Sigma, Y)$ (preserving thereby 
the special Lagrangian condition).  A theorem by McLean \cite{McL} states that 
first-order deformations of a special Lagrangian map $f\colon \Sigma\to Y$ are canonically identified with $H^1_{\rm DR}(\Sigma, {\mathbf R})$ and that all first-order deformations of $f\colon \Sigma\to Y$ can be extended to actual deformations implying that the moduli space ${\M}_{\rm slag}(\Sigma, Y)$ of special Lagrangian maps from $\Sigma$ to $Y$ is a smooth manifold of dimension $b_1(\Sigma)$.
McLean's result together with the generalized mirror conjecture of Kontsevich suggests \cite{Va99}
the identification of the number $h^1(\Sigma)$ of complex moduli (assuming here a pairing of the number $b_1(\Sigma)$ of real moduli with the same number of real moduli of the $U(1)$ bundle) with the number $h^1(C, \sEnd(E))$ of
moduli of the vector bundle over $C$ on the holomorphic side. So one expects that actually the dimensions of A-type and B-type D-brane moduli spaces agree.
\end{remark}

\subsection{T-duality as a relative \FM}

Let us now discuss how the relative Fourier-Mukai transformation acts on D-branes on elliptic Calabi-Yau threefolds.
It is known  that T-duality on the elliptic fiber maps in general
(the subscripts indicate whether fiber $\mathfrak f$ or base $B$ is contained (resp. contains) the wrapped world-volume)
\begin{eqnarray}\label{Dtual}
D6&\rightarrow& \tilde D4_B\cr
D4_B\rightarrow \tilde D6&,&
 D4_{\mathfrak f}\rightarrow\tilde D2_B \cr
D2_B\rightarrow \tilde D4_{\tilde{\mathfrak f}}&,&
 D2_{\mathfrak f}\rightarrow  \tilde D0\cr
D0&\rightarrow& \tilde D2_{\tilde{\mathfrak f}}
\end{eqnarray}

One can describe T-duality on the $T^2$ fiber maps given in \ref{Dtual} at the sheaf level.
For this let us consider the skyscraper
sheaf ${\bbC}(x)$ at a point $x$ of $X$. It is a WIT$_0$ sheaf and its FM
transform $\bS^0({\bbC}(x))$ is a torsion-free rank one sheaf $\Lcl_x$ on the
fiber of $X$ over $p(x)$, because with the identification $X\simeq\tilde X$
the point
$x$ corresponds precisely to $\Lcl_x$ (see \cite{BBHM98} or \cite{HMP02}) as we
expect from \ref{Dtual}
and thus we see $D0\rightarrow \tilde D2_{\tilde{\mathfrak f}}$.

For the topological invariants we have indeed $n=x=a=0, S=\eta=0, s=1$ and then
\begin{equation}\label{fmaptw}\ch_i(\bS^0({\bbC}(x)))=0, \quad i=0,1,3, \qquad
\ch_2(\bS^0({\bbC}(x)))=\mathfrak f
\end{equation}

If we start with $\Oc_\Theta$; proceeding as in (3.16) of
\cite{HMP02} we
have\footnote{The formulae differ from those in \cite{HMP02} because we are using
a different Poincar\'e sheaf}

\begin{eqnarray}\label{fmsigma}
\bS^0(\Oc_\Theta)&=&\Oc_X\,,\qquad \bS^1(\Oc_\Theta)=0\cr \bS^0(\Oc_X)&=&0\,,\quad\qquad \bS^1(\Oc_X)=\Oc_\Theta\otimes p^\ast
\omega_B
\end{eqnarray}
Then $\Oc_\Theta$ transforms to the structure sheaf of $X$ and $\Oc_X$ transforms to a line bundle on $\Theta$ as we expect from
\ref{Dtual}
since $D4_B\leftrightarrow \tilde D6$. We have as before the transformations at the
cohomology level;
$$
n=0, \quad x=1, \quad S=0, \quad \eta=\frac 12 c_1, \quad
a= 0, \quad s=\frac 16 \Theta c_1^2
$$
and then we get
$$
\ch_0(\bS^0(\Oc_\Theta))=1,\quad
\ch_i(\bS^0(\Oc_\Theta))=0,\quad i=1,2,3
$$
 
Finally, let us consider a sheaf $\Fc$ on $B$; by \ref{fmsigma} we have
\begin{eqnarray}
 \bS^0(\Oc_\sigma\otimes
p^\ast \Fc)=p^\ast \Fc&,& \quad \bS^1(\Oc_\Theta\otimes
p^\ast \Fc)=0\cr 
\bS^0(p^\ast \Fc)=0 \;\quad &,& \quad \bS^1(p^\ast \Fc)=\Oc_\Theta\otimes p^\ast \Fc\otimes p^\ast \omega_B
\end{eqnarray}
Then, a sheaf $\Oc_\Theta\otimes p^\ast \Fc=\sigma_\ast \Fc$ supported on a curve
$\tilde C$ in $B$ embedded in $X$ via the section $\sigma$ transforms to a  sheaf on the
elliptic surface supported on the inverse image of $\tilde C$ in $X$ and
vice versa.
This is what we expected form the map $D2_B\leftrightarrow \tilde D4_{\tilde{\mathfrak f}} $ of \ref{Dtual}.

Then at the sheaf level we have the
relations \ref{Dtual} appropriate for the fiberwise T-duality on D-branes
\begin{eqnarray}\label{Dtuals}
D4_B&\rightarrow& \tilde{D6}\cr
D2_B&\rightarrow& \tilde{D4}_{\tilde{\mathfrak f}}\cr
D0&\rightarrow& \tilde{D2}_{\tilde{\mathfrak f}}
\end{eqnarray}

\subsubsection*{Adiabatic character of T-duality}

T-duality on fibers has an adiabatic character, that is, 
by using a decomposition of the cohomology into base and fiber parts the operation of the fiberwise duality on the cohomology will be seen to take the form
one gets from an adiabatic extension of the same operation on the
cohomology of a torus of complex dimension 1 (an elliptic non-singular curve),  fulfilling the expectations from the interpretation as T-duality on D-branes. 

Let us see that the action of the \FM\ in cohomology, that is $\ch(\bS(E))$ and $\ch(\hat\bS(E))$ as described in (\ref{e:chVS}) and (\ref{e:chVSS}), once an appropriated twisted charge is introduced, shows the desired adiabatic character of T-duality.


We now modify the action of \FM\ in cohomology by twisting with an appropriate charge. To this end, we introduce the \emph{effective charge} of a D-brane state $G\in D(X)$ by 
\begin{equation}
 Q(G)=\ch(G)\cdot\sqrt{\op\td (X)}
\label{e:effectcharge}
\end{equation}
in agreement with \eqref{eq:charge} and consider the so-called ${\mathbf f}$-map.
This is the map ${\mathbf f}\colon H(X,\bbQ)\break \to H(X,\bbQ)$ given by
$$
{\mathbf f}(x)=\hat\pi_\ast (\pi^\ast(x)\cdot Z)\,,\text{where\ } Z= \sqrt{\hat\pi^\ast\op\td (X)}\cdot \ch(\Pc)\cdot \sqrt{\pi^\ast\op\td (X)}
$$

If $x=0$ ($G$ has degree zero on
fibers) and its Chern characters $\ch_1(G)$ and $\ch_2(G)$ belong respectively to $\bbQ\Theta \oplus H^2(B, \bbQ)$ and $H^2(B, \bbQ)\Theta \oplus \bbQ$,\footnote{For many families of elliptic Calabi-Yau threefolds for which a mirror family has been constructed (cf. \cite{Can94, Can94a}, all elements in $H^2(X,\bbQ)$ and $H^4(X,\bbQ)$ are of this form} one sees from (\ref{e:chVS}) that the effective charge of $G$ transforms to
$$
\bar{\mathbf f}(Q(G))=M\cdot Q(G) 
$$
where
$$
M=\begin{pmatrix}
0&1&0&0&0&0\\ -1&0&0&0&0&0\\ 0&0&0&1&0&0\\  0&0&-1&0&0&0\\
0&0&0&0&0&1\\  0&0&0&0&-1&0
\end{pmatrix}
$$
Then the map between the effective charges of the D-brane state defined
by $G$ and its FM transform
\emph{can  be exhibited as the transformation matrix given by the
adiabatic extension
$M$ of the usual T-duality matrix} 
$$
\begin{pmatrix} 0&1\\ -1&0
\end{pmatrix}
$$ 
on the fiber. The proof that the latter matrix is the matrix for T-duality on one non-singular fiber, that is, for the \FM\ on a non-singular elliptic curve was stated in that form by the first time in \cite{Bri98}.

\subsection{D-branes and homological mirror symmetry}

We recall here that Kontsevich proposed a homological mirror symmetry \cite{Kon95} for a pair $(X,Y)$ of mirror dual Calabi-Yau manifolds; it is conjectured
that there exists a categorical equivalence between the bounded derived 
category $D(X)$ and Fukaya's $A_{\infty}$ category ${\Fc}(Y)$ which is defined by using the symplectic structure on $Y$. An object of ${\Fc}(Y)$ is a special Lagrangian submanifolds with a flat $U(1)$ bundle on it. If we consider a family of manifolds $Y$ the object of ${\Fc}(Y)$ undergoes monodromy transformations. On the other side, the object of $D(X)$ is a complex of coherent sheaves on $X$ and under the categorical equivalence between $D(X)$ and ${\Fc}(Y)$ the monodromy (of three-cycles) is mapped to certain self-equivalences in $D(X)$. 
We listed in subsection \ref{ss:fmfs} some examples of $D(X)$ self-equivalences provided by \FMF s and mentioned Orlov's theorem \ref{th:orlov} saying that any self-equivalence of $D(X)$ is a \FMF.

Now since all elements in ${D}(X)$ may be represented by suitable
complexes of vector bundles on $X$, we can consider the topological 
K-group and the image $K_{\rm hol}(X)$ of ${D}(X)$. The Fourier-Mukai
transform $\bS^{\Ec}\colon {D}(X)\to {D}(X)$ induces then a corresponding 
automorphism $K_{\rm hol}(X)\to K_{\rm hol}(X)$ and also an automorphism
on $H^{\rm even}(X,{\bbQ})$ if one uses the Chern character ring homomorphism
${\rm ch}:K(X)\to H^{\rm even}(X, {\bbQ})$. With this in mind one can now introduce
various kernels and their associated monodromy transformations. 

For instance, let $D$ be a divisor in $X$ and consider the kernel ${\Oc}_{\Delta}(D)$
with $\Delta$ being the diagonal in $X\times X$; the corresponding Fourier-Mukai
transform acts on an object $G\in {D}(X)$ as twisting by a line bundle $G\otimes {\Oc}(D)$, this automorphism is then identified with the monodromy about the large complex
structure limit point (LCSL-point) in the complex structure moduli space ${\Mc}_{\bbC}(Y)$.

Furthermore, considering the kernel given by the ideal sheaf  $\mathcal{I}_\Delta$ on $\Delta$ one 
has that the effect of $\bS^{\mathcal{I}_\Delta}$ on $H^{\rm even}(X)$ can be 
expressed by taking the Chern character ring homomorphism:
\begin{equation}
\ch (\bS^{{\mathcal I}_\Delta}(G))=\ch_0(\bS^{\Oc_{X\times X}}(G))-\ch  (G)
=\left(\int \ch(G)\cdot \op\td (X)\right)-\ch(G)
\label{coni}
\end{equation}

Kontsevich proposed that this automorphism should reproduce the monodromy
about the principal component of the discriminant of the mirror family $Y$. At the principal
component we have vanishing $S^3$ cycles (and the conifold singularity) thus this monodromy may
be identified with the Picard-Lefschetz formula.

Now given a pair of mirror dual Calabi-Yau threefolds, like the examples given by Candelas and others \cite{Can91,Can97,Can94,Can94a,Can94b}, and using the fact that 
$A$-type and $B$-type D-branes get exchanged under mirror symmetry, we can make Kontsevich's correspondence between automorphisms of $D(X)$ and monodromies 
of three-cycles explicit. 

For this we first choose a basis for the three-cycles $\Sigma_i\in H^3(Y,{\Z})$ such that the intersection form takes the canonical form $\Sigma_i\cdot \Sigma_j=\delta_{j,i+b_{2,1}+1}=\eta_{i,j}$ for $i=0,...,b_{2,1}$. It follows that a three-brane wrapped about the cycle $\Sigma=\sum_i n_i\Sigma^i$ has an (electric,magnetic) charge vector ${\bf n}=(n_i)$. The periods of the holomorphic three-form $\Omega$
are then given by
\begin{equation}
\Pi_i=\int_{\Sigma_i}\Omega
\label{per}
\end{equation}  
and can be used to provide projective coordinates on the complex structure moduli space;
more precisely, if we choose a symplectic basis $(A_i,B_j)$ of $H_2(Y,{\Z})$ then the 
$A_i$ periods serve as projective coordinates and the $B_j$ periods satisfy 
the relations $\Pi^j=\eta_{i,j}\partial {F}/\partial \Pi^i$ where ${F}$ is the prepotential which has near the
large radius limit the asymptotic form (cf. \cite{Can94a, HKTY, Ho98})
$$
{F}=\frac 16 \sum_{abc}k_{abc}t_at_bt_c+\frac 12 \sum_{ab}c_{ab}t_at_b-
\sum_a \frac {c_2(X)J_a}{24} t_a+ \frac{\zeta(3)}{2(2\pi i)^3} \chi(X)+{\Oc}(q).
$$
where $\chi(X)$ is the Euler characteristic of $X$, $c_{ab}$ are rational constants (with $c_{ab}=c_{ba}$) reflecting an $Sp(2h^{11}+2)$ ambiguity and $k_{abc}$ is the classical triple intersection number given by
$$
k_{abc}=\int_X J_a\wedge J_b\wedge J_c.
$$
The periods determine the central charge $Z({\bf n})$ of a three-brane wrapped about the cycle $\Sigma=\sum_in_i[\Sigma_i]$ 
\begin{equation}
Z({\mathbf n})=\int_{\Sigma}\Omega=\sum_i n_i\Pi_i,
\label{cent}
\end{equation}
following the conventions used in the literature, we will write $\sum_in_i\Pi_i=n_6\Pi_1+n_4^1\Pi_2+
n_4^2\Pi_3+n_0\Pi_4+n_2^1\Pi_5+n_2^2\Pi_6$. The $\Pi$'s are given by the associated period vector
$$
\begin{pmatrix} \Pi_1\\ \Pi_2\\ \Pi_3 \\ \Pi_4 \\ \Pi_a \end{pmatrix}=
\begin{pmatrix}  \frac16 k_{abc}t_at_bt_c+\frac{c_2(X)J_a}{24}t_a\\
-\frac12 k_{abc}+c_{ab}t_b+\frac{c_2(X)J_a}{24}\\
1\\ t_a \end{pmatrix} \,.
$$

On the other side, the central charge associated with an object $E$ of $D(X)$ 
is given by \cite{HaMo98} 
\begin{equation}
Z(E)=-\int_X e^{-t_aJ_a} \ch(E) (1+\frac{c_2(X)}{24}).
\label{cenmu}
\end{equation}
The two central charges are to be identified under mirror symmetry thus by comparing the expressions \ref{cent} with \ref{cenmu}, one obtains a map relating the Chern-classes of $E$ to the D-brane charges ${\mathbf n}$. We find 
\begin{align*}
\ch_0(E)&=n_6\\
\ch_1(E)&=n_4^aJ_a\\
\ch_2(E)&=n_2^b+c_{ab}n_4^a\\
\ch_3(E)&=-n_0-\frac{c_2(X)J_b}{12}n_4^b.
\end{align*}

If we insert the expressions for $\ch(E)$ in \ref{coni}, we find the linear tranformation acting on ${\mathbf n}$
\begin{equation}n_6\to n_6+n_0
\label{shf}
\end{equation}
which agrees with the monodromy transformation about the conifold locus.

Similarly one finds that the monodromy transformation about the large complex structure limit point
corresponds to automorphisms 
\begin{equation}[E]\to [E\otimes {\Oc}_X(D)]
\label{lcsl}
\end{equation}
where $D$ is the associated divisor defining the large radius limit in the K\"ahler moduli space.
Using the central charge identification, the automorphism/monodromy correspondence has been 
made explicit \cite{DiRo,ACHY04,ACHY01} for various dual pairs of mirror Calabi-Yau threefolds (given as hypersurfaces in weighted projective spaces).  This identification provides evidences for Kontsevich's proposal of homological mirror symmetry. 

\subsubsection*{Fiberwise  T-duality as a Kontsevich monodromy}

Another striking consequence of the identification of T-duality for an elliptic Calabi-Yau threefold $b\colon X\to B$ with the relative \FM\ 
$\bS\colon D(X)\to D(X)$ with kernel the Poincar\'e relative sheaf $\Pc$, is that we can describe T-duality as a composition of self-equivalences of $D(X)$ that correspond to Kontsevich monodromies. 
We can then derive the matrix of the action of T-duality on cohomology.

As a consequence of the expression of the Poincar\'e sheaf (\ref{e:poin}), the relative \FM\ or T-duality can be expressed as a composition
$$
\bS=\bS^{\Oc_\Delta(2c_1)}\circ \bS^{\Oc_\Delta(\Theta)}\circ \bS^{j_\ast \Ic}\circ \bS^{j_\ast \Oc_\Delta(\Theta)}
$$
where $j\colon X\fp X\hookrightarrow X\times X$ is the immersion and $\Ic$ is the ideal of the ``relative'' diagonal $X\hookrightarrow
X\times_B X$. All the \gif s in the above formula correspond to Kontsevich monodromies, with the possible  exception of $\bS_{j_\ast {\mathcal I}}$. 
Since $\mathcal{I}=$ ideal of ``relative'' diagonal $X\hookrightarrow
X\times_B X$ and $\mathcal{I}_\Delta=$ ideal of diagonal $X\hookrightarrow
X\times X$, $\bS_{j_\ast {\mathcal I}}$ is not Kontsevich's conifold
$\bS_{\mathcal{I}_\Delta}$, but it is somehow similar.

\subsection{Application to heterotic string theory}

A compactification of the ten-dimen\-sional heterotic string is 
given by a holomorphic, stable $G$-bundle $V$ over a Calabi-Yau manifold $X$. The 
Calabi-Yau condition, the holomorphy and stability of $V$ are a direct consequence of the 
required supersymmetry in the uncompactified space-time. We assume that the underlying
ten-dimensional space $M_{10}$ is decomposed as $M_{10}=M_{4}\times X$ where $M_{4}$ 
(the uncompactified space-time) denotes the four-dimensional Minkowski space and $X$ a six-dimensional compact space
given by a Calabi-Yau threefold. Now let us be more precise: supersymmetry 
requires that the connection $A$ on $V$ satisfies 
$$F_A^{2,0}=F_A^{0,2}=0, \qquad F^{1,1}\wedge J^2=0\,.$$
It follows that the connection has to be a holomorphic connection on a holomorphic 
vector bundle and in addition to satisfy the Donaldson-Uhlenbeck-Yau
equation that has a unique solution if the vector bundle is $\mu$-stable. 
The first topological condition the vector bundle has to satisfy, is 
$$c_1(V)=0\quad ({\rm mod}\ 2)$$
to ensure that the bundle $V$ admits a spin-structure.

In addition to $X$ and $V$ one has to specify a $B$-field on $X$ of field strength $H$. In order 
to get an anomaly free theory, the Lie group $G$ is fixed to be either $E_8\times E_8$ or 
$Spin(32)/{{\Z}_2}$ or one of their subgroups and $H$ has to satisfy the identity 
$$dH=\tr R\wedge R-\Tr F\wedge F$$ 
where $R$ and $F$ are the associated curvature forms of the spin connection on $X$ and the gauge connection on $V$. Also $\tr$ refers to the trace of the composite endomorphism of the tangent bundle to $X$ and $\Tr$ denotes the trace in the adjoint representation of $G$. For any closed four-dimensional submanifold $X_4$ of the ten-dimensional space-time $M_{10}$, the four form
$\tr R\wedge R-\Tr F\wedge F$ must have trivial cohomology since
$$\int_{X_4}dH=\int_{X_4} \tr R\wedge R-\Tr F\wedge F=0.$$
Thus the second
topological condition $V$ has to satisfy is
$$c_2(TX)=c_2(V).$$
A physical interpretation of the third Chern-class can be given as a result of the decomposition of the ten-dimensional space-time into a four-dimensional flat Minkowski space and $X$. 
The decomposition of the corresponding ten-dimensional Dirac operator with values in $V$ shows that: \emph{massless four-dimensional fermions are in one to one correspondence with zero modes of
the Dirac operator $D_V$ on $X$.} 

More precisely, the spectrum of charged matter is directly related
to properties of $X$ and $V$. So, let us obtain the spectrum of 
massless fermions! We start in ten dimensions with the Dirac equation
$iD_{10}\Psi=0=i(D_4+D_X)\Psi$ further making the ansatz $\Psi=\psi(x)\phi(y)$
where $x$ and $y$ are coordinates on $X$ respectively $M_4$. If $\psi$ is an 
eigenspinor of eigenvalue $m$ then it follows $iD_X\psi=m\psi$ and so we get
$(iD_4+m)\phi=0$.
So we learn that $\psi$ looks like a fermion of mass $m$, to a four-dimensional
observer. Thus, massless four-dimensional fermions are in one to one
correspondence with zero modes of the Dirac operator on $X$. The charged 
four-dimensional fermions are obtained from ten-dimensional ones,
which transform under the adjoint of $E_8$. 

Now since massless fermions in four dimensions are related to the zero 
modes of the Dirac operator on $X$, they can be related to
the cohomology groups $H^k(X,V)$. The index of the 
Dirac operator can be written as
$${\rm index}(D)=\int_X \td (X)ch(V)={\frac{1}{ 2}}\int_X c_3(V)\,,$$
equivalently, one can write the index as ${\rm index}(D)=\sum_{i=0}^3(-1)^k \dim H^k(X,V)$.
For stable vector bundles one has $H^0(X,V)=H^3(X,V)=0$ and so
$$\dim H^2(X,V)-\dim H^1(X,V)={\frac{1}{ 2}}\int_X c_3(V)$$
whose absolute value gives the net-number of fermion generations $N_{\rm gen}$. 
It has been observed in nature that $N_{\rm gen}=3$ thus one would like to find vector bundles with $c_3(V)=\pm 6$.  In case $V=TX$,  one has to search for Calabi-Yau threefolds of Euler characteristic 
$\pm 6$. This inspired an earlier work by Tian and Yau \cite{Tyau}.
\begin{note}
For a detailed introduction to string theory and further aspects of the 
heterotic string we refer to \cite{GSW1,GSW2,Polch1} and \cite{Polch2}. For a recent 
discussion of world-sheet stability issues see \cite{SiWi,BeaWi} and \cite{BaSe}. Some aspects
of vector bundles and bundle cohomology, in the context of string theory, have been studied for
example in \cite{AsM97,AsDo,DiacIon,DLO} and in \cite{DoOv}. 
\end{note}

The inclusion of background five-branes changes the topological constrain \cite{SW,DMW,FMW97}. Now from classical Maxwell theory it is known that in four dimensions the magnetic dual
of the electron is the monopol. In the same way one obtains that the magnetic dual of the fundamental 
heterotic string in ten dimensions is the heterotic five-brane. Various five-brane solutions 
of the heterotic string equations of motion have been discussed in \cite{CHS,CHSt}: the gauge five-brane, the symmetric five-brane and the neutral five-brane. It has been shown that the gauge and symmetric five-brane solution involve finite size instantons of an unbroken
non-Abelian gauge group. In contrast, the neutral five-branes can be
interpreted as zero size instantons of the $SO(32)$ heterotic string \cite{WIT96}.

The magnetic five-brane contributes a source term to the Bianchi identity for the three-form $H$,
$$dH=\tr R\wedge R-\Tr F\wedge F-n_5\sum_{\rm five-branes}\delta_5^{(4)}$$
and integration over a four-cycle in $X$ gives
$$c_2(TX)=c_2(V)+[W]\,.$$
The new term $\delta_5^{(4)}$ is a current that integrates to one in the direction 
transverse to a single five-brane whose class is denoted by $[W]$. The class $[W]$ is the Poincar\'e
dual of an integer sum of all these sources and thus $[W]$ should be a integral class, representing 
a class in $H_2(X,{\Z})$.  $[W]$ can be further specified taking into account
that supersymmetry requires that five-branes are wrapped on holomorphic curves thus $[W]$ 
must correspond to the homology class of holomorphic curves. This fact constraints $[W]$ to be
an algebraic class. Further, algebraic classes include negative classes, however, these lead to negative magnetic charges, which are un-physical, and so they have to be excluded. This constraints $[W]$ to be an effective class. This is actualy the unique constrain because, as mentioned above, supersymmetry 
implies that $TX$ and $V$ are holomorphic bundles and since the characteristic classes 
of holomorphic bundles are algebraic, it follows that $[W]=c_2(TX)-c_2(V)$ is algebraic.
Thus for a given Calabi-Yau threefold $X$ the effectivity of $[W]$ constraints the choice of vector bundles $V$.  

In summary, we are looking for stable holomorphic vector bundles on Calabi-Yau threefolds whose characteristic classes satisfy the constraints 
$$
c_1(V)=0\  (\rm{mod}\ 2), \quad
[W]=c_2(TX)+ch_2(V)\,.
$$
However, to follow the spirit of the present paper we will restrict to the discussion of vector bundles on elliptically fibered Calabi-Yau threefolds. 

Three approaches to construct holomorphic vector bundles on elliptically fibered Calabi-Yau threefolds, with structure group the complexification
$G_{\bbC}$ of a compact Lie group $G$, have been introduced in \cite{FMW97}. The parabolic bundle approach applies for any simple $G$. One considers deformations of certain minimally unstable $G$-bundles corresponding to special maximal parabolic subgroups of $G$. The spectral cover approach (i.e., a relative Fourier-Mukai transformation) applies for $SU(n)$ and $Sp(n)$ bundles. The del Pezzo surface approach applies for $E_6$, $E_7$ and $E_8$ bundles and uses the relation between subgroups of $G$ and
singularities of del Pezzo surfaces. Various aspects of these approaches have been further explored
in  \cite{BJPS,FMW99,FMW98,Don98,AsDo,DLO,DOPW,Bri98,BM2,BBHM98,BBH97a,Mac96,And,Cu98,DiacIon}. 

In what follows we require, as in other parts of this work, that our elliptically fibered Calabi-Yau threefold $\pi\colon X\to B$ 
has a section $\sigma$ (in addition to the smoothness
of $B$ and $X$). Remember than this (and the Calabi-Yau condition) restricts the base $B$ to be a Hirzebruch surface
$({ F}_m$, $m\geq 0)$, a del Pezzo surface $(dP_k, k=0,...,8)$, a rational elliptic surface $(dP_9)$, blown-up Hirzebruch surfaces or an Enriques surface \cite{DLO,MV96b}. 

Now following our above discussion we can define a sheaf on $X$ in two ways
\begin{align*}
V&={\bf \Phi}^0(i_\ast \Lcl) & \pi_{C*}\Lcl&=\sigma^\ast V\\
\tilde V&=\hat{\bf \Phi}^0(i_\ast \Lcl) & \pi_{C*}\Lcl&=\sigma^\ast \tilde V\otimes \omega_B
\end{align*}
where $i\colon C\to X$ is the closed immersion of $C$ into $X$ where $\tilde V$ is related to $V$ by 
\cite{AnHR} $\tilde V=\tau^\ast V\otimes \pi^\ast \omega_B$. We can now determine the topological invariants of the Fourier-Mukai transform of a general complex ${\G}$ in the derived category using the expressions
derived in \ref{e:chE}. So if we start with the sheaf $E=i_\ast \Lcl$  with Chern characters given by 
\begin{align*}
\ch_0(i_\ast \Lcl)=0,\ 
\ch_1(i_\ast \Lcl)=n\Theta+\pi^\ast {\eta},\
\ch_2(i_\ast \Lcl)=\Theta \pi^\ast \eta_E+a_E F,\
\ch_3(i_\ast \Lcl)={s_E}
\end{align*}
with $\eta_E, {\eta}\in H^2(B,\bbQ)$, then the Chern characters of the
Fourier-Mukai transform $V={\bf \Phi}^0(i_\ast \Lcl)$ of $i_\ast \Lcl$ are given by
\begin{align*}
\ch_0(V)&=n\\
\ch_1(V)&=\pi^\ast (\eta_E-{\frac{1} {2}}n c_1)\\
\ch_2(V)&=(-\pi^\ast {\eta})\Theta+({s_E} -{\frac{1}{2}}\pi^\ast \eta_E c_1\Theta
+{\frac{1}{12}}n c_1^2\Theta)F\\
\ch_3(V)&=-a_E +\frac12 \Theta c_1\pi^\ast {\eta}.
\end{align*}
Now $\eta_E$, $a_E$ and $s_E$ are not completely arbitrary. 
Also we have not given an explicit expression for $c_1(\Lcl)$ so far. 
For this we analyze the Grothendieck-Riemann-Roch theorem applied
to the $n$-sheeted cover $\pi_C\colon C\rightarrow B$ which gives
$${\ch(\pi_{C*}\Lcl)\td(B)=\pi_{C*}(\ch(\Lcl)\td(C))}$$
and so we find 
$$
 c_1(\sigma^\ast V)+\frac 12 c_1(B)=\pi_{C*}\big(c_1(\Lcl)+\frac{c_1(C)}{2}\big)\,.
$$
For $(1,1)$ classes $\alpha$ on $B$ we have $\pi_{C*}\pi_C^\ast \alpha=n\alpha$ and
$\sigma^\ast $ applied to $V$ gives $c_1(\sigma^\ast V)=\eta_E- \frac12 nc_1(B)$
so we get
$$
\pi_{C*}(c_1(\Lcl))=\pi_{C*}\big(-\frac{c_1(C)}{2}+\frac{\pi^\ast _C\eta_E}{n}\big)
$$
which gives
$$
c_1(\Lcl)=-\frac{c_1(C)}{ 2}+ \frac{\pi^\ast _C\eta_E}{n}+\gamma\,,
$$
where $\gamma\in H^{1,1}(C,{\Z})$ is some cohomology class satisfying $\pi_{C*}\gamma=0\in H^{1,1}(B,{\Z})$. The general solution for $\gamma$ has been derived in \cite{FMW97} and is given by
$\gamma=\lambda(n\Theta_{\vert_{C}}-\pi_C^\ast \eta+n\pi_C^\ast c_1(B))$ with $\lambda$ some rational number which we will specify below. Let us also note $\gamma$ restricted to $S=C\cap \Theta$
is given by $\gamma_{\vert_{S}}=-\lambda\pi^\ast \eta(\pi^\ast \eta-n\pi^\ast c_1(B))\Theta$.

Having fixed $c_1(\Lcl)$ we can now go on and determine $a_E$ and $s_E$ in terms of 
$\eta_E$. For this apply the Grothendieck-Riemann-Roch theorem to $i\colon C\rightarrow X$ 
which gives ${\ch(i_\ast \Lcl)\td(X)=i_\ast (\ch(\Lcl)\td(C))}.$
We note that $i_\ast (1)=C$ and using the fact that $i_\ast (c_1(B)\gamma)=0$ a simple computation 
gives 
\begin{align*}a_E&=\gamma_{\vert_{S}}+\frac 1n\eta_E\eta\\
s_E&=\frac1{24}{nc_1^2}+\frac1{2n}{\eta_E^2}-\varpi\end{align*}
where $\varpi$ is given by
$${\varpi=-\frac1{24}c_1(B)^2(n^3-n)+\frac12\big(\lambda^2-\frac14\big)n\eta(\eta-nc_1(B))}\,.$$

A priori ${\bf \Phi}^0(i_\ast \Lcl)$ gives $U(n)$ vector bundles on $X$ whose properties
have been analyzed in detail in \cite{AnHR04}. To make contact with the work of \cite{FMW97} let us describe the reduction from $U(n)$ to $SU(n)$ and thus recover the second Chern class of an $SU(n)$ vector bundles originally computed in \cite{FMW97} and the third Chern class evaluated  in \cite{And, DiacIon}.

In order to describe the reduction to $SU(n)$ we specify
the class $\eta_E=\frac12 nc_1$ giving $c_1(V)=0$. If we insert this into the above expressions for $a_E$ and $s_E$ we find the new expressions 
\begin{align*}
a_E&=\gamma_{\vert_{S}}+\frac12 c_1\eta\\
s_E&=\frac16 n c_1^2\Theta-\varpi\,,\end{align*}
and so we find the Chern-classes of an $SU(n)$ vector bundle $V$ on the elliptic fibered Calabi-Yau threefold 
$${r(V)= n, \quad c_1(V)= 0, \quad c_2(V)= \pi^\ast (\eta)\Theta+\pi^\ast (\varpi), \quad c_3(V)=-2\gamma_{\vert_{S}}}$$
in agreement with \cite{FMW97, And, DiacIon}. 

Now using the fact that on an elliptically fibered Calabi-Yau threefold we can decompose
$W=W_B+a_f \mathfrak f$ where $W_B$ is the class of a curve in $B$ and $\mathfrak f$ the fiber of $X$, we 
find 
$$
W_B=\Theta\pi^\ast (12c_1(B)-\eta), \quad a_f=c_2(B)+11c_1(B)^2-\varpi\,.
$$

\begin{remark} Note that in \cite{FMW97} it has been shown that for a vector bundle $V$ with structure group $E_8\times E_8$, the number of five-branes $a_f$ agrees with the number $N$ of three-branes required for anomaly cancellation in $F$-theory on a Calabi-Yau fourfold $Y$ given by $\frac{\chi(Y)}{24}$. If the structure group $G$ of $V$ is contained in $E_8$ then the observed physical gauge group in four-dimensions corresponds to the commutant $\tilde{H}$ of $G$ in $E_8$, typically being of ADE-type. This leads to a generalized physical set-up: The heterotic string compactified on $X$ and a pair of $G$-bundles is dual to $F$-theory compactified on a Calabi-Yau fourfold $Y$ with section and a section $\theta$ of ADE singularities. The generalized set-up has been analyzed in \cite{FMW97} and \cite{BJPS,AC,ACu,ACC}. In particular it has been shown that $\frac{\chi(Y)}{24}=a_{E_8}+a_f$ where $a_{E_8}$ corresponds to the number of five-branes associated to the $E_8$ vector bundle and $a_f$ as given above.
\end{remark}
\begin{remark} The discussion of stability of $V$ and $W$ depends on the properties of the defining data $C$ and $\Lcl$. If $C$ is irreducible and $\Lcl$ a line bundle over $C$ then $V$ and $W$ will be vector bundles stable with respect to
\begin{equation}\label{pol}
J=\epsilon J_0+\pi^\ast H_B, \qquad {\epsilon} > 0
\end{equation}
if $\epsilon$ is sufficiently small (cf. \cite[Theorem 7.1]{FMW99} where the statement is proven under the additional assumption that the restriction of $V$ to the generic fiber is regular and semistable). Here $J_0$ refers to some arbitrary K\"ahler class
on $X$ and $H_B$ a K\"{a}hler class on the base $B$. It implies that the bundle $V$ can be taken 
to be stable with respect to $J$ while keeping the volume
of the fiber $\mathfrak f$ of $X$ arbitrarily small compared to the volumes of effective curves
associated with the base. That $J$ is actually a good polarization can be seen by assuming 
$\epsilon=0$. Now one observes that  ${\pi}^\ast H_B$ is not a K\"{a}hler class on $X$ since
its integral is non-negative on each effective curve $C$ in $X$, however, there is one curve, the fiber $\mathfrak f$, where the integral vanishes. This means that ${\pi}^\ast H_B$ is on the boundary of
the K\"{a}hler cone and to make $V$ stable, one has to move slightly into the interior of the K\"{a}hler cone, that is, into the chamber which is closest to the boundary point ${\pi}^\ast H_B$. 
Also we note that although ${\pi}^\ast H_B$ is in the boundary of the K\"ahler cone, we can still define the slope $\mu_{{\pi}^\ast H_B}(V)$ with respect to it. Since $({\pi}^\ast H_B)^2$
is some positive multiple of the class of the fiber $\mathfrak f$, semi-stability with respect to ${\pi}^\ast H_B$ is implied by semi-stability of the restrictions $V {\vert}_\mathfrak f$ to the fibers.
Assume that $V$ is not stable with respect to $J$, then there is a
destabilizing sub-bundle $V' \subset V$ with $\mu_J(V') \ge \mu_J(V)$.
But semi-stability along the fibers says that $\mu_{{\pi}^\ast H_B}(V') \le
\mu_{{\pi}^\ast H_B}(V)$. If we had equality, it would follow that $V'$
arises by the spectral construction from a proper sub-variety of the
spectral cover of $V$, contradicting the assumption that this cover is
irreducible. So we must have a strict inequality $\mu_{{\pi}^\ast H_B}(V')
<\mu_{{\pi}^\ast H_B}(V)$. Now taking $\epsilon$ small enough, we can
also ensure that $\mu_{J}(V') < \mu_{J}(V)$ thus $V'$ cannot destabilize $V$.
\end{remark}

Let us now consider the case that $C$ is flat over $B$. If $C$ is not irreducible than there may exist line bundles such that $V=\Phi^0(i_\ast \Lcl)$ is not stable with respect to the polarization given by \ref{pol}, however,
the condition one has to impose to the spectral data in order that $V$ is a stable sheaf on $X$ with respect to \ref{pol}, has been derived in \cite{AnHR03} (cf. Remark \ref{r:preservation}). Actually, if $C$ is flat over $B$ and $\Lcl$ is a pure dimension
sheaf on $C$ than $V$ is stable with respect to $\bar{J}=\bar\epsilon \Theta+\pi^\ast  H_B$ for sufficiently small $\bar\epsilon$ if and only if $i_\ast \Lcl$ is stable with respect to this polarization. Let us note here that
stability with respect to \ref{pol} for $\epsilon$ sufficiently small is equivalent to stability with respect to $\bar J$ for sufficiently small $\bar\epsilon$ if we take $J_0=a\Theta+b\pi^\ast H_B$ for some positive $a$, $b$. 
Furthermore, note that if $C$ is irreducible and $\Lcl$ is a line bundle the latter condition is automatically satisfied. 

Moreover, $V$ and $W$ are simultaneously stable with respect to $\tilde J$. This is not a surprise because from $ch_i(\tau^\ast  V)=ch(V)$ we know that $V$ and $\tau^\ast V$ are simultaneously stable and from $W=\tau^\ast V\otimes \pi^\ast \omega_B^{-1}$ we know that stability is the same for $W$ and $\tau^\ast V$.

Note that assuming $C$ is flat and $\Lcl$ a pure dimension one sheaf, one finds a larger class of stable sheaves $V$ (and $W$) than originally constructed in \cite{FMW97}, because we do not need an irreducible spectral cover and that the restriction of the bundle (sheaf) to the generic fiber is regular. 

\subsubsection*{Moduli of vector bundles}

The number of moduli of $V$ can be determined in two ways, depending
on whether one works with $V$ directly or with its spectral cover data
$(C,{\mathcal N})$. In the direct approach one is restricted to
so called $\tau$-\emph{invariant} bundles and therefore to a rather special point
in the moduli space, whereas the second approach is not restricted to such a
point. We give a brief review of both approaches; the first approach was originally
introduced in \cite{FMW97} making concrete earlier observations in \cite{Wi86}. The issue
of $\tau$-\emph{invariance} has been also addressed in \cite{DoOv}.

The first approach starts with the index of the $\bar\partial$ operator
with values in $\sEnd(V)$. The index can be evaluated using $$\operatorname{index}(\bar\partial)=\sum_{i=0}^{3}(-1)^i \dim H^i(X,\sEnd(V)).$$
As this index vanishes by Serre duality on the Calabi-Yau threefold,
one has to introduce
a further twist to get a non-trivial index problem. This is usually
given if the Calabi-Yau space
admits a discrete symmetry group \cite{Wi86}. In case of elliptically
fibered Calabi-Yau manifolds
one has such a group $G$ given by the involution $\tau$ coming from the
``sign flip'' in the elliptic
fibers. One assumes that this symmetry can be lifted to an action
on the bundle at least at some
point in the moduli space \cite{FMW97}. In particular the action of $\tau$
lifts to an action on the adjoint bundle $ad(V)$, the
traceless endomorphisms of $\sEnd(V)$. It follows that the index of the
$\bar\partial$ operator generalizes to a character valued index where
for each $g\in G$ one defines
$\operatorname{index}(g)=\sum_{i=0}^3(-1)^{i+1} \operatorname{Tr}_{H^i(X,ad(V))}g$ where
$\operatorname{Tr}_{H^i(X,ad(V))}$
refers to a trace in the vector space $H^i(X,ad(V))$. The particular
form of this index for elliptic Calabi-Yau threefolds has been
determined in \cite{FMW97} (with $g=1+\frac{\tau}{2}$) one finds $\operatorname{index}(g)=
\sum_{i=0}^3(-1)^{i+1} \dim H^i(X, ad(V))_e$ where the subscript
``e'' indicates the projection
onto the even subspace of $H^i(X, ad(V))$. One can compute this index
using a fixed point theorem
as shown in \cite{FMW97}.

The second approach makes intuitively clear where
the moduli of $V$ are coming from, namely, the number of parameters
specifying the spectral cover $C$ and by the dimension of the space of
holomorphic
line bundles $\Lcl$ on $C$. The first number is given by the dimension of
the linear system
$|C|=|n\Theta+\eta|$. The second number is given by the dimension of
the Picard group $Pic(C)=H^1(C,{\Oc}^\ast _C)$ of $C$.
One thus expects the moduli of $V$ to be given by \cite{BDO}
$$
h^1(X, \sEnd(V))=\dim |C|+\dim Pic(C).
$$
If $C$ is an
irreducible, effective, positive divisor in $X$ one can evaluate $h^1(X, \sEnd(V))$
using the Riemann-Roch theorem.

\subsection{Feedback to derived categories}

The interpretation of  (topological) D-branes on a Calabi-Yau threefold as objects of the derived category, has brought new interest to the old problem of ascertaining to what extent the derived category determines the geometry of the variety. 

A lot of work in that direction has been done and a general answer could be that the derived category $D(X)$ contains plenty of information about the variety $X$ itself. The first of the number of milestones in the path that leads  from the derived category of an algebraic variety to the variety itself,  is due to Bondal and Orlov \cite{BO01}:
\begin{prop} Assume that $X$ is a smooth projective variety and let $K_X$ be the canonical divisor.  If either $K_X$ or $-K_X$ is ample ($X$ is of general type or Fano), then $X$ can be
reconstructed from $D(X)$.
\qed\end{prop}

The result is not true for other kinds of algebraic surfaces. Mukai had proven by that time \cite{Muk81} that there exist non isomorphic abelian varieties and also non isomorphic K3 surfaces having equivalent derived categories. In the case of K3 surfaces, Orlov \cite{Or97}  gave a much precise statement, that be thought as a Torelli type theorem for K3's. He proved that two complex K3-surfaces have equivalent derived categories if and only if the trascendental lattices of their cohomology spaces are Hogde-isometric. 

Du to the fact that every isomorphism of derived categories is a \FMF, as we reported in Theorem \ref{th:orlov}, two varieties with isomorphic derived categories are also known as \emph{Mukai partners} or \emph{mirror partners}.  After Mukai, the problem of finding Fourier-Mukai partners has been considered by Bridgeland-Maciocia \cite{BM2} and  Kawamata \cite{Kaw02b}; they have proved  that if $X$ is a smooth projective surface, then there is a finite number of surfaces $Y$ (up to isomorphisms) such that  $D(X)\simeq D(Y)$. 
There are recent works about counting the number of Mukai partners and determining the structure of the finite set of such partners (cf. \cite{HLOY03,Ste} for K3 surfaces or \cite{Ue04} for elliptic surfaces).
For Calabi-Yau threefolds, Bridgeland  \cite{Bri02} has proven the following result:
\begin{prop} If two Calabi-Yau threefolds $X$ and $Y$ are birational, then they have equivalent derived categories, $D(X)\simeq D(Y)$.
\label{p:birational}
\qed\end{prop} 
If a given Calabi-Yau manifold $X$ undergoes a ``flop-transition'' (i.e., if a rational curve $C$ is
blown-down and then blown-up in a transverse direction) one typically finds that the 
topological type of $X$ changes. In particular, if one starts from $X$ and a collection of
holomorphic curves {$C^i$} on $X$ then the second Chern-class $c_2(X)$ of $X$ changes 
according to a theorem by Tian and Yau \cite{Tyau} as
$c_2(X')=c_2(X)+2\sum_i\int_D[C^i]$
where $D$ is an arbitrary divisor and $[C^i]\in H^4(X)$ which is Poincare dual of $C^i$.
From the physics point of view such transitions have been studied in \cite{AsGrMo}
(for a review see also \cite{Gree}) where it has been shown that the conformal field theory 
associated to the singular target space $X_{sing}$ is perfectly well defined (i.e., one can 
``smoothly'' go from $X$ to $X'$). However, if one includes D-branes one expects certain 
jumps. For instance, the Calabi-Yau manifold $X$ can be considered (to some approximation) as the moduli space of $0$-branes thus it must undergo some transition during the flop \cite{Do01a}. These transitions have been analyzed (using $\Pi$ -stability) in \cite{As03}, showing that only 0-branes
associated to $C$ are affected by the flop transition. Moreover, from the physics of the B-model
one expects that the derived categories of $X$ and $X'$ are equal, in agreement with Proposition \ref{p:birational}.

One can make the following conjecture (cf. also \cite{Ca1}):
\begin{conj}
If we have  two Calabi-Yau threefolds $X$ and $Y$ with $D(X)\simeq D(Y)$, then  $X$ is deformation equivalent to a birational model of  $Y$.
\qed
\end{conj}
This conjecture has great interest in string theory as well, since whenever two Calabi-Yau threefolds have the same derived category, they should have the same D-branes. 

The conjecture is still unproven. C\u ald\u araru \cite{Ca1} has found some explicit models of  Fourier-Mukai partners for three-dimensional Calabi-Yau threefolds. In order to give evidences for the conjecture (or in the contrary to disprove it),  it will be very interesting to exhibit new ones. 

When the problem of reconstructing the variety from the derived category is considered for varieties other than Calabi-Yau threefolds, there are important contributions due to Kawamata. He proved that if $X$, $Y$ are smooth projective varieties with $D(X)\simeq D(Y)$, then $n=\dim X=\dim Y$ and if moreover $\kappa(X)=n$ (that is, $X$ is of general type), then there exist birational morphisms $f\colon Z\to X$, $g\colon Z\to Y$ such that $f^\ast K_X\sim g^\ast K_Y$  \cite{Kaw02b}. 

\appendix
\section{Simpson stability for pure sheaves}\label{A:Simpson}

Here we give a brief account of Simpson stability for pure sheaves, introduced to construct moduli spaces of stable and semistable sheaves which may fail to be torsion-free. The interested reader is referred to \cite{Simp96a} for details.

\subsection{Simpson's slope and reduced Hilbert polynomial}
The notions of rank, slope, actually the whole Hilbert polynomial, can be extended to arbitrary sheaves. The definition depends on the choice of a polarization (or an ample divisor in  more algebraic geometry language). Take then a projective scheme $Y$ with a polarization $H$, that we take very ample by the sake of simplicity. This means that there is a closed immersion of $Y$ into a projective space $\Ps^N$ such that $H$ is the intersection of $Y$ with a hyperplane class.  For every sheaf $\Ec$ on $Y$ there is a polynomial $P(\Ec,n)$ with rational coefficients and degree $s=\dim \Supp(\Ec)$ such that
$$
P(\Ec,n)=\chi(Y,\Ec(n))=\sum_{i\ge 0} \dim H^i(Y,\Ec(n))\,, \quad (\Ec(n)=\Ec\otimes\Oc_Y(nH))
$$
This is the \emph{Hilbert polynomial} of $\Ec$ and can be written in the form
$$
P(\Ec,n)= \frac{r(\Ec)}{s!} n^s + \frac{d(\Ec)}{(s-1)!} n^{s-1}+\dots
$$
where $r(\Ec)$ and $d(\Ec)$ are integer numbers. These numbers are close relatives to the  rank  $\rk(\Ec)$ and the degree $\deg(\Ec)$ with respect to $H$ when $Y$ is irreducible and $\Ec$ is torsion-free as we shall see in the following subsection.

Simpson defined the \emph{reduced Hilbert polynomial} and the \emph{slope} of $\Ec$ as the polynomial $p_S(\Ec)$ defined by
$$
p_S(\Ec,n)=\frac{P(\Ec,n)}{r(\Ec)}\,, \quad \mu_S(\Ec)=\frac{d(\Ec)}{r(\Ec)}
$$

\subsection{Pure sheaves and stability}

When $Y$ is irreducible, a sheaf $\Ec$ on $Y$ is torsion-free precisely when it has no subsheaves supported by a subvariety of smaller dimension. The latter property has a sense for any $Y$, and this was taken as the notion that substitutes torsion-freeness when we are on non irreducible varieties. The right definition is
 
\begin{defin}
A coherent sheaf $\Ec$ is pure of dimension 
$s=s(\Ec)$ if the support of any non-zero subsheaf $0\to \Fc\to \Ec$ has dimension $s$ as well.
\end{defin}

With this definition, torsion-free sheaves on connected varieties can be understood as pure sheaves of maximal dimension.

We define $p\le q$ for two rational polynomials whenever $p(n)\le q(n)$ for $n\gg 0$;  then Simpson definition of semistability and stability is as follows
\begin{defin}
A coherent sheaf $\Ec$ on $Y$ is (Gieseker) semistable   if it is pure of dimension $s=s(\Ec)$ and for every non-zero subsheaf 
$0\to \Fc\to \Ec$ one has
$$
p_S(\Fc)\le p_S(\Ec)\,.
$$
A sheaf is (Gieseker) stable if $p_S(\Fc)<p_S(\Ec)$ for every non-zero subsheaf 
$0\to \Fc\to \Ec$ with  $r(\Fc)<r(\Ec)$.
\end{defin}
For slope stability we have
\begin{defin}
A coherent sheaf $\Ec$ on $Y$ is $\mu_S$-semistable  if it is pure of dimension $s=s(\Ec)$ and for every non-zero subsheaf 
$0\to \Fc\to \Ec$ one has
$$
\mu_S(\Fc)\le \mu_S(\Ec)\,.
$$
A sheaf is $\mu_S$-stable if $\mu(\Fc)<\mu(\Ec)$ for every non-zero subsheaf 
$0\to \Fc\to \Ec$ with  $r(\Fc)<r(\Ec)$.
\end{defin}

A direct computation shows that we have implications
$$
\xymatrix{
\text{$\mu_S$-stable} \ar@{=>}[r] \ar@{=>}[d] & \text{Gieseker stable}\ar@{=>}[d] \\
\text{$\mu_S$-semistable}&\text{Gieseker semistable}\ar@{=>}[l]
}
$$

We may wonder whether those stability conditions are equivalent to the usual ones for a torsion-free sheaf $\Ec$ on an irreducible projective variety $Y$. In this case, Gieseker (semi)stability is defined as above though using
$$
p(\Ec,n)=\frac{P(\Ec,n)}{\rk(\Ec)}
$$
as reduced Hilbert plynomial and also $\mu$-(semi)stability is defined as above using the ordinary slope
$$
\mu(\Ec)=\frac{\deg(\Ec)}{\rk(\Ec)}\,.
$$

The answer is that Gieseker (semi)stability for $p_S(\Ec)$ and for $p(\Ec)$ are equivalent and also that $\mu_S$ and $\mu$-(semi)stability are equivalent. The reason is that we have 
\begin{equation}
r(\Ec)=\rk(\Ec)\cdot \deg (Y)\,,\quad d(\Ec)=\deg(\Ec)+ \rk(\Ec) C
\label{e:rank2}
\end{equation}
where $\deg(Y)$ is the degree of $Y$ in $\Ps^N$ and $C$ is a constant. For instance, if $Y$ is a smooth surface and $\Ec$ is locally free then the Hilbert polynomial is
$$
P(\Ec,n)=\frac{H^2\cdot \rk(\Ec)}{2} n^2 + (\deg(\Ec)-\frac12 \rk(\Ec) H\cdot K_X)n+\chi(Y,\Ec)
$$
so that
$$
r(\Ec)=\rk(\Ec) \cdot H^2=\rk(\Ec)\cdot \deg (Y)\,,\quad
d(\Ec)= \deg(\Ec)-\frac12 \rk(\Ec) H\cdot K_X\,,
$$
in agreement with (\ref{e:rank2}).

Now, (\ref{e:rank2}) gives the hint of a sensible definition for the rank of a sheaf on its support
\begin{defin} Let $Y$ be a projective scheme and $H$ polarization in $Y$. The polarized rank of a coherent sheaf $\Ec$ on $Y$ is the rational number
$$
\rk_H(\Ec)=\frac{r(\Ec)}{\deg(\Supp(\Ec))}
$$
where $\Supp(\Ec)$ is the support of $\Ec$.
\label{d:polrank}
\end{defin}
Note that if $Y$ is irreducible and smooth and $\Supp(\Ec)$ is different from $Y$, then the ordinary rank of $\Ec$ (even if defined as the Chern character $\ch_0$ as in (\ref{e:rank}) )is zero, but the polarized rank may be not. If $\Supp(\Ec)$ is irreducible, then $\rk_H(\Ec)=\rk(\rest{\Ec}{\Supp(\Ec)})$.

Simpson constructed moduli spaces for stable and semistable sheaves both for $\mu_S$ and Gieseker stability. Let us talk simply about $\mu_S$-stability because the analogous results for Gieseker stability are also true.

We start by recalling the definition of S-equivalence: it is defined to ensure the existence a coarse moduli space of semistable sheaves with prescribed topological invariants. Actually the moduli space parametrizes S-equivalence classes of semistable sheaves rather than semistable sheaves themselves. 

The precise definition requires the notion of the Jordan-H\"older filtration: every semistable sheaf $\Fc$ has a filtration $\Fc=\Fc_m\supset \Fc_{m-1}\supset\dots\supset\Fc_0=0$ whose quotients $\Fc_i/\Fc_{i-1}$ are stable with the same slope as $\Fc$. The Jordan-H\"older filtration is not unique but the graduate $G(\Fc)=\oplus_i \Fc_i/\Fc_{i-1}$ is uniquely determined by the sheaf. Two semistable sheaves $\Fc$, $\G$ are then called S-equivalent if  $G(\Fc)\simeq G(\G)$.

We then see that two stable sheaves are S-equivalent only when they are isomorphic.

The final existence result is
\begin{thm} {\ }
\begin{enumerate}
\item
There exists a coarse moduli scheme $M^{ss}(Y,p(n))$ for the moduli problem of S-equivalence classes of  semistable pure sheaves with fixed reduced Hilbert polynomial $p(n)$.
\item The moduli scheme $M^{ss}(Y,p(n))$ is projective.
\item The closed points of $M^{ss}(Y,p(n))$ represent S-equivalence classes of semiestable sheaves on $X$ with reduced Hilbert polynomial  $p(n)$
$\Ec\sim\Ec'$ if
$G(\Ec)\simeq G(\Ec')$.
\item The exists an open subscheme $M^s(Y,p(n))\subseteq M^{ss}(Y,p(n))$  whose points represent the isomorphism classes of stable  sheaves.
\end{enumerate}\qed\label{t:fundamental}
\end{thm}

Analogous results hold for families, that is for projective morphisms $Y\to B$ with a relative polarization and flat sheaves which are (semi)stable on the fibers.

The first construction of the moduli space of semistable torsion-free sheaves are due to Mumford \cite{GIT}, Narasimhan and Seshadri \cite{NaSe65,Se67}  for curves (cf. also \cite{New78}), to Gieseker for surfaces
\cite{Gie77} and to Maruyama in arbitrary dimension \cite{Mar77,Mar78}. Simpson construction \cite{Simp96a} though still based as the other ones (with the exception of  \cite{NaSe65}) on Geometric Invariant Theory, is much simpler and works for singular varieties as well.

\medskip\noindent {\bf Acknowledgments.} We thank U.
Bruzzo, C. Bartocci, J.M. Mu\~noz Porras and H. Kurke for useful
discussions and suggestions. 



\begin{thebibliography}{100}

\bibitem{AK80}
{\sc A.~B. Altman and S.~L. Kleiman}, {\em Compactifying the {P}icard scheme},
  Adv. in Math., 35 (1980), pp.~50--112.

\bibitem{And}
{\sc B.~Andreas}, {\em On vector bundles and chiral matter in {$N=1$} heterotic
  compactifications}, J. High Energy Phys.,  (1999), pp.~Paper 11, 11 pp.
  (electronic).

\bibitem{AC}
{\sc B.~Andreas and G.~Curio}, {\em Three-branes and five-branes in {$N=1$}
  dual string pairs}, Phys. Lett. B, 417 (1998), pp.~41--44.

\bibitem{ACu}
\leavevmode\vrule height 2pt depth -1.6pt width 23pt, {\em On discrete twist
  and four-flux in {$N=1$} heterotic/{F}-theory compactifications}, Adv. Theor.
  Math. Phys., 3 (1999), pp.~1325--1413.

\bibitem{ACC}
\leavevmode\vrule height 2pt depth -1.6pt width 23pt, {\em Horizontal and
  vertical five-branes in heterotic/{F}-theory duality}, J. High Energy Phys.,
  (2000), pp.~Paper 13, 20.

\bibitem{AnHR}
{\sc B.~Andreas and D.~Hern{\'a}ndez~Ruip{\'e}rez}.
\newblock In preparation.

\bibitem{AnHR04}
\leavevmode\vrule height 2pt depth -1.6pt width 23pt, {\em {$U(n)$} vector
  bundles on {C}alabi-{Y}au threefolds for string theory compactifications}.
\newblock {\tt hep-th/0410170}.

\bibitem{AnHR03}
\leavevmode\vrule height 2pt depth -1.6pt width 23pt, {\em Comments on {$N=1$}
  heterotic string vacua}, Adv. Theor. Math. Phys., 7 (2003), pp.~751--786.

\bibitem{ACHY04}
{\sc B.~Andreas, S.-T. Yau, G.~Curio, and D.~Hern{\'a}ndez~Ruip{\'e}rez}, {\em
  Fourier-mukai transform and mirror symmetry for {D}-branes on elliptic
  {C}alabi-{Y}au}.
\newblock {\tt math.AG/0012196}.

\bibitem{ACHY01}
\leavevmode\vrule height 2pt depth -1.6pt width 23pt, {\em Fibrewise
  {$T$}-duality for {D}-branes on elliptic {C}alabi-{Y}au}, J. High Energy
  Phys.,  (2001), pp.~Paper 20, 13.

\bibitem{As04a}
{\sc P.~S. Aspinwall}, {\em D-branes on {C}alabi-{Y}au manifolds}.
\newblock {\tt hep-th/0403166}.

\bibitem{As04}
\leavevmode\vrule height 2pt depth -1.6pt width 23pt, {\em D-branes,
  $\pi$-stability and $\theta$-stability}.
\newblock {\tt hep-th/0407123}.

\bibitem{As03}
\leavevmode\vrule height 2pt depth -1.6pt width 23pt, {\em A point's point of
  view of stringy geometry}, J. High Energy Phys.,  (2003), pp.~002, 15.

\bibitem{AsDo}
{\sc P.~S. Aspinwall and R.~Y. Donagi}, {\em The heterotic string, the tangent
  bundle and derived categories}, Adv. Theor. Math. Phys., 2 (1998),
  pp.~1041--1074.

\bibitem{AsDo02}
{\sc P.~S. Aspinwall and M.~R. Douglas}, {\em D-brane stability and monodromy},
  J. High Energy Phys.,  (2002), pp.~no. 31, 35.

\bibitem{AsGrMo}
{\sc P.~S. Aspinwall, B.~R. Greene, and D.~R. Morrison}, {\em Calabi-{Y}au
  moduli space, mirror manifolds and spacetime topology change in string
  theory}, Nuclear Phys. B, 416 (1994), pp.~414--480.

\bibitem{AsM97}
{\sc P.~S. Aspinwall and D.~R. Morrison}, {\em Point-like instantons in {$K3$}
  orbifolds}, Nuclear Phys. B, 503 (1997), pp.~533--564.

\bibitem{At57}
{\sc M.~F. Atiyah}, {\em Vector bundles over an elliptic curve}, Proc. London
  Math. Soc. (3), 7 (1957), pp.~414--452.

\bibitem{BBH97a}
{\sc C.~Bartocci, U.~Bruzzo, and D.~Hern{\'a}ndez~Ruip{\'e}rez}, {\em A
  {F}ourier-{M}ukai transform for stable bundles on {K}3 surfaces}, J. Reine
  Angew. Math., 486 (1997), pp.~1--16.

\bibitem{BBHJ}
{\sc C.~Bartocci, U.~Bruzzo, D.~Hern{\'a}ndez~Ruip{\'e}rez, and M.~Jardim},
  {\em Nahm and {F}ourier-{M}ukai transforms in geometry and mathematical
  physics}.
\newblock To appear in Progress in Mathematical Physics, Birkha\"user, 2005.

\bibitem{BBHM98}
{\sc C.~Bartocci, U.~Bruzzo, D.~Hern{\'a}ndez~Ruip{\'e}rez, and J.~M.
  Mu{\~n}oz~Porras}, {\em Mirror symmetry on {K}3 surfaces via
  {F}ourier-{M}ukai transform}, Comm. Math. Phys., 195 (1998), pp.~79--93.

\bibitem{BBHM02}
\leavevmode\vrule height 2pt depth -1.6pt width 23pt, {\em Relatively stable
  bundles over elliptic fibrations}, Math. Nachr., 238 (2002), pp.~23--36.

\bibitem{BaSe}
{\sc A.~Basu and S.~Sethi}, {\em World-sheet stability of {$(0,2)$} linear
  sigma models}, Phys. Rev. D (3), 68 (2003), pp.~025003, 8.

\bibitem{BeaWi}
{\sc C.~Beasley and E.~Witten}, {\em Residues and world-sheet instantons}, J.
  High Energy Phys.,  (2003), pp.~065, 39 pp. (electronic).

\bibitem{BNR}
{\sc A.~Beauville, M.~S. Narasimhan, and S.~Ramanan}, {\em Spectral curves and
  the generalised theta divisor}, J. Reine Angew. Math., 398 (1989),
  pp.~169--179.

\bibitem{Can94b}
{\sc P.~Berglund, P.~Candelas, X.~de~la Ossa, and et~al.}, {\em Periods for
  {C}alabi-{Y}au and {L}andau-{G}inzburg vacua}, Nuclear Phys. B, 419 (1994),
  pp.~352--403.

\bibitem{BJPS}
{\sc M.~Bershadsky, A.~Johansen, T.~Pantev, and V.~Sadov}, {\em On
  four-dimensional compactifications of {$F$}-theory}, Nuclear Phys. B, 505
  (1997), pp.~165--201.

\bibitem{BVS96}
{\sc M.~Bershadsky, C.~Vafa, and V.~Sadov}, {\em D-branes and topological field
  theories}, Nuclear Phys. B, 463 (1996), pp.~420--434.

\bibitem{BO01}
{\sc A.~I. Bondal and D.~O. Orlov}, {\em Reconstruction of a variety from the
  derived category and groups of autoequivalences}, Compositio Math., 125
  (2001), pp.~327--344.

\bibitem{Bri03pp}
{\sc T.~Bridgeland}, {\em Stability conditions on {K}3 surfaces}.
\newblock {\tt math.AG/0307164}.

\bibitem{Bri02pp}
\leavevmode\vrule height 2pt depth -1.6pt width 23pt, {\em Stability conditions
  on triangulated categories}.
\newblock {\tt math.AG/0212237}.

\bibitem{Bri98}
\leavevmode\vrule height 2pt depth -1.6pt width 23pt, {\em Fourier-{M}ukai
  transforms for elliptic surfaces}, J. Reine Angew. Math., 498 (1998),
  pp.~115--133.

\bibitem{Bri99}
\leavevmode\vrule height 2pt depth -1.6pt width 23pt, {\em Equivalences of
  triangulated categories and {F}ourier-{M}ukai transforms}, Bull. London Math.
  Soc., 31 (1999), pp.~25--34.

\bibitem{Bri02}
\leavevmode\vrule height 2pt depth -1.6pt width 23pt, {\em Flops and derived
  categories}, Invent. Math., 147 (2002), pp.~613--632.

\bibitem{BM2}
{\sc T.~Bridgeland and A.~Maciocia}, {\em Fourier-{M}ukai transforms for {$K3$}
  and elliptic fibrations}, J. Algebraic Geom., 11 (2002), pp.~629--657.

\bibitem{BDO}
{\sc E.~Buchbinder, R.~Donagi, and B.~A. Ovrut}, {\em Vector bundle moduli and
  small instanton transitions}, J. High Energy Phys.,  (2002), pp.~no. 54, 44.

\bibitem{Ca1}
{\sc A.~C{\u{a}}ld{\u{a}}raru}, {\em Fiberwise stable bundles on elliptic
  threefolds with relative {P}icard number one}, C. R. Math. Acad. Sci. Paris,
  334 (2002), pp.~469--472.

\bibitem{CHSt}
{\sc C.~G. Callan, Jr., J.~A. Harvey, and A.~Strominger}, {\em Worldbrane
  actions for string solitons}, Nuclear Phys. B, 367 (1991), pp.~60--82.

\bibitem{CHS}
\leavevmode\vrule height 2pt depth -1.6pt width 23pt, {\em Worldsheet approach
  to heterotic instantons and solitons}, Nuclear Phys. B, 359 (1991),
  pp.~611--634.

\bibitem{Can94}
{\sc P.~Candelas, X.~de~la Ossa, A.~Font, S.~Katz, and D.~R. Morrison}, {\em
  Mirror symmetry for two-parameter models. {I}}, Nuclear Phys. B, 416 (1994),
  pp.~481--538.

\bibitem{Can97}
\leavevmode\vrule height 2pt depth -1.6pt width 23pt, {\em Mirror symmetry for
  two parameter models. {I}}, in Mirror symmetry, II, vol.~1 of AMS/IP Stud.
  Adv. Math., Amer. Math. Soc., Providence, RI, 1997, pp.~483--543.

\bibitem{Can91}
{\sc P.~Candelas, X.~C. de~la Ossa, P.~S. Green, and L.~Parkes}, {\em A pair of
  {C}alabi-{Y}au manifolds as an exactly soluble superconformal theory},
  Nuclear Phys. B, 359 (1991), pp.~21--74.

\bibitem{Can94a}
{\sc P.~Candelas, A.~Font, S.~Katz, and D.~R. Morrison}, {\em Mirror symmetry
  for two-parameter models. {II}}, Nuclear Phys. B, 429 (1994), pp.~626--674.

\bibitem{Cu98}
{\sc G.~Curio}, {\em Chiral matter and transitions in heterotic string models},
  Phys. Lett. B, 435 (1998), pp.~39--48.

\bibitem{De}
{\sc P.~Deligne}, {\em Courbes elliptiques: formulaire d'apr\`es {J}. {T}ate},
  in Modular functions of one variable, IV (Proc. Internat. Summer School,
  Univ. Antwerp, Antwerp, 1972), Springer, Berlin, 1975, pp.~53--73. Lecture
  Notes in Math., Vol. 476.

\bibitem{DiacIon}
{\sc D.-E. Diaconescu and G.~Ionesei}, {\em Spectral covers, charged matter and
  bundle cohomology}, J. High Energy Phys.,  (1998), pp.~Paper 1, 15 pp.
  (electronic).

\bibitem{DiRo}
{\sc D.-E. Diaconescu and C.~R{\"o}melsberger}, {\em D-branes and bundles on
  elliptic fibrations}, Nuclear Phys. B, 574 (2000), pp.~245--262.

\bibitem{Don98}
{\sc R.~Donagi}, {\em Taniguchi lectures on principal bundles on elliptic
  fibrations}, in Integrable systems and algebraic geometry (Kobe/Kyoto, 1997),
  World Sci. Publishing, River Edge, NJ, 1998, pp.~33--46.

\bibitem{DLO}
{\sc R.~Donagi, A.~Lukas, B.~A. Ovrut, and D.~Waldram}, {\em Holomorphic vector
  bundles and non-perturbative vacua in {M}-theory}, J. High Energy Phys.,
  (1999), pp.~Paper 34, 46 pp.\ (electronic).

\bibitem{DoOv}
{\sc R.~Donagi, B.~A. Ovrut, T.~Pantev, and D.~Waldram}, {\em Standard-model
  bundles}, Adv. Theor. Math. Phys., 5 (2001), pp.~563--615.

\bibitem{DOPW}
\leavevmode\vrule height 2pt depth -1.6pt width 23pt, {\em Standard models from
  heterotic {M}-theory}, Adv. Theor. Math. Phys., 5 (2001), pp.~93--137.

\bibitem{Do01a}
{\sc M.~R. Douglas}, {\em D-branes, categories and {$N=1$} supersymmetry}, J.
  Math. Phys., 42 (2001), pp.~2818--2843.
\newblock Strings, branes, and M-theory.

\bibitem{Do01b}
\leavevmode\vrule height 2pt depth -1.6pt width 23pt, {\em D-branes on
  {C}alabi-{Y}au manifolds}, in European Congress of Mathematics, Vol. II
  (Barcelona, 2000), vol.~202 of Progr. Math., Birkh\"auser, Basel, 2001,
  pp.~449--466.

\bibitem{Do02}
\leavevmode\vrule height 2pt depth -1.6pt width 23pt, {\em Dirichlet branes,
  homological mirror symmetry, and stability}, in Proceedings of the
  International Congress of Mathematicians, Vol. III (Beijing, 2002), Beijing,
  2002, Higher Ed. Press, pp.~395--408.

\bibitem{DFR}
{\sc M.~R. Douglas, B.~Fiol, and C.~R\"omelsberger}, {\em Stability and bps
  branes}.
\newblock {\tt hep-th/0002037}.

\bibitem{DMW}
{\sc M.~J. Duff, R.~Minasian, and E.~Witten}, {\em Evidence for
  heterotic/heterotic duality}, Nuclear Phys. B, 465 (1996), pp.~413--438.

\bibitem{Fr95}
{\sc R.~Friedman}, {\em Vector bundles and {${\rm SO}(3)$}-invariants for
  elliptic surfaces}, J. Amer. Math. Soc., 8 (1995), pp.~29--139.

\bibitem{FM}
{\sc R.~Friedman and J.~W. Morgan}, {\em Smooth four-manifolds and complex
  surfaces}, vol.~27 of Ergebnisse der Mathematik und ihrer Grenzgebiete (3)
  [Results in Mathematics and Related Areas (3)], Springer-Verlag, Berlin,
  1994.

\bibitem{FMW97}
{\sc R.~Friedman, J.~W. Morgan, and E.~Witten}, {\em Vector bundles and {${\rm
  F}$} theory}, Comm. Math. Phys., 187 (1997), pp.~679--743.

\bibitem{FMW98}
\leavevmode\vrule height 2pt depth -1.6pt width 23pt, {\em Principal
  {$G$}-bundles over elliptic curves}, Math. Res. Lett., 5 (1998), pp.~97--118.

\bibitem{FMW99}
\leavevmode\vrule height 2pt depth -1.6pt width 23pt, {\em Vector bundles over
  elliptic fibrations}, J. Algebraic Geom., 8 (1999), pp.~279--401.

\bibitem{Ful}
{\sc W.~Fulton}, {\em Intersection theory}, vol.~2 of Ergebnisse der Mathematik
  und ihrer Grenzgebiete (3) [Results in Mathematics and Related Areas (3)],
  Springer-Verlag, Berlin, 1984.

\bibitem{Gie77}
{\sc D.~Gieseker}, {\em On the moduli of vector bundles on an algebraic
  surface}, Ann. of Math. (2), 106 (1977), pp.~45--60.

\bibitem{GSW2}
{\sc M.~B. Green, J.~H. Schwarz, and E.~Witten}, {\em Superstring theory.
  {V}ol. 1}, Cambridge Monographs on Mathematical Physics, Cambridge University
  Press, Cambridge, second~ed., 1988.
\newblock Introduction.

\bibitem{GSW1}
\leavevmode\vrule height 2pt depth -1.6pt width 23pt, {\em Superstring theory.
  {V}ol. 2}, Cambridge Monographs on Mathematical Physics, Cambridge University
  Press, Cambridge, second~ed., 1988.
\newblock Loop amplitudes, anomalies and phenomenology.

\bibitem{Gree}
{\sc B.~R. Greene}, {\em String theory on {C}alabi-{Y}au manifolds}, in Fields,
  strings and duality (Boulder, CO, 1996), World Sci. Publishing, River Edge,
  NJ, 1997, pp.~543--726.

\bibitem{EGAIII-1}
{\sc A.~Grothendieck}, {\em \'{E}l\'ements de g\'eom\'etrie alg\'ebrique.
  {III}. \'{E}tude cohomologique des faisceaux coh\'erents. {I}}, Inst. Hautes
  \'Etudes Sci. Publ. Math.,  (1961), p.~167.

\bibitem{Hart66}
{\sc R.~Hartshorne}, {\em Residues and duality}, With an appendix by P.
  Deligne. Lecture Notes in Mathematics, vol. 20, Springer-Verlag, Berlin,
  1966.

\bibitem{Hart77}
\leavevmode\vrule height 2pt depth -1.6pt width 23pt, {\em Algebraic geometry},
  Graduate Texts in Mathematics, vol. 52, Springer-Verlag, New York, 1977.

\bibitem{HaMo98}
{\sc J.~A. Harvey and G.~Moore}, {\em On the algebras of {BPS} states}, Comm.
  Math. Phys., 197 (1998), pp.~489--519.

\bibitem{HMP02}
{\sc D.~Hern{\'a}ndez~Ruip{\'e}rez and J.~M. Mu{\~n}oz~Porras}, {\em Stable
  sheaves on elliptic fibrations}, J. Geom. Phys., 43 (2002), pp.~163--183.

\bibitem{Hit87}
{\sc N.~J. Hitchin}, {\em The self-duality equations on a {R}iemann surface},
  Proc. London Math. Soc. (3), 55 (1987), pp.~59--126.

\bibitem{Ho01}
{\sc P.~Horja}, {\em Derived category automorphisms from mirror symmetry}.
\newblock To appear in Duke Mathematical Journal, {\tt math.AG/0103231}.

\bibitem{Ho98}
{\sc S.~Hosono}, {\em G{KZ} systems, {G}r\"obner fans, and moduli spaces of
  {C}alabi-{Y}au hypersurfaces}, in Topological field theory, primitive forms
  and related topics (Kyoto, 1996), vol.~160 of Progr. Math., Birkh\"auser
  Boston, Boston, MA, 1998, pp.~239--265.

\bibitem{HKTY}
{\sc S.~Hosono, A.~Klemm, S.~Theisen, and S.-T. Yau}, {\em Mirror symmetry,
  mirror map and applications to complete intersection {C}alabi-{Y}au spaces},
  Nuclear Phys. B, 433 (1995), pp.~501--552.

\bibitem{HLOY03}
{\sc S.~Hosono, B.~H. Lian, K.~Oguiso, and S.-T. Yau}, {\em Fourier-{M}ukai
  partners of a {$K3$} surface of {P}icard number one}, in Vector bundles and
  representation theory (Columbia, MO, 2002), vol.~322 of Contemp. Math., Amer.
  Math. Soc., Providence, RI, 2003, pp.~43--55.

\bibitem{Illu90}
{\sc L.~Illusie}, {\em Cat\'egories d\'eriv\'ees et dualit\'e: travaux de
  {J}.-{L}. {V}erdier}, Enseign. Math. (2), 36 (1990), pp.~369--391.

\bibitem{JM}
{\sc M.~Jardim and A.~Maciocia}, {\em A {F}ourier-{M}ukai approach to spectral
  data for instantons}, J. Reine Angew. Math., 563 (2003), pp.~221--235.

\bibitem{Kaw02b}
{\sc Y.~Kawamata}, {\em {$D$}-equivalence and {$K$}-equivalence}, J.
  Differential Geom., 61 (2002), pp.~147--171.

\bibitem{kod}
{\sc K.~Kodaira}, {\em On compact analytic surfaces. {II}, {III}}, Ann. of
  Math. (2) 77 (1963), 563--626; ibid., 78 (1963), pp.~1--40.

\bibitem{Kon95}
{\sc M.~Kontsevich}, {\em Homological algebra of mirror symmetry}, in
  Proceedings of the International Congress of Mathematicians, Vol.\ 1, 2
  (Z\"urich, 1994), Basel, 1995, Birkh\"auser, pp.~120--139.

\bibitem{LB}
{\sc H.~Lange and C.~Birkenhake}, {\em Complex abelian varieties}, vol.~302 of
  Grundlehren der Mathematischen Wissenschaften [Fundamental Principles of
  Mathematical Sciences], Springer-Verlag, Berlin, 1992.

\bibitem{Mac96}
{\sc A.~Maciocia}, {\em Generalized {F}ourier-{M}ukai transforms}, J. Reine
  Angew. Math., 480 (1996), pp.~197--211.

\bibitem{Macri}
{\sc E.~Macr\`\i}, {\em Some examples of moduli spaces of stability conditions
  on derived categories}.
\newblock {\tt math.AG/0411613}.

\bibitem{Mar77}
{\sc M.~Maruyama}, {\em Moduli of stable sheaves. {I}}, J. Math. Kyoto Univ.,
  17 (1977), pp.~91--126.

\bibitem{Mar78}
\leavevmode\vrule height 2pt depth -1.6pt width 23pt, {\em Moduli of stable
  sheaves. {II}}, J. Math. Kyoto Univ., 18 (1978), pp.~557--614.

\bibitem{McL}
{\sc R.~P. McLean}, {\em Deformations of calibrated submanifolds}.
\newblock PhD thesis, Duke University.

\bibitem{MM97}
{\sc R.~Minasian and G.~Moore}, {\em {$K$}-theory and {R}amond-{R}amond
  charge}, J. High Energy Phys.,  (1997), pp.~Paper 2, 7 pp. (electronic).

\bibitem{Mir83}
{\sc R.~Miranda}, {\em Smooth models for elliptic threefolds}, in The
  birational geometry of degenerations (Cambridge, Mass., 1981), vol.~29 of
  Progr. Math., Birkh\"auser Boston, Mass., 1983, pp.~85--133.

\bibitem{MV96b}
{\sc D.~R. Morrison and C.~Vafa}, {\em Compactifications of {$F$}-theory on
  {C}alabi-{Y}au threefolds. {II}}, Nuclear Phys. B, 476 (1996), pp.~437--469.

\bibitem{Muk81}
{\sc S.~Mukai}, {\em Duality between {$D(X)$} and {$D(\hat X)$} with its
  application to {P}icard sheaves}, Nagoya Math. J., 81 (1981), pp.~153--175.

\bibitem{Muk87a}
\leavevmode\vrule height 2pt depth -1.6pt width 23pt, {\em On the moduli space
  of bundles on {$K3$} surfaces. {I}}, in Vector bundles on algebraic varieties
  (Bombay, 1984), Tata Inst. Fund. Res. Stud. Math., vol. 11, Tata Inst. Fund.
  Res., Bombay, 1987, pp.~341--413.

\bibitem{Mum74}
{\sc D.~Mumford}, {\em Abelian varieties}, Tata Inst. Fund. Res. Stud. Math.,
  vol. 5, Tata Inst. Fundamental Res., Bombay, 1970.

\bibitem{GIT}
{\sc D.~Mumford, J.~Fogarty, and F.~Kirwan}, {\em Geometric invariant theory},
  vol.~34 of Ergebnisse der Mathematik und ihrer Grenzgebiete (2) [Results in
  Mathematics and Related Areas (2)], Springer-Verlag, Berlin, third~ed., 1994.

\bibitem{NaSe65}
{\sc M.~S. Narasimhan and C.~S. Seshadri}, {\em Stable and unitary vector
  bundles on a compact {R}iemann surface}, Ann. of Math. (2), 82 (1965),
  pp.~540--567.

\bibitem{New78}
{\sc P.~E. Newstead}, {\em Introduction to moduli problems and orbit spaces},
  vol.~51 of Tata Institute of Fundamental Research Lectures on Mathematics and
  Physics, Tata Institute of Fundamental Research, Bombay, 1978.

\bibitem{Or97}
{\sc D.~O. Orlov}, {\em Equivalences of derived categories and ${K}3$
  surfaces}, J. Math. Sci. (New York), 84 (1997), pp.~1361--1381.
\newblock Algebraic geometry, 7.

\bibitem{Polch1}
{\sc J.~Polchinski}, {\em String theory. {V}ol. {I}}, Cambridge Monographs on
  Mathematical Physics, Cambridge University Press, Cambridge, 1998.
\newblock An introduction to the bosonic string.

\bibitem{Polch2}
\leavevmode\vrule height 2pt depth -1.6pt width 23pt, {\em String theory.
  {V}ol. {II}}, Cambridge Monographs on Mathematical Physics, Cambridge
  University Press, Cambridge, 1998.
\newblock Superstring theory and beyond.

\bibitem{Rim}
{\sc D.~S. Rim}, {\em Formal deformation theory}, in Groupes de monodromie en
  g\'eom\'etrie alg\'ebrique. {I}, Springer-Verlag, Berlin, 1972, pp.~viii+523.
\newblock S\'eminaire de G\'eom\'etrie Alg\'ebrique du Bois-Marie 1967--1969
  (SGA 7 I), Dirig\'e par A. Grothendieck. Avec la collaboration de M. Raynaud
  et D. S. Rim, Lecture Notes in Mathematics, Vol. 288.

\bibitem{SW}
{\sc N.~Seiberg and E.~Witten}, {\em Comments on string dynamics in six
  dimensions}, Nuclear Phys. B, 471 (1996), pp.~121--134.

\bibitem{Sen}
{\sc A.~Sen}, {\em Stable non-{BPS} states in string theory}, J. High Energy
  Phys.,  (1998), pp.~Paper 7, 21 pp.\ (electronic).

\bibitem{FAC}
{\sc J.-P. Serre}, {\em Faisceaux alg\'ebriques coh\'erents}, Ann. of Math.
  (2), 61 (1955), pp.~197--278.

\bibitem{Ser65}
\leavevmode\vrule height 2pt depth -1.6pt width 23pt, {\em Alg\`ebre locale.
  {M}ultiplicit\'es}, vol.~11 of Cours au Coll\`ege de France, 1957--1958,
  r\'edig\'e par Pierre Gabriel. Seconde \'edition, 1965. Lecture Notes in
  Mathematics, Springer-Verlag, Berlin, 1965.

\bibitem{Se67}
{\sc C.~S. Seshadri}, {\em Space of unitary vector bundles on a compact
  {R}iemann surface}, Ann. of Math. (2), 85 (1967), pp.~303--336.

\bibitem{SiWi}
{\sc E.~Silverstein and E.~Witten}, {\em Criteria for conformal invariance of
  {$(0,2)$} models}, Nuclear Phys. B, 444 (1995), pp.~161--190.

\bibitem{Simp96a}
{\sc C.~T. Simpson}, {\em Moduli of representations of the fundamental group of
  a smooth projective variety. {I}}, Inst. Hautes \'Etudes Sci. Publ. Math.,
  (1994), pp.~47--129.

\bibitem{Ste}
{\sc P.~Stellari}, {\em Some remarks about the {F}{M}-partners of {K}3 surfaces
  with small {P}icard number}, Geom. Dedicata, 108 (2004), pp.~1--14.

\bibitem{Th00}
{\sc R.~P. Thomas}, {\em Derived categories for the working mathematician}.
\newblock {\tt math.AG/0001045}.

\bibitem{Tyau}
{\sc G.~Tian and S.-T. Yau}, {\em Three-dimensional algebraic manifolds with
  {$C\sb 1=0$} and {$\chi=-6$}}, in Mathematical aspects of string theory (San
  Diego, Calif., 1986), vol.~1 of Adv. Ser. Math. Phys., World Sci. Publishing,
  Singapore, 1987, pp.~543--559.

\bibitem{Tu94a}
{\sc L.~W. Tu}, {\em Semistable bundles over an elliptic curve}, Adv. Math., 98
  (1993), pp.~1--26.

\bibitem{Ue04}
{\sc H.~Uehara}, {\em An example of {F}ourier-{M}ukai partners of minimal
  elliptic surfaces}, Math. Res. Lett., 11 (2004), pp.~371--375.

\bibitem{Va99}
{\sc C.~Vafa}, {\em Extending mirror conjecture to {C}alabi-{Y}au with
  bundles}, Commun. Contemp. Math., 1 (1999), pp.~65--70.

\bibitem{Ver}
{\sc J.-L. Verdier}, {\em Cat\'egories d\'eriv\'ees. quelques r\'esultats (etat
  0)}, in Cohomologie \'etale, S\'eminaire de G\'eom\'etrie Alg\'ebrique du
  Bois-Marie SGA 4${1\over 2}$, Avec la collaboration de J. F. Boutot, A.
  Grothendieck, L. Illusie et J. L. Verdier, Springer-Verlag, Berlin, 1977,
  pp.~262--311. Lecture Notes in Mathematics, Vol. 569.

\bibitem{Ver96}
\leavevmode\vrule height 2pt depth -1.6pt width 23pt, {\em Des cat\'egories
  d\'eriv\'ees des cat\'egories ab\'eliennes}, Ast\'erisque,  (1996),
  pp.~xii+253 pp. (1997).
\newblock With a preface by Luc Illusie, Edited and with a note by Georges
  Maltsiniotis.

\bibitem{Wi86}
{\sc E.~Witten}, {\em New issues in manifolds of {${\rm SU}(3)$} holonomy},
  Nuclear Phys. B, 268 (1986), pp.~79--112.

\bibitem{WIT96}
\leavevmode\vrule height 2pt depth -1.6pt width 23pt, {\em Small instantons in
  string theory}, Nuclear Phys. B, 460 (1996), pp.~541--559.

\bibitem{Wit98}
\leavevmode\vrule height 2pt depth -1.6pt width 23pt, {\em D-branes and
  {$K$}-theory}, J. High Energy Phys.,  (1998), pp.~Paper 19, 41 pp.\
  (electronic).

\bibitem{Wit96a}
\leavevmode\vrule height 2pt depth -1.6pt width 23pt, {\em Dynamics of quantum
  field theory}, in Quantum fields and strings: a course for mathematicians,
  Vol. 1, 2 (Princeton, NJ, 1996/1997), Amer. Math. Soc., Providence, RI, 1999,
  pp.~1119--1424.

\bibitem{Yo01}
{\sc K.~Yoshioka}, {\em Moduli spaces of stable sheaves on abelian surfaces},
  Math. Ann., 321 (2001), pp.~817--884.

\end{thebibliography}

\end{document}